\documentclass[12pt]{article}
\usepackage{amsmath,amssymb,enumerate}
\usepackage{graphicx,xcolor} 
\usepackage[all]{xy}
\usepackage{bm}
\newtheorem{tm}{Theorem}[section]
\newtheorem{lm}[tm]{Lemma}
\newtheorem{co}[tm]{Corollary}
\newtheorem{re}[tm]{Remark}

\newtheorem{pr}[tm]{Proposition}

\newtheorem{ex}[tm]{Example}

\makeatletter
\newcommand{\subscripts}[3]{%
  \@mathmeasure\z@\displaystyle{#2}%
  \global\setbox\@ne\vbox to\ht\z@{}\dp\@ne\dp\z@
  \setbox\tw@\box\@ne
  \@mathmeasure4\displaystyle{\copy\tw@_{#1}}%
  \@mathmeasure6\displaystyle{{#2}_{#3}}%
  \dimen@-\wd6 \advance\dimen@\wd4 \advance\dimen@\wd\z@
  \hbox to\dimen@{}\mathop{\kern-\dimen@\box4\box6}%
}
\makeatother
\newcommand{\qed}{~~\hbox{\rule{4pt}{8pt}}}

\def\supp{\mathop{\rm supp}\nolimits}

\newcommand{\A}{\mathcal{A}}
\newcommand{\T}{\mathcal{T}}
\newcommand{\g}{\frak{g}}
\newcommand{\p}{\frak{p}}
\newcommand{\m}{\frak{m}}
\newcommand{\bg}{\bm{\mathrm{g}}}
\newcommand{\nn}{\nonumber}
\newcommand{\Prob}{\bm{\mathrm{P}}}

\newcommand{\h}{\mathrm{H}}
\newcommand{\e}{\mathrm{e}}
\newcommand{\x}{\bm{\mathrm{x}}}
\newcommand{\y}{\bm{\mathrm{y}}}
\newcommand{\ve}{\varepsilon}
\newcommand{\dis}{\displaystyle}

\newcommand{\Ker}{\mathrm{Ker}\,}
\newcommand{\im}{\mathrm{Im}\,}
\newcommand{\del}{\partial}
\newcommand{\ol}{\overline}
\newcommand{\la}{\langle}
\newcommand{\La}{\langle\!\langle}
\newcommand{\ra}{\rangle}
\newcommand{\Ra}{\rangle\!\rangle}
\newcommand{\LA}{\longrightarrow}
\newcommand{\Span}{\mathrm{span}_{\mathbb{R}}}
\newcommand{\G}{G_{\Gamma}}

\newcommand{\z}{\bm{\mathrm{z}}}
\newcommand{\Hom}{\mathrm{Hom}}
\newcommand{\Dom}{\mathrm{Dom}}

\newcommand{\Log}{\text{\rm{\,log\,}}}
\newcommand{\Exp}{\text{\rm{exp\,}}}
\newcommand{\hol}{\text{-H\"ol}}

\makeatletter
 
 \@addtoreset{equation}{section}
\makeatother
\begin{document}
\setlength{\baselineskip}
{15.5pt}
\title{
Central limit theorems for
non-symmetric 
random walks on 
nilpotent covering graphs: Part I 
}
\author{\Large
{Satoshi Ishiwata\footnote{
Department of Mathematical Sciences, Faculty of Science, Yamagata University,
1-4-12, Kojirakawa, Yamagata 990-8560, Japan
(e-mail: {\tt ishiwata@sci.kj.yamagata-u.ac.jp})}
\footnote{Partially supported by Grant-in-Aid for Young Scientists (B)(25800034) 
and Grant-in-Aid for Scientific Research (C)(17K05215) from JSPS.}, 
Hiroshi Kawabi\footnote{
Department of Mathematics, Faculty of Economics, Keio University, 4-1-1, 
Hiyoshi, Kohoku-ku, Yokohama 223-8521, Japan (e-mail: {\tt kawabi@keio.jp})}
\footnote{Partially supported by Grant-in-Aid for Scientific Research 
(C)(26400134, 17K05300) from JSPS.}
}
{\Large{and Ryuya Namba\footnote{
Graduate School of Natural Sciences, Okayama University, 3-1-1, 
Tsushima-Naka, Kita-ku, Okayama 700-8530, Japan (e-mail: {\tt{rnamba@s.okayama-u.ac.jp}})}
\hspace{0.1mm} \footnote{Partially supported by JSPS Research Fellowships for Young Scientists 
(18J10225).}
}}
}
\date{\today}
%
\maketitle 
%
%
\begin{abstract}

In the present paper, we study
central limit theorems (CLTs)
for non-symmetric random walks on nilpotent covering graphs
from a point of view of discrete geometric analysis 
developed by Kotani and Sunada.
We establish a semigroup CLT for a non-symmetric random walk
on a nilpotent covering graph.
Realizing the nilpotent covering graph into a nilpotent
Lie group through a discrete harmonic map, we give a geometric
characterization of the limit semigroup on the nilpotent Lie group.
More precisely, we show that the limit semigroup is generated by the 
sub-Laplacian with a non-trivial drift on the nilpotent Lie group
equipped with the Albanese metric. 
The drift term arises from the non-symmetry of the random walk
and it vanishes when the random walk is symmetric. 
Furthermore, by imposing the 
``centered condition'', we establish 
a functional CLT (i.e., Donsker-type invariance principle)
in a H\"older space over the nilpotent Lie group.
The functional CLT is extended to the case where the realization 
is not necessarily harmonic. 
We also obtain an explicit representation of the limiting diffusion 
process on the nilpotent Lie group and discuss a relation 
with rough path theory. 
Finally, we give several examples of random walks on nilpotent covering
graphs with explicit computations.

\vspace{3mm}

\vspace{3mm}
\noindent
{\bf AMS Classification (2010):} 
60F17, 
60G50, 
60J10, 
22E25 
\end{abstract}


\section{Introduction}
There are many interests in the study of random walks 
on infinite graphs
in many branches of mathematics such as probability theory, 
harmonic analysis, geometry, 
graph theory and group theory.
Among these branches, the long time behavior of random walks on infinite graphs is  
one of the major themes.
For instance, 
a central limit theorem (CLT), that is,
a generalization of the Laplace--de Moivre theorem,
has been studied 
intensively and extensively 
in various settings.
These mathematical backgrounds basically motivate our study. 
For basic results on random walks, we refer to Spitzer \cite{Spitzer}, 
Woess \cite{Woess}, Lawler--Limic \cite{Lawler} and references therein.  

\begin{figure}[htp]
\begin{center}
\includegraphics[width=12cm]{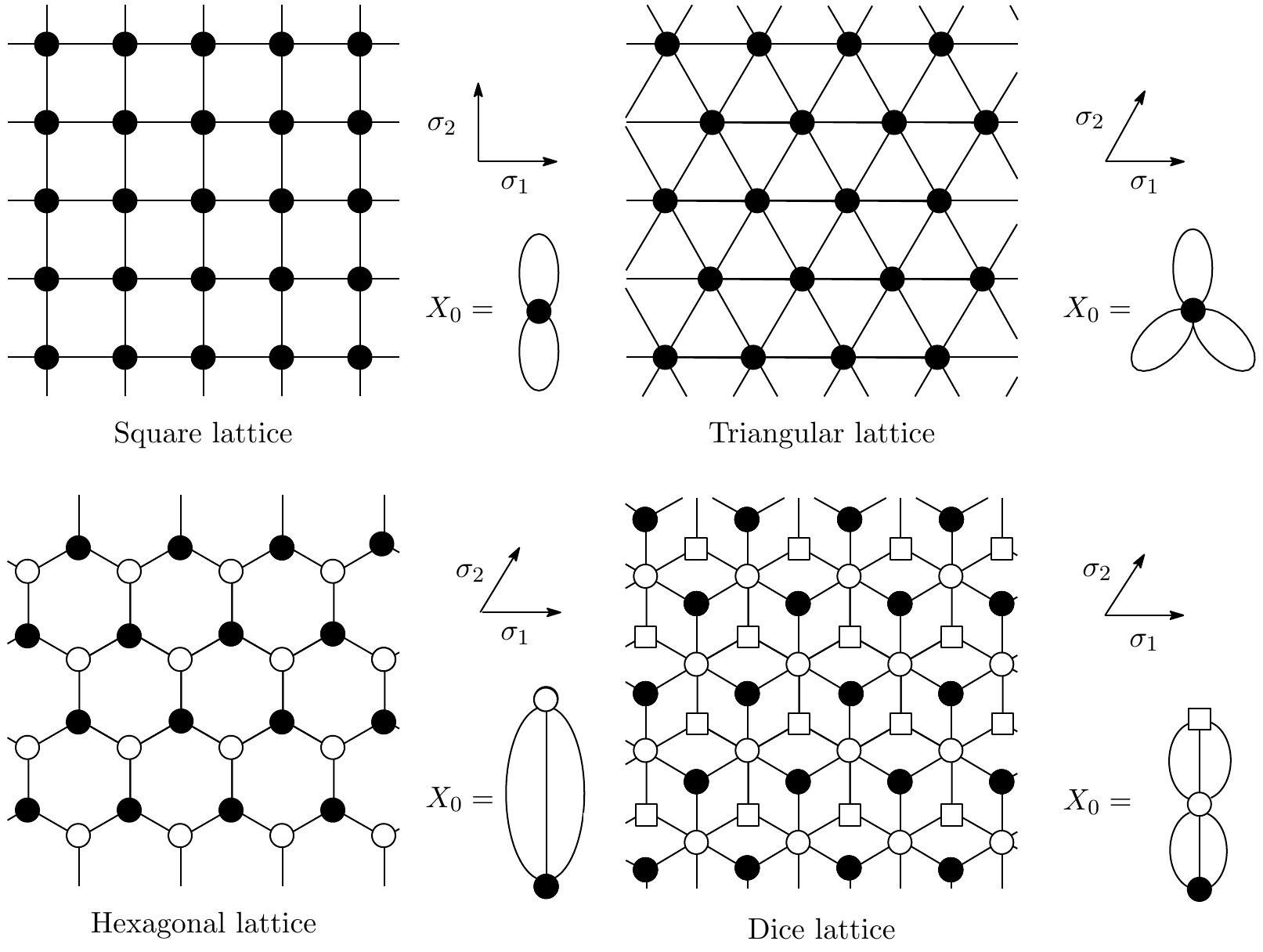}
\end{center}
\vspace{-0.5cm}
\caption{Crystal lattices with the covering transformation group 
$\Gamma=\langle \sigma_{1}, \sigma_{2} \rangle
\cong {\mathbb Z}^{2}$}
\label{crystal lattices}
\end{figure}

In these studies of random walks on infinite graphs, 
many authors have also discussed 
what kinds
of structures of underlying graphs affect the long time behavior of random walks.
It is known that geometric structures such as the {\it{periodicity}} 
of underlying graphs play important roles in them (cf.~Spitzer \cite{Spitzer}).
A typical example of periodic infinite graphs 
is a {\it{crystal lattice}}, that is, a covering graph $X$ of a finite graph $X_{0}$ 
whose covering transformation group $\Gamma$ is abelian. 
It is a generalization of the square lattice,
the triangular lattice, the hexagonal lattice, 
the dice lattice and so on (see Figure \ref{crystal lattices}). 
We remark that
the crystal lattice has 
inhomogeneous local structures though it has a periodic global structure.
Kotani and Sunada \cite{KS00-TAMS} introduced the notion of {\it{standard realization}}
of a crystal lattice $X$, which is a discrete
harmonic map from $X$ into the Euclidean space $\Gamma \otimes {\mathbb R}$
equipped with the {\it{Albanese metric}},
to characterize an equilibrium configuration of $X$. 
In a series of papers Kotani--Shirai--Sunada \cite{KSS}, Kotani \cite{Kotani} and 
Kotani--Sunada \cite{KS00-CMP, KS00-TAMS, KS06}, 
they developed 
a hybrid field 
of several traditional disciplines 
including graph theory, geometry, discrete group theory and 
probability theory.
Since 
this new field, called {\it discrete geometric analysis},
was introduced by Sunada,
it has been making new interactions with many other fields.
For example, Le Jan employs discrete geometric analysis
effectively in a series of recent studies of Markov loops (see e.g., \cite{Lejan1, Lejan2}). 
We refer to Sunada \cite{S, S2} for recent developments of discrete geometric analysis.
Especially, in \cite{KS00-CMP},  a geometric characterization of 
the diffusion semigroup appeared in the CLT-scaling limit 
of the symmetric random walk on $X$ was given
in terms of discrete geometric analysis. 
Ishiwata, Kawabi and Kotani \cite{IKK} generalized these results 
to the non-symmetric case and 
established two kinds of functional CLTs
(i.e., Donsker-type invariance principles)
for non-symmetric random walks on crystal lattices. 
We also refer to Guivar'ch \cite{Gui} and Kramli--Szasz \cite{KS83-PTRF} 
for related early works, Kotani \cite{Kotani contemp} and Kotani--Sunada \cite{KS06}
for a large deviation principle (LDP)
and Namba \cite{Namba} for yet another  
functional CLT for non-symmetric random walks on crystal lattices.

On the other hand, long time behaviors of symmetric or centered random walks on  
groups have been studied intensively and extensively. 
In particular, the notion of {\it{volume growth}} of groups
plays a key role
in the interface between probability theory and group theory.
Generally speaking, it is difficult to characterize
a finitely generated group itself in terms of its volume growth. 
We refer to Saloff-Coste \cite{SC-review} for basic problems and 
results for random walks on such groups including
ones of superpolynomial volume growth. 
On the contrary,
there is a remarkable theorem on a group of polynomial volume growth 
due to Gromov, which asserts that it is essentially characterized 
as a nilpotent group (cf.~Gromov \cite{Gromov} and Ozawa \cite{Ozawa}).
Hence, we find a large number of papers on long time behaviors of 
random walks on nilpotent groups. 
See e.g., Wehn \cite{Wehn}, Tutubalin \cite{Tutu}, Stroock--Varadhan \cite{SV}, 
Raugi \cite{Raugi}, Watkins \cite{Watkins}, Pap \cite{Pap} and Alexopoulos \cite{A3}
for related results on CLTs on nilpotent Lie groups, and
Breuillard \cite{Bre} 
for an overview of random walks on Lie groups.
We also refer to 
Alexopoulos \cite{A1, A2}, Breuillard \cite{Bre2}, 
Diaconis--Hough \cite{DH} and Hough \cite{Hough} for 
local CLTs on 
nilpotent Lie groups.


%

\begin{figure}[htp]
\begin{center}
\includegraphics[width=10cm]{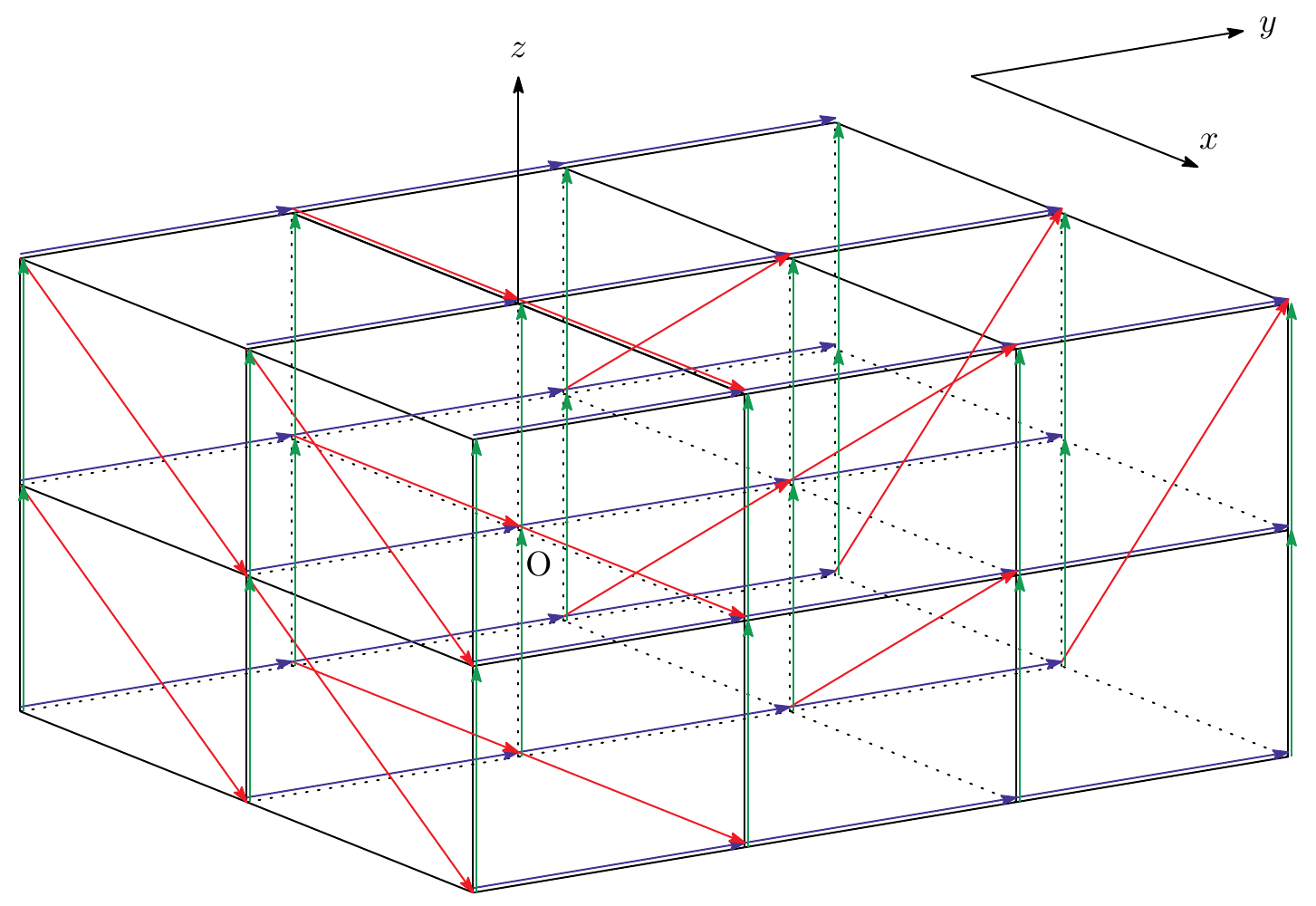}
\end{center}
\vspace{-0.5cm}
\caption{A part of the Cayley graph of $\Gamma=\mathbb{H}^3(\mathbb{Z})$ }
\label{Cayley graph}
\end{figure}

In view of these developments, we study the long time behavior of 
random walks on a covering graph $X$ whose
covering transformation group $\Gamma$ is a finitely generated group 
of polynomial volume growth. 
It is regarded as a generalization of a crystal lattice or
the Cayley graph of a finitely generated group of polynomial volume growth. 
A typical example of such a group is the 3-dimensional 
discrete Heisenberg group $\Gamma=\mathbb{H}^3(\mathbb{Z})$ 
(see Figure \ref{Cayley graph}). 
Thanks to Gromov's theorem mentioned above, $\Gamma$ has a finite extension of a torsion free
nilpotent subgroup $\widetilde{\Gamma} \lhd \Gamma$. 
Therefore, $X$ is regarded as a covering graph of the finite quotient graph 
$\widetilde{\Gamma} \backslash X$ 
whose covering transformation group is $\widetilde{\Gamma}$.
Throughout the present paper, we may assume that $X$
is a covering graph of a finite graph $X_0$ 
whose covering transformation group $\Gamma$ is a finitely generated, 
torsion free nilpotent group of step $r$ ($r \in \mathbb{N}$)
without loss of generality. 
We now mention a few related works. 
Ishiwata \cite{Ishiwata} discussed symmetric random walks on nilpotent covering graphs 
and extended the notion of standard realization of crystal lattices to the nilpotent case. 
Besides, in \cite{Ishiwata, Ishiwata2}, a semigroup CLT 
and a local CLT for symmetric random walks were obtained
by realizing the nilpotent covering graph $X$ into a nilpotent Lie group $G$
such that $\Gamma$ is isomorphic 
to a cocompact lattice in $G$ (cf.~Malc\'ev \cite{Malcev}).
We notice that, in spite of such developments, long time behaviors of 
non-symmetric random walks on
nilpotent covering graphs have not been studied sufficiently though
an LDP on nilpotent covering graphs was obtained in Tanaka \cite{Tanaka}. 


Under these circumstances, 
we establish CLTs
for non-symmetric random walks on a $\Gamma$-nilpotent covering graph $X$.
As an extension of the notion of standard realization introduced
in \cite{Ishiwata} to the non-symmetric case, 
we define the {\it modified standard realization} 
$\Phi_0$ from $X$ into a nilpotent Lie group 
$G=G_\Gamma$ whose Lie algebra is equipped with the Albanese metric.
Through the map $\Phi_0$, 
we obtain a semigroup CLT (Theorem \ref{CLT1}),
which means that the $n$-th iteration of the 
``transition shift operator'' converges to a diffusion semigroup 
on $G$ as $n \to \infty$ with a suitable scale change on $G$. 
The infinitesimal generator $-\A$ of the diffusion semigroup 
is the sub-Laplacian with a non-trivial drift $\beta(\Phi_0)$
affected by the non-symmetry of the given random walk. 
Furthermore, by imposing an additional natural condition  {\bf (C)}, 
we prove a functional CLT in a H\"older space over $G$ (Theorem \ref{FCLT1}), 
which is much stronger than Theorem \ref{CLT1}. 
Roughly speaking, we capture a $G$-valued diffusion process associated with $-\A$
through the CLT-scaling limit of the non-symmetric random walk on $X$.
We call the condition {\bf (C)}  the {\it centered condition}.
The functional CLT is also extended to the case where
the realization $\Phi : X \LA G$ is not necessarily harmonic (Theorem \ref{FCLT1-general})
under {\bf (C)}. 
In this case, several technical difficulties appear in the proof of the functional CLT. 
To overcome them,we take a modified harmonic realization $\Phi^{(\varepsilon)}_{0}: X \to G$
and show that the
(${\mathfrak g}^{(1)}$-){\it{corrector}},
the difference between 
$\Phi^{(\varepsilon)}$ and $\Phi^{(\varepsilon)}_{0}$ in 
the ${\mathfrak g}^{(1)}$-direction, is not so big.
This approach is the so-called
{\it{corrector method}} in the context of stochastic homogenization theory,
and it is effectively used in the study of random walks in random environments
(see e.g., Papanicolaou--Varadhan \cite{PV}, Kozlov \cite{Kozlov} 
and Kumagai \cite{Kumagai}).We then obtain that 
a sequence $\big\{ \tau_{n^{-1/2}} \big( \Phi(w_{[nt]}) \big);
0\leq t \leq 1 \big\}_{n=1}^{\infty}$ also converges 
in law 
to the diffusion process $(Y_{t})_{0\leq t \leq 1}$ as $n \to \infty$.
In a subsequent paper \cite{IKN}, we will consider
the {\it weakly asymmetric case} and establish
another kind of CLTs for a family of random walks 
on the nilpotent covering graph $X$
which interpolates between the original non-symmetric random walk
and the symmetrized one. We also capture a $G$-valued diffusion process 
different from the one obtained in the present paper.  
The comparison between these two diffusions will be given in Remark \ref{diffusions}. 


Let us give another motivation of the present paper from rough path theory. 
It is known that  rough path theory was first initiated by Lyons in \cite{Lyons}
to discuss line integrals and ordinary differential equations (ODEs) driven by an irregular path
such as a sample path of Brownian motion $B=(B_t)_{0 \leq t \leq 1}$ on $\mathbb{R}^d$. 
Rough path theory 
makes us possible to handle 
a Stratonovich type stochastic differential equation (SDE) 
driven by Brownian motion $B$ as
a deterministic ODE driven by standard 
{\it Brownian rough path} (i.e., Stratonovich enhanced Brownian motion) 
${\bf B}=(B, \mathbb{B})$,
where ${\bf B}$ is a couple of Brownian motion $B$ itself and its Stratonovich 
iterated integral $\mathbb{B}$.
Thus rough path theory provides a new insight to the usual SDE-theory and
it has developed rapidly in stochastic analysis.
For more details on an overview of rough path theory
and its applications to stochastic analysis, see
Lyons--Qian \cite{LQ}, Friz--Victoir \cite{FV} and Friz--Hairer \cite{FH}. 
In the rough path framework,
several authors have studied Donsker-type invariance principles. 
Among them, 
Breuillard--Friz--Huesmann \cite{BFH} first studied this problem
for Brownian rough path. 
Namely, they captured Stratonovich enhanced
Brownian motion ${\bf B}=(B, \mathbb{B})$ on $\mathbb{R}^d$ as the usual CLT-scaling limit 
of the natural rough path lift of an $\mathbb{R}^d$-valued random walk
with the centered condition. 
We also refer to Bayer--Friz \cite{BF} for applications to cubature and 
Chevyrev \cite{Chevyrev} for a recent study on an extension
to the case of L\'evy processes. 
Here we should note that there are good approximations to Brownian motion
which do not converge to ${\bf B}$ but instead to a {\it distorted Brownian rough path}
$\ol{{\bf B}}=(B, \mathbb{B}+{\bm \beta})$,
where ${\bm \beta}$ is an anti-symmetric perturbation of $\mathbb{B}$. 
For example, 
Friz--Gassiat--Lyons \cite{FGL} constructed
such a rough path called {\it magnetic Brownian rough path}
as the small mass limit of the natural rough path lift of 
a physical Brownian motion on $\mathbb{R}^d$ in a magnetic field.
Through this approximation, 
they showed an effect of the magnetic field 
appears explicitly in the anti-symmetric perturbation term ${\bm \beta}$. 
See also e.g., Lejay--Lyons \cite{LL} and Friz--Oberhauser \cite{FO}
for related results on this topic. 

In view of the background described above, 
we discuss a random walk approximation of the 
distorted Brownian rough path $\ol{{\bf B}}$ 
from a perspective of discrete geometric analysis.  
Since the unique {\it Lyons extension} of $\ol{{\bf B}}$ 
of order $r$ ($r \geq 2$) can be regarded as a diffusion process on 
a free step-$r$ nilpotent Lie group $\mathbb{G}^{(r)}(\mathbb{R}^d)$
(see Section 5 below for definition), 
we obtain such a diffusion process in Corollary \ref{RP-relation} 
through the CLT-scaling limit 
of a non-symmetric random walk on a nilpotent covering graph $X$
as a direct application of Theorem \ref{FCLT1}. 
Besides, we observe that the non-symmetry of 
the random walk on $X$ affects 
the anti-symmetric perturbation term of $\ol{{\mathbb{B}}}$ explicitly. 
Recently, Lopusanschi--Simon \cite{LS} proved a 
similar invariance principle for $\ol{{\bf B}}$ to ours.
However, they did not discuss an explicit relation between the perturbation term, 
called the {\it area anomaly}, 
and the non-symmetry of the given random walk. 
See also Lopusanschi--Orenshtein \cite{LO} for a related result. 
In view of that, Corollary \ref{RP-relation} gives a 
new approach to such an invariance principle in that 
we pay much attention to the non-symmetry of random walks on $X$.


The rest of the present paper is organized as follows:
We introduce our framework and state the main results in Section 2.
We make a preparation from nilpotent Lie groups, 
the Carnot--Carath\'eodory metric, homogeneous norms and 
discrete geometric analysis in Section 3. 
A relation between the $G$-valued Markov chain and the notion of 
modified harmonicity is also discussed. 
In the former part of Section 4,
we prove the first main result (Theorem \ref{CLT1}) and 
give several properties of the non-trivial drift $\beta(\Phi_0)$ 
(Proposition \ref{property-beta}).
Trotter's approximation theorem plays a crucial role in the proof of Theorem \ref{CLT1}. 
We then prove a functional CLT (Theorem \ref{FCLT1}) for
the non-symmetric random walk under the centered condition {\bf (C)} 
in the latter part of Section 4.
We show the tightness of the family of probability measures 
induced by the $G$-valued stochastic processes 
given by the geodesic interpolation of the given random walk (Lemma \ref{tightness1}).
In the case $r=2$, we prove it
by combining the modified harmonicity of $\Phi_0$ with standard martingale techniques.
On the other hand, the same argument is insufficient 
in the case $r \geq 3$. 
To handle the higher-step terms, 
we employ a novel pathwise argument inspired by
the proof of the Lyons extension theorem (cf.~Lyons \cite{Lyons}) in rough path theory.
However, we need a careful examination of the proof of Lyons' extension theorem
since rough path theory is build on free nilpotent Lie groups and 
our nilpotent Lie group $G$ is not necessarily free. 
As a consequence of Theorem \ref{CLT1}, 
the convergence of the finite dimensional distribution
of the stochastic process (Lemma \ref{CFDD1}) is proved. 
Moreover, a functional CLT in the case where the realization is non-harmonic 
(Theorem \ref{FCLT1-general}) is also proved by applying the corrector method described above. 
An explicit representation of the limiting diffusion process
is given in Section 5. 
We also discuss a relation between this diffusion process
and rough path theory by using this representation formula in the case where
$\Gamma$ is the free discrete nilpotent group over $\mathbb{Z}^d$. 
We give two examples of 
non-symmetric random walks on nilpotent covering graphs
with explicit calculations in Section 6. 
Finally, we give a comment on another approach to CLTs in the non-centered case
in Appendix A 
(see Theorems \ref{Another CLT} and \ref{Another FCLT}).

Throughout the present paper, $C$ denotes a positive constant that may change from line to line 
and $O(\cdot)$ stands for the Landau symbol. 
If the dependence
of $C$ and $O(\cdot)$ 
are significant, 
we denote them like $C(N)$ and $O_{N}(\cdot)$,
respectively.


\section{Framework and Results}

We introduce our framework and state the main results in this section.  
Let $\Gamma$ be a torsion free, finitely generated nilpotent group and 
$X=(V, E)$ a $\Gamma$-nilpotent covering graph, 
where $V$ is the set of all vertices and $E$ is the set of all oriented edges. 
The graph $X$ possibly have multiple edges or loops and is equipped 
with the discrete topology induced by the graph distance. 
For an edge $e \in E$, 
we denote by $o(e)$ and $t(e)$ the origin and the terminus of $e$, respectively. 
The inverse edge of $e \in E$ is defined by an edge, say $\ol{e}$, satisfying $o(\ol{e})=t(e)$
and $t(\ol{e})=o(e)$. 
We set $E_x=\{e \in E \, | \, o(e)=x\}$ for $x \in V$. 
A path $c$ in $X$ of length $n$ is a sequence $c=(e_1, e_2, \dots, e_n)$ of $n$ edges 
$e_1, e_2, \dots, e_n \in E$ with $o(e_{i+1})=t(e_i)$ for $i=1, 2, \dots, n-1$. 
We denote by $\Omega_{x, n}(X)$ 
the set of all paths in $X$ of length $n \in \mathbb{N} \cup \{\infty\}$ starting from $x \in V$. 
Put $\Omega_x(X)=\Omega_{x, \infty}(X)$ for simplicity. 

We introduce a {\it transition probability}, that is, 
a function $p : E \LA [0, 1]$ satisfying 
$$
\sum_{e \in E_x}p(e)=1 \qquad (x \in V) \quad \text{and} \quad p(e)+p(\ol{e})>0 \qquad (e \in E).
$$
Moreover, we impose that $p$ is invariant under the $\Gamma$-action, that is, 
$p(\gamma e)=p(e)$ for $\gamma \in \Gamma$ and $e \in E$.  
The value $p(e)$ represents the probability that a particle at the origin $o(e)$ 
moves to the terminus $t(e)$ along the edge $e \in E$ in a unit time. 
The random walk associated with $p$ is 
the $X$-valued time-homogeneous Markov chain 
$(\Omega_x(X), \mathbb{P}_x, \{w_n\}_{n=0}^\infty)$, where $\mathbb{P}_x$ is 
the probability measure on $\Omega_x(X)$ satisfying
$$
\mathbb{P}_x\big(\{c=(e_1, e_2, \dots, e_n, *, *, \dots)\}\big)
=p(e_1)p(e_2) \cdots p(e_n)
\qquad \big(c \in \Omega_x(X)\big)
$$
and 
$w_n(c):=o( e_{n+1})$ for $n \in \mathbb{N} \cup\{0\}$
and $c=(e_1, e_2, \dots, e_n, \dots) \in \Omega_x(X)$. 

We define the {\it transition operator} $L$ associated with the transition probability $p$ by
$$
Lf(x):=\sum_{e \in E_x}p(e)f\big(t(e)\big) \qquad (x \in V, \, f : V \LA \mathbb{R})
$$
and the $n$-step transition probability $p(n, x, y)$ by
$p(n, x, y):=L^n \delta_y(x)$ for $n \in \mathbb{N}$ and $x, y \in V$, 
where $\delta_y$ stands for the Dirac delta function at $y$ 
and $p(c)=p(e_1)p(e_2) \cdots p(e_n)$ for $c=(e_1, e_2, \dots, e_n) \in \Omega_{x, n}(X)$.  
Let $X_0=(V_0, E_0)=\Gamma \backslash X$ be the quotient graph, which is finite by definition.  
Then the random walk on $X_0$ is  
 induced through the covering map $\pi : X \LA X_0$. 
 We write $p : E_0 \LA [0, 1]$ the transition probability on $X_0$, by abuse of notation. 
For $n \in \mathbb{N}$ and $x, y \in V_0$, 
we also denote by $p(n, x, y)$ 
the $n$-step transition probability of the random walk on $X_0$.
In what follows, we assume that
the random walk on $X_0$ is {\it irreducible}. 
Namely, for $x, y \in V_0$, there exists $n=n(x, y) \in \mathbb{N}$ 
such that $p(n, x, y)>0$.
We then find 
a unique positive function $m : V_0 \LA (0, 1]$ which is called the 
{\it invariant measure} on $X_0$ satisfying
$$
\sum_{x \in V_0}m(x)=1 \quad\text{and}\quad
 m(x)=\sum_{e \in (E_0)_x}p(\ol{e})m\big(t(e)\big) \qquad (x \in V_0),
$$
thanks to the Perron-Frobenius theorem.
We set $\widetilde{m}(e):=p(e)m\big(o(e)\big)$ for $e \in E_0$. 
The random walk on $X_0$ is called ($m$-){\it symmetric} 
if $\widetilde{m}(e)=\widetilde{m}(\ol{e})$ for $e \in E_0$.
Otherwise, it is called ($m$-){\it non-symmetric}. 
We also write $m : V \LA (0, 1]$ for the $\Gamma$-invariant lift 
of $m : V_0 \LA (0, 1]$. Namely, $m$ takes the same value 
for every fiber $\pi^{-1}(x_0) \, (x_0 \in V_0)$ and satisfies
$m(\gamma x)=m(x)$ for $\gamma \in \Gamma$ and $x \in V$.
We denote by $\h_1(X_0, \mathbb{R})$ and $\h^1(X_0, \mathbb{R})$ 
 the first homology group 
and the first cohomology group of $X_0$, respectively. 
We define the {\it homological direction} of the given random walk on $X_0$ by
$$
\gamma_p:=\sum_{e \in E_0}\widetilde{m}(e)e \in \h_1(X_0, \mathbb{R}).
$$
It is clear that the random walk on $X_0$ is $(m$-)symmetric
if and only if $\gamma_p=0$. 
In this sense, $\gamma_p$ gives the homological drift of given random walk on $X_0$. 

On the other hand, we provide a continuous state space 
in which the $\Gamma$-nilpotent 
covering graph $X$ is properly realized. 
There exists a connected and simply connected nilpotent Lie group $(G, \cdot )=\G$ 
such that $\Gamma$ is isomorphic to a cocompact lattice in $G$
by applying Malc\'ev's theorem (cf.~Malc\'ev \cite{Malcev}).
A piecewise smooth $\Gamma$-equivariant map 
$\Phi : X \LA G$ is called a {\it periodic realization} of $X$.
Let $(\g, [\cdot, \cdot])$ be the Lie algebra of $G$, which is regarded as the tangent space 
at the unit $\bm{1}_G$. 
Since the exponential map $\Exp : \g \LA G$ is a diffeomorphism,
global coordinate systems on $G$ are induced through the exponential map.  
We write $\Log : G \LA \g$ for the inverse map of $\Exp : \g \LA G$. 

We construct a new product $*$ on $G$ in the following manner. 
Set
$\frak{n}_1:=\frak{g}$ and $\frak{n}_{k+1}:=[\frak{g}, \frak{n}_k]$ for $k \in \mathbb{N}$.
Since $\g$ is nilpotent, 
we find an integer $r \in \mathbb{N}$ such that
$\frak{g}=\frak{n}_1 \supset \dots 
\supset \frak{n}_r \supsetneq \frak{n}_{r+1}=\{\bm{0}_{\frak{g}}\}$.
The integer $r$ is called the {\it step number} of $\g$ or $G$. 
We define the subspace $\frak{g}^{(k)}$ of $\frak{g}$ by
$\frak{n}_k=\frak{g}^{(k)} \oplus \frak{n}_{k+1}$ for $k=1, 2, \dots, r$.
Then the Lie algebra $\frak{g}$ is decomposed as 
$\frak{g}=\frak{g}^{(1)} \oplus \g^{(2)} \oplus \dots \oplus \frak{g}^{(r)}$ 
and each  $Z \in \frak{g}$ is uniquely written as $Z=Z^{(1)} + Z^{(2)}+\dots + Z^{(r)}$, 
where 
$Z^{(k)} \in \frak{g}^{(k)}$ for $k=1, 2, \dots, r$.
Define a map $\tau^{(\g)}_\ve : \g \LA \g$ by
$$
\tau^{(\g)}_\ve(Z):=\ve Z^{(1)} + \ve^2 Z^{(2)} + \dots + \ve^r Z^{(r)} \qquad (\ve \geq 0, \, Z \in \frak{g})
$$
and also define a Lie bracket product $[\![ \cdot, \cdot ]\!]$ on $\frak{g}$ by
$$
[\![ Z_1, Z_2 ] \!]:=\lim_{\ve \searrow 0} \tau^{(\g)}_\ve
 \big[\tau^{(\g)}_{1/\ve}(Z_1), \tau^{(\g)}_{1/\ve}(Z_2)\big] 
\qquad (Z_1, Z_2 \in \frak{g}).
$$ 
We introduce a map $\tau_\ve : G \LA G$, called the {\it dilation operator} on $G$, by
$$
\tau_\ve (g):=\Exp \big( \tau^{(\g)}_\ve \big( \Log(g)\big)\big) \qquad (\ve \geq 0, \, g \in G),
$$
which gives scalar multiplications on $G$. 
We note that $\tau_\ve$ may not be a group homomorphism, though 
it is a diffeomorphism on $G$. 
By making use of the dilation map $\tau_\ve$, a Lie group product $*$ on $G$ 
is defined as follows:
$$
g * h :=\lim_{\ve \searrow 0} 
\tau_{\ve}\big( \tau_{1/\ve}(g) \cdot \tau_{1/\ve}(h)\big) \qquad (g, h \in G).
$$
The Lie group $G_\infty=(G, *)$ is called the {\it limit group} of $(G, \cdot)$. 
It is a {\it stratified Lie group} of step $r$
in the sense that $(\g, [\![\cdot, \cdot]\!])$ is decomposed as
$\g=\bigoplus_{k=1}^r \g^{(k)}$
satisfying
$
[\![\frak{g}^{(k)}, \frak{g}^{(\ell)}]\!]
\subset \frak{g}^{(k+\ell)}$ unless $k+j > r$
and the subspace $\frak{g}^{(1)}$ generates $\frak{g}$. 
The relation between these two Lie group products is given in the next section. 
We endow $G$ with the so-called {\it Carnot--Carath\'eodory metric} $d_{\mathrm{CC}}$, 
which is an intrinsic metric defined by
\begin{equation}\label{def-dcc}
d_{\mathrm{CC}}(g, h):=\inf\Big\{ \int_0^1 \|\dot{w}_t\|_{g_0} \, dt \, 
\Big| \, \begin{matrix} w \in \mathrm{Lip}([0, 1]; G), \, w_0=g, \, w_1=h, \\ 
\text{$w$ is tangent to $\g^{(1)}$}\end{matrix} \Big\}
\end{equation}
for $g, h \in G$, where we write $\mathrm{Lip}([0, 1]; G)$ 
for the set of all Lipschitz continuous 
paths and $\|\cdot\|_{g_0}$ for the norm on $\g^{(1)}$ induced by the Albanese metric (see Section 2.3 for details).

Let $\pi_1(X_0)$ be the fundamental group of $X_0$. 
Then we have a canonical surjective homomorphism 
$\rho : \pi_1(X_0) \LA \Gamma$
by the general theory of covering spaces. 
This map gives rise to a surjective homomorphism 
$\rho : \h_1(X_0, \mathbb{Z}) \LA \Gamma/[\Gamma, \Gamma]$ 
and we have
a surjective linear map 
$\rho_{\mathbb{R}} : \h_1(X_0, \mathbb{R}) \LA \frak{g}^{(1)}$
by extending it linearly. 
We call $\rho_{\mathbb{R}}(\gamma_p) \in \frak{g}^{(1)}$
the {\it asymptotic direction}. 
Note that $\gamma_p=0$ implies 
$\rho_{\mathbb{R}}(\gamma_p)=\bm{0}_{\frak{g}}$.
However, the converse does not always hold. 
We induce a special flat metric $g_0$ 
on $\frak{g}^{(1)}$, which is
called the {\it Albanese metric} associated with the transition probability $p$
by using the discrete Hodge-Kodaira theorem (cf.~Kotani--Sunada \cite[Lemma 5.2]{KS06}).
The construction of the metric is given in the next section. 
A periodic realization $\Phi_0 : X \LA G$ is called  ($p$-){\it modified harmonic} if
\begin{equation}\label{m-harmonicity}
\sum_{e \in E_x}p(e)\Log\Big(\Phi_0\big(o(e)\big)^{-1} \cdot  \Phi_0\big(t(e)\big)\Big)\Big|_{\g^{(1)}}
=\rho_{\mathbb{R}}(\gamma_p) \qquad (x \in V).
\end{equation}
Such $\Phi_0$ is uniquely determined up to $\g^{(1)}$-translation. 
The modified harmonicity describes the most natural realization of
the nilpotent covering graph $X$ in the geometric point of view. 
If we equip $\g^{(1)}$ with the Albanese metric $g_0$, 
the modified harmonic realization 
$\Phi_0 : X \LA G$ is called the {\it modified standard realization}.

For a metric space $\T$, we denote by $C_\infty(\T)$ 
the Banach space of continuous functions $f : \T \LA \mathbb{R}$ vanishing at infinity
with the uniform topology $\|\cdot \|_\infty^{\T}$. 
For $q >1$, we define
$$
C_{\infty, q}(X \times \mathbb{Z})
:=\big\{ f=f(x, z) : X \times \mathbb{Z} \LA \mathbb{R} \, \big| \, 
f(\cdot, z) \in C_\infty(X), \, \|f\|_{\infty, q}<\infty\big\},
$$
where $\|f\|_{\infty, q}$ is a norm on $C_{\infty, q}(X \times \mathbb{Z})$ given by
$$
\|f\|_{\infty, q}:=\frac{1}{C_q}\sum_{z \in \mathbb{Z}}\frac{\|f(\cdot, z)\|_{\infty}^X}{1+|z|^q},
\qquad C_q:=\sum_{z \in \mathbb{Z}}\frac{1}{1+|z|^q}<\infty.
$$
Then we see that $(C_{\infty, q}(X \times \mathbb{Z}), \|\cdot\|_{\infty, q})$
is a Banach space. 
We introduce the 
{\it transition-shift operator} 
$\mathcal{L}_p : C_{\infty, q}(X \times \mathbb{Z}) \LA C_{\infty, q}(X \times \mathbb{Z})$ by
\begin{equation}\label{transition-shift}
\mathcal{L}_p f(x, z):=\sum_{e \in E_x}p(e)f\big(t(e), z+1\big) \qquad (x \in V, \, z \in \mathbb{Z})
\end{equation}
and the {\it approximation operator} 
$\mathcal{P}_\ve : C_\infty(G) \LA C_{\infty, q}(X \times \mathbb{Z})$ by
\begin{equation}\label{scaling}
\mathcal{P}_\ve f(x, z):=f\Big( \tau_\ve \big(\Phi_0(x)*\exp(-z\rho_{\mathbb{R}}(\gamma_p))\big)\Big) 
\qquad (0 \leq \ve \leq 1, \, x \in V, \, z \in \mathbb{Z}).
\end{equation}
We extend each $Z \in \frak{g}$ as a left invariant vector field
$Z_*$ on $G$ as follows:
$$
Z_*f(g)=\frac{d}{d\ve}\Big|_{\ve=0}f\big( g * \Exp(\ve Z)\big) 
\qquad \big(f \in C^\infty(G), \, g \in G\big).
$$ 
We put 
$$
\beta(\Phi_0)
:=\sum_{e \in E_0}\widetilde{m}(e)
\Log \Big( \Phi_0\big(o(\widetilde{e})\big)^{-1} 
\cdot \Phi_0\big(t(\widetilde{e})\big)
\cdot \exp(-\rho_{\mathbb{R}}(\gamma_p))\Big)\Big|_{\frak{g}^{(2)}},
$$
where $\widetilde{e}$ stands for a lift of $e \in E_0$ to $X$. 
We note that $\gamma_p=0$ implies $\beta(\Phi_0)=\bm{0}_{\g}$. 
However, even if $\rho_{\mathbb{R}}(\gamma_p) =\bm{0}_{\g}$,
the quantity $\beta(\Phi_0)$ does not vanish in general. 
Furthermore, 
$\beta(\Phi_0)$ does not depend on $\g^{(2)}$-components of 
the modified harmonic realization $\Phi_0 : X \LA G$, though it has the ambiguity 
in the components corresponding to $\g^{(2)} \oplus \g^{(3)} \oplus \cdots \oplus \g^{(r)}$.
See Proposition \ref{property-beta} for details and Section 6.2 
for a concrete example.

Then the first main result is as follows:

\begin{tm}\label{CLT1}
For $q>4r+1$, the following hold:

\vspace{1mm}
\noindent
{\rm(1)} 
For $0 \leq s \leq t$ and $f \in C_\infty(G)$, we have
\begin{equation}\label{semigroup CLT1}
\lim_{n \to \infty}\Big\| \mathcal{L}_p^{[nt]-[ns]}\mathcal{P}_{n^{-1/2}} f 
- \mathcal{P}_{n^{-1/2}} \e^{-(t-s)\A}f\Big\|_{\infty, q}=0,
\end{equation}
where $(\e^{-t\A})_{t \geq 0}$ is the $C_0$-semigroup with the infinitesimal generator 
$\A$ on $C_0^\infty(G)$ defined by
\begin{equation}\label{generator}
\A:=-\frac{1}{2}\sum_{i=1}^{d_1}V_{i*}^2 - \beta(\Phi_0)_*,
\end{equation}
where $\{V_1, V_2, \dots, V_{d_1}\}$ denotes an orthonormal basis of $(\g^{(1)}, g_0)$. 

\vspace{2mm}
\noindent
{\rm (2)} Let $\mu$ be a Haar measure on $G$. Fix $z \in \mathbb{Z}$. 
Then, for any sequence $\{x_n\}_{n=1}^\infty \subset V$ with 
$$
\lim_{n \to \infty}\tau_{n^{-1/2}}\Big( \Phi_0(x_n) * \exp\big(-z\rho_{\mathbb{R}}(\gamma_p)\big)\Big)=g \in G
$$
and for any $f \in C_{\infty}(G)$, we have
\begin{equation}\label{semigroup CLT1-2}
\lim_{n \to \infty}\mathcal{L}_p^{[nt]}\mathcal{P}_{n^{-1/2}}f(x_n, z)
=\e^{-t\A}f(g):=\int_G \mathcal{H}_t(h^{-1} *g)f(h) \, \mu(dh) \qquad (t>0),
\end{equation}
where $\mathcal{H}_t(g)$ is a fundamental solution to $\del u/\del t+\A u=0.$
\end{tm}

Fix a reference point $x_* \in V$ with $\Phi_0(x_*)=\bm{1}_G$ and put 
$\xi_n(c):=\Phi_0\big(w_n(c)\big)$
for $n \in \mathbb{N}\cup\{0\}$ and $c \in \Omega_{x_*}(X)$.
We then have a $G$-valued random walk $(\Omega_{x_*}(X), \mathbb{P}_{x_*}, \{\xi_n\}_{n=0}^\infty)$
starting from $\bm{1}_G$.
For $t \geq 0$, we define a map $\mathcal{X}_t^{(n)} : \Omega_{x_*}(X) \LA G$ by
$$
\mathcal{X}_t^{(n)}(c):=\tau_{n^{-1/2}}\Big( \xi_{[nt]}(c)*\exp\big(-[nt]\rho_{\mathbb{R}}(\gamma_p)\big)\Big) 
\qquad \big(n \in \mathbb{N}, \, c \in \Omega_{x_*}(X)\big).
$$
Denote by
$\mathcal{D}_n$ the partition $\{t_k=k/n \, | \, k=0, 1,  \dots, n\}$
of $[0, 1]$ for $n \in \mathbb{N}$. 
We define a $G$-valued continuous stochastic process $(\mathcal{Y}_t^{(n)})_{0 \leq t \leq 1} $
by the geodesic interpolation of $\{\mathcal{X}_{t_k}^{(n)}\}_{k=0}^n$with respect to 
$d_{\mathrm{CC}}$. 
It is worth noting that (\ref{semigroup CLT1-2}) implies 
\begin{equation}\label{semigroup CLT1-3}
\lim_{n \to \infty}\sum_{c \in \Omega_{x_*}(X)}p(c)f\big(\mathcal{X}_t^{(n)}(c)\big)
=\int_G \mathcal{H}_t(h^{-1})f(h) \, \mu(dh) \qquad \big( f \in C_\infty(G)\big).
\end{equation}

We now consider an SDE
\begin{equation}\label{SDE}
dY_t = \sum_{i=1}^{d_1}V_{i*}(Y_t) \circ dB_t^i + \beta(\Phi_0)_*(Y_t) \, dt,\qquad Y_0=\bm{1}_G,
\end{equation}
where $(B_t)_{0 \leq t \leq 1}=(B_t^1, B_t^2, \dots, B_t^{d_1})_{0 \leq t \leq 1}$ 
is an $\mathbb{R}^{d_1}$-valued
standard Brownian motion with $B_0=\bm{0}$. 
Let $(Y_t)_{0 \leq t \leq 1}$ be the $G$-valued diffusion process which solves (\ref{SDE}). 
In Proposition \ref{SDE-solution} below,
we prove that the infinitesimal generator of $(Y_t)_{0 \leq t \leq 1}$ coincides with 
$-\A$ defined by (\ref{generator}). 
Let $C_{\bm{1}_G}([0, 1]; G)$ be the set of all continuous paths $w : [0, 1] \LA G$
such that $w_0=\bm{1}_G$ and $\mathrm{Lip}([0, 1]; G) \subset C_{\bm{1}_G}([0, 1]; G)$ 
the set of all Lipschitz continuous paths. 
For $\alpha<1/2$, we define the $\alpha$-H\"older distance $\rho_\alpha$ 
on $C_{\bm{1}_G}([0, 1]; G)$ by
$$
 \rho_\alpha(w^1, w^2):=\sup_{0 \leq s < t \leq 1}
\frac{d_{\mathrm{CC}}(u_s, u_t)}{|t-s|^\alpha},
\qquad u_t:=(w_t^1)^{-1}*w_t^2 \qquad (0 \leq t \leq 1). 
$$
We set
$C_{\bm{1}_G}^{0, \alpha\hol}([0, 1]; G)
:=\ol{\mathrm{Lip}([0, 1]; G)}^{\rho_\alpha}$,
which is a Polish space
(cf.~Friz--Victoir \cite[Section 8]{FV}).   
Let $\Prob^{(n)}$ be the image measure on $C^{0, \alpha\hol}_{\bm{1}_G}([0, 1]; G)$
induced by $\mathcal{Y}_\cdot^{(n)}$ for $n \in \mathbb{N}$.

We now in a position to present a functional CLT, the second main theorem,
for the non-symmetric random walk $\{w_n\}_{n=0}^\infty$ on $X$. 

\begin{tm}\label{FCLT1} 
We assume the centered condition
{\bf (C):} 
$\rho_{\mathbb{R}}(\gamma_p)=\bm{0}_{\g}$.
Then the sequence $(\mathcal{Y}_t^{(n)})_{0 \leq t \leq 1}$
converges in law to the $G$-valued diffusion process $(Y_t)_{0 \leq t \leq 1} $
in $C^{0, \alpha\text{\normalfont{-H\"ol}}}_{\bm{1}_G}([0, 1]; G)$ as $n \to \infty$ 
for all $\alpha<1/2$.
\end{tm}

Finally, we generalize Theorem \ref{FCLT1} to the case
where the realization is not necessarily modified harmonic. 
We take a  periodic realizations $\Phi : X \LA G$
such that $\Phi(x_*)=\bm{1}_G$ for some base point $x_* \in V$.
On the other hand, we may take the 
modified harmonic realization $\Phi_0 : X \LA G$   
such that $\Phi_0(x)^{(i)}=\Phi(x)^{(i)}$ for $x \in V$
and $i=2, 3, \dots, r$ without loss of generality. 
We now define the ($\g^{(1)}$-){\it corrector} $\mathrm{Cor}_{\g^{(1)}} : X \LA \g^{(1)}$ by
\begin{equation}\label{corrector}
\mathrm{Cor}_{\g^{(1)}}(x):=\log\big(\Phi(x)\big)\big|_{\g^{(1)}} - 
\log\big(\Phi_0(x)\big)\big|_{\g^{(1)}} \qquad (x \in V).
\end{equation}
By periodicities of $\Phi$ and $\Phi_0$, we easily see that 
the set $\{\mathrm{Cor}_{\g^{(1)}}(x) \, | \, x \in V\}$ is a finite set.  
In particular, we find a positive constant $M>0$ such that 
$\max_{x \in \mathcal{F}}\|\mathrm{Cor}_{\g^{(1)}}(x)\|_{\g^{(1)}} \leq M$.

Let $(\ol{\mathcal{Y}}_t^{(n)})_{0 \leq t \leq 1} \, (n \in \mathbb{N})$
be the $G$-valued stochastic processes defined by just replacing 
$\Phi_0$ by $\Phi$ in the definition of $(\mathcal{Y}_t^{(n)})_{0 \leq t \leq 1}$. 
Thanks to several properties of $\mathrm{Cor}_{\g^{(1)}}$, 
we establish the following functional CLT.

\begin{tm}\label{FCLT1-general}
Assume the centered condition {\bf (C)}.
The sequence $\{\ol{\mathcal{Y}}_t^{(n)}\}_{n=1}^\infty$ converges in law
to the $G$-valued diffusion process $(Y_t)_{0 \leq t \leq 1}$ in 
$C_{\bm{1}_G}^{0, \alpha{{\normalfont \hol}}}([0, 1]; G)$ as $n \to \infty$. 
\end{tm}

Let us make comments on our main theorems.
As is emphasized in Breuillard \cite[Section 6]{Bre},
the situation of the non-centered case $\rho_{\mathbb{R}}(\gamma_p) \neq \bm{0}_{\g}$
 is quite different from the centered case $\rho_{\mathbb{R}}(\gamma_p) = \bm{0}_{\g}$
and thus some technical difficulties 
arise to obtain CLTs.
That is why there are few papers which discuss CLTs for non-centered 
random walks on nilpotent Lie groups. 
We obtain, in Theorem \ref{CLT1}, 
a semigroup CLT for the non-centered random walk $\{\xi_n\}_{n=0}^\infty$ on $G$
with a canonical dilation $\tau_{n^{-1/2}}$, while
Cr\'epel--Raugi \cite{CR} and 
Raugi \cite{Raugi} proved similar CLTs for the random walk to (\ref{semigroup CLT1-3}) 
with spatial scalings whose orders are  higher than $\tau_{n^{-1/2}}$.
On the other hand, in the present paper, 
we need to assume the centered condition {\bf (C)}
to obtain a functional CLT (Theorem \ref{FCLT1}) 
for $\{\xi_n\}_{n=0}^\infty$ in the H\"older topology, stronger than
the uniform topology in $C_{\bm{1}_G}([0, 1]; G)$. 
In Appendix A, we mention a method to reduce 
the non-centered case $\rho_{\mathbb{R}}(\gamma_p) \neq \bm{0}_{\g}$
to the centered case by employing a 
measure-change technique based on Alexopoulos \cite{A2}.



\section{Preparations}

\subsection{Limit groups}

Let us review some properties of the limit group. 
For more details, see e.g., Alexopoulos \cite{A1}
and Ishiwata \cite{Ishiwata}. 
We also refer to  
Cr\'epel--Raugi \cite{CR} and 
Goodman \cite{Goodman} for related topics. 
 Let $(G, \cdot)$ be a connected and simply connected nilpotent Lie group 
 of step $r$ and $(\g, [\cdot, \cdot])$ the corresponding Lie algebra. 
 Then the limit group 
 $G_\infty=(G,*)$ of $(G, \cdot)$ is a stratified Lie group of step $r$ 
 and its Lie algebra coincides with $(\g, [\![\cdot, \cdot]\!])$. 
 Namely, 
 the Lie algebra $\frak{g}=\frak{g}^{(1)} \oplus \g^{(2)} \oplus \dots \oplus \frak{g}^{(r)}$
satisfies that $[\![ \g^{(i)}, \g^{(j)}]\!] \subset \g^{(i+j)}$ whenever $i+j \leq r$
and the subspace $\frak{g}^{(1)}$ generates  $\frak{g}$. 
It should be noted that the dilation map $\tau_\ve : G \LA G$ 
is a group automorphism on $(G, *)$ (see \cite[Lemma 2.1]{Ishiwata}). 
We also note that 
the exponential map $\exp : \g_\infty \LA G_\infty$
coincides with the original exponential map $\exp : \g \LA G$. 
Furthermore, for any $g \in G$, 
the inverse element of $g$ in $(G, \cdot)$ 
coincides with the inverse element in $(G, *)$.

We set $d_k=\dim_{\mathbb{R}}\g^{(k)}$ for $k=1, 2, \dots, r$ and $N=d_1+d_2+\cdots+d_r$. 
For $k=1, 2, \dots, r$, we
denote by $\{X_1^{(k)}, X_2^{(k)}, \dots, X_{d_k}^{(k)}\}$
a basis of the subspace $\g^{(k)}$.
We introduce several kinds of 
global coordinate systems in $G$ through $\exp : \g \LA G$.
We write  
$g^{(k)}=(g_1^{(k)}, g_2^{(k)}, \dots, g_{d_k}^{(k)}) \in \mathbb{R}^{d_k}$ for $k=1, 2, \dots, r$. 
We identify the nilpotent Lie group $G$ with $\mathbb{R}^N$
as a differentiable manifold by

\vspace{1mm}
\noindent
$\bullet$
{\it canonical $(\cdot)$-coordinates of the first kind} :
$$
\begin{aligned}
\mathbb{R}^N \ni (g^{(1)}, g^{(2)}, \dots,g^{(r)})\longmapsto 
&\,\,g=\Exp\Big( \sum_{k=1}^r \sum_{i=1}^{d_k}g_i^{(k)}X_i^{(k)}\Big) \in G,
\end{aligned}
$$

\noindent
$\bullet$
{\it canonical $(\cdot)$-coordinates of the second kind} :
$$
\begin{aligned}
\mathbb{R}^N \ni &(g^{(1)}, g^{(2)}, \dots,g^{(r)})\\
\longmapsto 
&\,\,g=\Exp\big( g_{d_r}^{(r)}X_{d_r}^{(r)}\big) \cdot \Exp\big( g_{d_r-1}^{(r)}X_{d_r-1}^{(r)}\big)\cdot
\cdots \cdot\Exp\big( g_{1}^{(r)}X_{1}^{(r)}\big)\\
&\hspace{1cm} \cdot \Exp\big( g_{d_{r-1}}^{(r-1)}X_{d_{r-1}}^{(r-1)}\big) 
\cdot \Exp\big( g_{d_{r-1}-1}^{(r-1)}X_{d_{r-1}-1}^{(r-1)}\big) \cdot\cdots
\cdot \Exp\big( g_{1}^{(r-1)}X_{1}^{(r-1)}\big) \\
&\hspace{1cm}\cdots \cdot \Exp\big( g_{d_1}^{(1)}X_{d_1}^{(1)}\big) \cdot \Exp\big( g_{d_1-1}^{(1)}X_{d_1-1}^{(1)}\big)
\cdot \cdots \cdot\Exp\big( g_{1}^{(1)}X_{1}^{(1)}\big) \in G,
\end{aligned}
$$

\noindent
$\bullet$
{\it canonical $(*)$-coordinates of the second kind} :
$$
\begin{aligned}
\mathbb{R}^N \ni &(g_*^{(1)}, g_*^{(2)}, \dots,g_*^{(r)})\\
\longmapsto 
&\,\,g=\Exp\big( g_{d_r*}^{(r)}X_{d_r}^{(r)}\big) * \Exp\big( g_{d_r-1*}^{(r)}X_{d_r-1}^{(r)}\big)*
\cdots *\Exp\big( g_{1*}^{(r)}X_{1}^{(r)}\big)\\
& \hspace{1cm}* \Exp\big( g_{d_{r-1}*}^{(r-1)}X_{d_{r-1}}^{(r-1)}\big) 
* \Exp\big( g_{d_{r-1}-1*}^{(r-1)}X_{d_{r-1}-1}^{(r-1)}\big) *\cdots 
*\Exp\big( g_{1*}^{(r-1)}X_{1}^{(r-1)}\big)\\
&\hspace{1cm}*\cdots * \Exp\big( g_{d_1*}^{(1)}X_{d_1}^{(1)}\big) * \Exp\big( g_{d_1-1*}^{(1)}X_{d_1-1}^{(1)}\big)
* \cdots *\Exp\big( g_{1*}^{(1)}X_{1}^{(1)}\big) \in G_\infty.
\end{aligned}
$$
We give the relations between the deformed product and the given product 
on $G$ as an easy application of the Campbell--Baker--Hausdorff (CBH) formula 
\begin{equation}\label{CBH-formula}
\log\big(\exp(Z_1) \cdot \exp(Z_2)\big)
=Z_1+Z_2+\frac{1}{2}[Z_1, Z_2]+\cdots \qquad(Z_1, Z_2 \in \g). 
\end{equation}
The following is straightforward from the definition of the deformed product.
\begin{equation}\label{rel-1}
\Log(g * h)\big|_{\g^{(k)}}=\Log(g \cdot h)\big|_{\g^{(k)}} \qquad (g, h \in G, \, k=1, 2).
\end{equation}
We notice that the relation above does not hold in general for $k=3, 4, \dots, r$. 
The following identities give us a comparison between $(\cdot)$-coordinates and $(*)$-coordinates. 
For $g \in G$, we have the following.
\begin{align}\label{rel-2}
g_{i*}^{(k)}&=g_i^{(k)} \quad (i=1, 2, \dots, d_k, \, k=1, 2), \\
g_{i*}^{(k)}&=g_{i}^{(k)}+\sum_{0 < |K| \leq k-1}C_K \mathcal{P}^K(g)
\qquad (i=1, 2, \dots, d_k, \, k=3, 4, \dots, r) \label{rel-3}
\end{align}
for some constant $C_K$, where $K$ stands for a multi-index 
$\big((i_1, k_1), (i_2, k_2), \dots, (i_\ell, k_\ell)\big)$ 
with length $|K|:=k_1+k_2+\dots+k_\ell$
and $\mathcal{P}^K(g):=g_{i_1}^{(k_1)} \cdot g_{i_2}^{(k_2)} \cdots g_{i_\ell}^{(k_\ell)}$. 
The invariances (\ref{rel-1}) and (\ref{rel-2}) play an important role 
to obtain main results. 
For $g, h \in G$, we also have
\begin{align}\label{rel-4}
(g*h)_{i*}^{(k)}&=(g\cdot h)_i^{(k)} \qquad (i=1, 2, \dots, d_k, \, k=1, 2), \\
(g*h)_{i*}^{(k)}&=(g\cdot h)_{i}^{(k)}+\sum_{\substack{ |K_1|+|K_2|\leq k-1 \\ |K_2|>0}}
C_{K_1, K_2} \mathcal{P}^{K_1}_*(g)\mathcal{P}^{K_2}(g \cdot h)\nn\\
&\hspace{4cm}
\qquad (i=1, 2, \dots, d_k, \, k=3, 4, \dots, r) \label{rel-5}
\end{align}
by using (\ref{rel-2}) and (\ref{rel-3}). 
See \cite[Section 2]{Ishiwata} for more details.

\subsection{Carnot--Carath\'eodory metric and homogeneous norms}

As is well-known, a nilpotent Lie group $G$ is a candidate
of the typical sub-Riemannian manifolds, 
which is a certain generalization of a Riemannian manifold.  
The notion of the Carnot--Carath\'eodory metric naturally appears
when we investigate distances between two points in $G$.
It is an important intrinsic metric in this context and  
is  degenerate in the sense that we only go along curves 
which are tangent to a ``horizontal subspace''
of the tangent space of $G$. 
We discuss several properties of 
the Carnot--Carath\'eodory metric on a nilpotent Lie group $G$ in this subsection.   
Note that the definition of such an intrinsic metric in more general setting 
is found in some references. 
See e.g., 
Varopoulos--Saloff-Coste--Coulhon \cite{VSC}. 

Recall that the {\it Carnot--Carath\'eodory metric} on $G$ is defined by (\ref{def-dcc}).
We know that the subspace $\frak{g}^{(1)}$ satisfies 
the so-called {\it H\"ormander condition} in $\frak{g}$,
that is,  $L_{\g^{(1)}}(g)=T_gG$ for any $g \in G$, 
where $L_{\g^{(1)}}(g)$ denotes the evaluation of $\g^{(1)}$ at $g \in G$. 
The Carnot-Carath\'eodory metric is then well-defined in the sense that
$d_{\mathrm{CC}}(g, h)<\infty$ for every $g, h \in G$, 
thanks to the H\"ormander condition on $\frak{g}^{(1)}$. 
Furthermore, the topology induced by the Carnot-Carath\'eodory 
metric $d_{\mathrm{CC}}$ coincides with the original one of $G$. 
We emphasize that $d_{\mathrm{CC}}$ 
is behaved well under dilations. More precisely, we have
\begin{equation}\label{CC-dilation}
d_{\mathrm{CC}}\big(\tau_\ve(g), \tau_\ve(h)\big)=\ve d_{\mathrm{CC}}(g, h) 
\qquad (\ve \geq 0, \, g, h \in G).
\end{equation}

We now present the notion of 
{\it homogeneous norm} on $G$.  
The one-parameter group of dilations $(\tau_\ve)_{\ve \geq 0}$ 
allows us to consider
scalar multiplications on nilpotent Lie groups. 
We replace the usual Euclidean norms 
by the following functions. 
A continuous function $\|\cdot\| : G \LA [0, \infty)$ is called a homogeneous norm on $G$ if
(i) $\|g\|=0$ if and only if $g=\bm{1}_G$ and 
(ii) 
$\|\tau_\ve g\|=\ve \|g\|$ for $\ve \geq 0$ and $g \in G$.
One of the typical examples of homogeneous norms is given by the 
Carnot--Carath\'eodory metric $d_{\mathrm{CC}}$. 
We define a continuous function $\|\cdot\|_{\mathrm{CC}} : G \LA [0, \infty)$ by
$\|g\|_{\mathrm{CC}}:=d_{\mathrm{CC}}(\bm{1}_G, g)$
for $g \in G$. 
Then 
$\|\cdot\|_{\mathrm{CC}}$ is a homogeneous norm on $G$ in view of (\ref{CC-dilation}). 
Another basic homogeneous norm is given in the following way.  
We denote by $\{X_1^{(k)}, X_2^{(k)}, \dots, X_{d_k}^{(k)}\}$ a basis in $\g^{(k)}$ for $k=1, 2, \dots, r$. 
We introduce a norm $\|\cdot\|_{\g^{(k)}}$ on $\g^{(k)}$ by the usual Euclidean norm.
If $Z \in \g$ is decomposed as $Z=Z^{(1)}+Z^{(2)}+\cdots +Z^{(r)} \, (Z^{(k)} \in \g^{(k)})$, 
we define a function $\|\cdot\|_{\g} : \g \LA [0, \infty)$ by
$\|Z\|_{\g}:=\sum_{k=1}^r \|Z^{(k)}\|_{\g^{(k)}}^{1/k}.$
We set
$\|g\|_{\Hom}:=\|\Log(g)\|_{\g}$ for $g \in G$.
We then observe that $\|\cdot\|_{\Hom}$ is a homogeneous norm on $G$. 
The homogenuity (ii) leads to the most important fact 
that all homogeneous norms on $G$ are equivalent, 
which is similar to the case of norms on the Euclidean space. 
More precisely, we have the following:

\begin{pr}\label{homogeneous equiv}{\bf (cf.~Goodman~\cite{Goodman})}
If $\|\cdot\|_1$ and $\|\cdot\|_2$ are two homogeneous norms on $G$, then there exists a constant
$C>0$ such that
$C^{-1}\|g\|_1 \leq \|g\|_2 \leq C\|g\|_1$ for $g \in G$.
\end{pr}
This proposition plays a crucial role to obtain Theorem \ref{FCLT1}. 
For more details on homogeneous norms, we also refer to 
Bonfiglioli--Lanconelli--Uguzzoni \cite{BLU}.

\subsection{Discrete geometric analysis}
We present some basics of discrete geometric analysis
on graphs due to Kotani--Sunada \cite{KS06} or Sunada \cite{S0, S, S2}.
We consider a finite graph $X_0=(V_0, E_0)$ 
and an irreducible random walk on $X_0$ associated 
with a non-negative transition probability 
$p : E_0 \LA [0, 1]$. 
We define the 0-chain group, 1-chain group, 0-cochain group and 1-cochain group by
$$
\begin{aligned}
C_0(X_0, \mathbb{R})&:=\Big\{ \sum_{x \in V_0}a_x x \, \Big| \, a_x \in \mathbb{R}\Big\}, &
C_1(X_0, \mathbb{R})&:=\Big\{ \sum_{e \in E_0}a_e e \, \Big| \, a_e \in \mathbb{R}, \, \ol{e}=-e\Big\}, \\
C^0(X_0, \mathbb{R})&:=\{f : V_0 \LA \mathbb{R}\}, &
C^1(X_0, \mathbb{R})&:=\{\omega : E_0 \LA \mathbb{R} \, | \, \omega(\ol{e})=-\omega(e)\},
\end{aligned}
$$
respectively. 
An element of $C^1(X_0, \mathbb{R})$ is called a 1-form on $X_0$. 
The boundary operator $\del : C_1(X_0, \mathbb{R}) \LA C_0(X_0, \mathbb{R})$ 
and the difference operator $d : C^0(X_0, \mathbb{R}) \LA C^1(X_0, \mathbb{R})$ 
are defined by $\del(e)=t(e)-o(e)$ for $e \in E_0$ and 
$df(e)=f\big(t(e)\big)-f\big(o(e)\big)$ for $e \in E_0$, respectively. 
Then, the first homology group $\h_1(X_0, \mathbb{R})$ and 
the first cohomology group $\h^1(X_0, \mathbb{R})$ are defined by 
$\Ker(\del) \subset C_1(X_0, \mathbb{R})$ and $C^1(X_0, \mathbb{R})/\im(d)$, respectively. 
We write $L : C^0(X_0, \mathbb{R}) \LA C^0(X_0, \mathbb{R})$ for 
the transition operator associated with $p$. 
We define a special 1-chain  by
$$
\gamma_p:=\sum_{e \in E_0}\widetilde{m}(e)e \in C_1(X_0, \mathbb{R}).
$$
It is easily seen that $\del(\gamma_p)=0$ so that $\gamma_p \in \h_1(X_0, \mathbb{R})$. 
Furthermore, it is clear that 
the random walk on $X_0$ is ($m$-)symmetric if and only if $\gamma_p=0$.
The 1-cycle $\gamma_p$ is called the {\it homological direction} of the given random walk on $X_0$. 
A simple application of the ergodic theorem leads to
the law of large numbers on $C_1(X_0, \mathbb{R})$.
$$
\lim_{n \to \infty}\frac{1}{n}(e_1+e_2+\dots+e_n)=\gamma_p, 
\quad \mathbb{P}_x\text{-a.e.} \, c=(e_1, e_2, \dots, e_n, \dots) \in \Omega_x(X_0).
$$
A 1-form $\omega \in C^1(X_0, \mathbb{R})$ is said to be {\it modified harmonic} if
\begin{equation}\label{form-harmonicity}
\sum_{e \in (E_0)_x}p(e)\omega(x)-\la \gamma_p, \omega \ra=0 \qquad (x \in V_0),
\end{equation}
where
$\la \gamma_p, \omega \ra:=
{}_{C_1(X_0, \mathbb{R})}\la \gamma_p, \omega \ra_{C^1(X_0, \mathbb{R})}$ 
is constant as a function on $V_0$. 
We denote by $\mathcal{H}^1(X_0)$ the space of modified harmonic 1-forms
and equip it with the inner product
$$
\La \omega_1, \omega_2 \Ra_p
:=\sum_{e \in E_0}\widetilde{m}(e)\omega_1(e)\omega_2(e)
 - \la \gamma_p, \omega_1 \ra \la \gamma_p, \omega_2 \ra 
 \qquad \big(\omega_1, \omega_2 \in \mathcal{H}^1(X_0)\big)
$$
associated with the transition probability $p$. 
We may identify $\h^1(X_0, \mathbb{R})$
with $\mathcal{H}^1(X_0)$ 
by the discrete Hodge-Kodaira theorem (cf.~\cite[Lemma 5.2]{KS06}).
We induce an inner product 
from $\h^1(X_0, \mathbb{R})$ by using this identification. 

Let $\Gamma$ be a torsion free, finitely generated nilpotent group of step $r$.
Then a $\Gamma$-nilpotent covering graph $X=(V, E)$ 
is defined by the $\Gamma$-covering of $X_0$.
Let $p : E \LA [0, 1]$ and $m : V \LA (0, 1]$ be the $\Gamma$-invariant lifts of  
$p : E_0 \LA [0, 1]$ and $m : V_0 \LA (0, 1]$, respectively. 
Denote by $\widehat{\pi} : G \LA G/[G, G]$ the canonical projection. 
Since $\Gamma$ is a cocompact lattice in $G$,
the subset $\widehat{\pi}(\Gamma) \subset G/[G, G]$
 is also a lattice in $G/[G, G] \cong \g^{(1)}$ (cf.~Malc\'ev \cite{Malcev} and Raghunathan \cite{Rag}). 
We take the canonical surjective homomorphism 
$\rho : \h_1(X_0, \mathbb{Z}) \LA \widehat{\pi}(\Gamma) \cong \Gamma/[\Gamma, \Gamma]$
and its realification is denoted by
$\rho_{\mathbb{R}} : \h_1(X_0, \mathbb{R})  \LA \widehat{\pi}(\Gamma) \otimes \mathbb{R}.$
We identify $\Hom(\widehat{\pi}(\Gamma), \mathbb{R})$ with a subspace 
of $\h^1(X_0, \mathbb{R})$ by using the transposed map ${}^t \rho_{\mathbb{R}}$.
We restrict  $\La \cdot, \cdot \Ra_p$ on $\h^1(X_0, \mathbb{R})$ 
to the subspace $\Hom(\widehat{\pi}(\Gamma), \mathbb{R})$ 
and take it up the dual inner product $\la \cdot , \cdot \ra_{alb}$ 
on $\widehat{\pi}(\Gamma) \otimes \mathbb{R}$. 
Then, a flat metric $g_0$ on $\g^{(1)}$ is induced 
and we call it the {\it Albanese metric} on $\g^{(1)}$. 
This procedure can be summarized as follows:

$$
\xymatrix{ 
(\frak{g}^{(1)}, g_0)  \ar @{<->}[d]^{\mathrm{dual}}
&\hspace{-1.5cm}\cong 
& \hspace{-1.5cm}\qquad\widehat{\pi}(\Gamma) \otimes \mathbb{R}
 \ar @{<<-}[r]^{\rho_{\mathbb{R}}} \ar @{<->}[d]^{\mathrm{dual}} 
 & \h_1(X_0, \mathbb{R}) \ar @{<->}[d]^{\mathrm{dual}} &\\
\Hom(\frak{g}^{(1)}, \mathbb{R}) 
&\hspace{-1.5cm}\cong 
& \hspace{-1.5cm}\quad\Hom(\widehat{\pi}(\Gamma), \mathbb{R})  
\ar @{^{(}->}[r]_{{}^t \rho_{\mathbb{R}}} & \h^1(X_0, \mathbb{R})   
&\hspace{-0.9cm}\cong  \big(\mathcal{H}^1(X_0), \La \cdot , \cdot \Ra_p\big).}
$$

A map $\Phi : X \LA G$ is said to be a {\it periodic realization} 
of $X$ when it satisfies
$\Phi(\gamma x)=\gamma \cdot \Phi(x)$ for $\gamma \in \Gamma$ and $x \in V.$
Fix a reference point $x_* \in V$
and define a special realization $\Phi_0 : X \LA G$ by
\begin{equation}\label{Alb map}
{}_{\Hom(\g^{(1)}, \mathbb{R})}\big\la \omega, 
\Log\big(\Phi_0(x)\big)\big|_{\g^{(1)}} \big\ra_{\g^{(1)}}=\int_{x_*}^x \widetilde{\omega} 
\qquad \big(\omega \in \Hom(\g^{(1)}, \mathbb{R}), \, x \in V\big),
\end{equation}
where $\widetilde{\omega}$ is the lift of 
$\omega={}^t \rho_{\mathbb{R}}(\omega) \in \h^1(X_0, \mathbb{R})$ to $X$. 
Here 
$
\int_{x_*}^x \widetilde{\omega}=\int_c \widetilde{\omega}=\sum_{i=1}^n \widetilde{\omega}(e_i)
$
for a path $c=(e_1, e_2, \dots, e_n)$ with $o(e_1)=x_*$ and $t(e_n)=x$. 
We note that this line integral does not depend on the choice of a path $c$. 
The following lemma asserts that such $\Phi_0$ enjoys the {\it modified harmonicity} 
in the sense of (\ref{m-harmonicity}).

\begin{lm} 
The periodic realization $\Phi_0 : X \LA G$ 
defined by {\normalfont{(\ref{Alb map})}} is the modified harmonic realization, that is, 
$$
\sum_{e \in E_x}p(e)\Log\Big(\Phi_0\big(o(e)\big)^{-1} \cdot  \Phi_0\big(t(e)\big)\Big)\Big|_{\g^{(1)}}
=\rho_{\mathbb{R}}(\gamma_p) \qquad (x \in V).
$$
\end{lm} 

\noindent  
{\bf Proof.} For each 
$\omega = {}^t \rho_{\mathbb{R}}(\omega) 
\in \h^1(X_0, \mathbb{R}) \cong \mathcal{H}^1(X_0)$ and $x \in V$, 
Equation (\ref{Alb map}) yields
$$
\begin{aligned}
&{}_{\Hom(\g^{(1)}, \mathbb{R})}\Big\la \omega, 
\sum_{e \in E_x}p(e)
\Log\Big(\Phi_0\big(o(e)\big)^{-1} \cdot  
              \Phi_0\big(t(e)\big)\Big)\Big|_{\g^{(1)}} \Big\ra_{\g^{(1)}}\\
&=\sum_{e \in E_x}p(e) 
{}_{\Hom(\g^{(1)}, \mathbb{R})}\Big\la 
\omega, \Log\big(\Phi_0\big(t(e)\big)\big)\big|_{\g^{(1)}} -
              \Log \big(\Phi_0\big(o(e)\big)\big)\Big|_{\g^{(1)}} \Big\ra_{\g^{(1)}}\\
&=\sum_{e \in E_x}p(e)\widetilde{\omega}(e)\\
&=-(\delta_p\omega)\big(\pi(x)\big)\\
&=\la \gamma_p, \omega \ra
={}_{\Hom(\g^{(1)}, \mathbb{R})}\big\la \omega, \rho_{\mathbb{R}}(\gamma_p) \big\ra_{\g^{(1)}}.
\end{aligned}
$$
This gives the desired equation (\ref{m-harmonicity}). \qed

\vspace{1mm}

\subsection{The Markov chain on $G$}

Let us consider a time-homogeneous Markov chain 
$(\Omega_x(X), \mathbb{P}_x, \{w_n\}_{n=0}^\infty)$ 
with values in a $\Gamma$-nilpotent covering graph $X$. 
We denote by $\Phi : X \LA G$ a $\Gamma$-equivariant realization of $X$. 
We then have the $G$-valued Markov chain 
$(\Omega_x(X), \mathbb{P}_x, \{\xi_n\}_{n=0}^\infty)$ defined by
$\xi_n(c):=\Phi\big(w_n(c)\big)$ for $n \in \mathbb{N} \cup \{0\}$ and 
$c \in \Omega_x(X)$,
through the map $\Phi$.
This gives rise to the $\g$-valued random walk 
$\Xi_n(c):=\Log\big(\xi_n(c)\big)=\Log\big(\Phi\big(w_n(c)\big) \big)$ for
$n \in \mathbb{N} \cup\{0\}$ and $c \in \Omega_x(X).$
We obtain the following law of large numbers on $\g^{(1)}$ 
by the ergodic theorem.
\begin{equation}\label{LLN}
\lim_{n \to \infty}\frac{1}{n}\Xi_n(\cdot)\big|_{\g^{(1)}}
=\rho_{\mathbb{R}}(\gamma_p), \quad \mathbb{P}_x\text{-a.s.}
\end{equation}

It is known that the notion of martingales plays a crucial role in the theory of stochastic processes. 
We give a certain characterization of modified harmonic realizations 
in view of martingale theory. 
Let $\pi_n : \Omega_x(X) \LA \Omega_{x, n}(X) \, (n \in \mathbb{N} \cup \{0\})$ be a projection 
defined by
$\pi_n(c):=(e_1, e_2, \dots, e_n)$ for $c=(e_1, e_2, \dots, e_n, \dots) \in \Omega_x(X)$.
Denote by $\{\mathcal{F}_n\}_{n=0}^\infty$ the filtration 
such that $\mathcal{F}_0=\{\emptyset, \Omega_x(X)\}$ and
$\mathcal{F}_n:=\sigma\big( \pi_n^{-1}(A) \, \big| \, 
A \subset \Omega_{x, n}(X)\big)$ for $n \in \mathbb{N}$. 
We mention that $\mathcal{F}_n$ is a sub-$\sigma$-algebra of 
$\mathcal{F}_\infty:=\bigvee_{n=0}^\infty \mathcal{F}_n$ for $n \in \mathbb{N}$. 
We will use the following lemma in the proof of Lemma \ref{sublemma1}. 

\begin{lm} \label{martingale}
Let $\{X_1^{(1)}, X_2^{(1)}, \dots, X_{d_1}^{(1)}\}$ be a basis of $\g^{(1)}$. 
Then a periodic realization $\Phi_0 : X \LA G$ is 
the modified harmonic realization
if and only if the stochastic process 
$$
\big\{ \Xi_n\big|_{X_i^{(1)}} - n \rho_{\mathbb{R}}(\gamma_p)
\big|_{X_i^{(1)}}\big\}_{n=0}^\infty \qquad (i=1, 2, \dots, d_1),
$$
with values in $\mathbb{R}$, is an $\{\mathcal{F}_n\}$-martingale. 
\end{lm}

\noindent
{\bf Proof.} Suppose that $\Phi_0$ is modified harmonic.  
For $n \in \mathbb{N}$ and $A \in \mathcal{F}_n$, we have
$$
\begin{aligned}
& \mathbb{E}^{\mathbb{P}_{x}}\Big[ \Xi_{n+1}\big|_{X_i^{(1)}} 
- (n+1)\rho_{\mathbb{R}}(\gamma_p)\big|_{X_i^{(1)}} \, ; \, A \Big]\\
&=\sum_{c \in \Omega_x(X)}\!\!\!p(c)
\Big\{\Log\Big(\Phi_0\big(t(e_{n+1})\big)\Big)\big|_{X_i^{(1)}} 
- (n+1)\rho_{\mathbb{R}}(\gamma_p)\big|_{X_i^{(1)}} \Big\}\bm{1}_{A}(c)\\
&=\sum_{c' \in \Omega_{x, n}(X)}\!\!\!\!\!p(c')\sum_{e \in E_{t(c')}}\!\!p(e)
\Bigg[\Big\{\Log\Big(\Phi_0\big(t(e)\big)\Big)\big|_{X_i^{(1)}} 
- \rho_{\mathbb{R}}(\gamma_p)\big|_{X_i^{(1)}}\Big\}
- n\rho_{\mathbb{R}}(\gamma_p)\big|_{X_i^{(1)}} \Bigg]\bm{1}_{A}(c'),
\end{aligned}
$$
where $\mathbb{E}^{\mathbb{P}_x}$ stands for the expectation 
with respect to the probability measure $\mathbb{P}_x$. 
In terms of the modified harmonicity of $\Phi_0$, this is equal to 
$$
\begin{aligned}
& \sum_{c' \in \Omega_{x, n}(X)}
p(c')\Big\{\log\Big(\Phi_0\big(o(e_{n+1})\big)\Big)\big|_{X_i^{(1)}} 
- n\rho_{\mathbb{R}}(\gamma_p)\big|_{X_i^{(1)}}\Big\}\bm{1}_{A}(c')\\
&=\mathbb{E}^{\mathbb{P}_{x}}\Big[ \Xi_{n}\big|_{X_i^{(1)}} 
- n\rho_{\mathbb{R}}(\gamma_p)\big|_{X_i^{(1)}} \, ; \, A \Big]
\end{aligned}
$$
Thus it follows that the process 
$\big\{ \Xi_n\big|_{X_i^{(1)}} 
- n \rho_{\mathbb{R}}(\gamma_p)\big|_{X_i^{(1)}}\big\}_{n=0}^\infty $ 
is an $\{\mathcal{F}_n\}$-martingale. 
The converse is obvious from the argument above. \qed



\section{Proof of main results}

The aim of this section is to prove main results
(Theorems \ref{CLT1}, \ref{FCLT1} and \ref{FCLT1-general}). 
In what follows, we set 
$$
\begin{aligned}
d\Phi_0(e)&=\Phi_0\big(o(e)\big)^{-1} \cdot \Phi_0\big(t(e)\big) \qquad (e \in E), \\
\|d\Phi_0\|_\infty&=\max_{e \in E_0}\Big\{ 
\big\|\log\big(d\Phi_0(\widetilde{e})\big)\big|_{\g^{(1)}}\big\|_{\g^{(1)}}
+ \big\|\log\big(d\Phi_0(\widetilde{e})\big)\big|_{\g^{(2)}}\big\|_{\g^{(2)}}^{1/2}\Big\},
 \end{aligned}
$$
where $\widetilde{e}$ stands for a lift of $e \in E_0$ to $X$. 
We should mention that 
\begin{equation}\label{N^k}
\Big( \Phi_0(x)^{-1} \cdot \Phi_0\big(t(c)\big)\Big)_i^{(k)} =O(N^k)
\end{equation}
for $x \in V, c \in \Omega_{x, N}(X),\, i=1, 2, \dots, d_k$ and $k=1, 2, \dots, r$. 
We also write 
$
\rho=\rho_{\mathbb{R}}(\gamma_p)$ and 
$e^{z \rho}=\exp\big( z\rho_{\mathbb{R}}(\gamma_p)\big) \, (z \in \mathbb{R})$
for brevity. 
We give an important property of the family of approximation operators 
$(\mathcal{P}_\ve)_{0 \leq \ve \leq 1}$ defined by (\ref{scaling}). 

\begin{lm} 
Let $q>1$. Then 
$\big(\big(C_{\infty, q}(X \times \mathbb{Z}), \|\cdot\|_{\infty, q}; \mathcal{P}_\ve\big)\big)_{0 \leq \ve \leq 1}$ 
is a family of Banach spaces approximating to 
the Banach space $\big(C_\infty(G), \|\cdot\|_\infty^G\big)$ in the sense of 
Trotter {\normalfont \cite{Trotter}:}
$$
\|\mathcal{P}_\ve f\|_{\infty, q} \leq \|f\|_\infty^G \quad \text{and}
\quad \lim_{\ve \searrow 0}
\|\mathcal{P}_\ve f\|_{\infty, q} =\|f\|_\infty^G \qquad \big(f \in C_\infty(G)\big).
$$
\end{lm}

\noindent
{\bf Proof.} The former assertion follows from 
$$
\|\mathcal{P}_\ve f\|_{\infty, q} = \frac{1}{C_q}\sum_{z \in  \mathbb{Z}}
\frac{\| f(\cdot, z)\|_\infty}{1+|z|^{q}} \leq  \frac{1}{C_q}\sum_{z \in  \mathbb{Z}}
\frac{\|f\|_\infty}{1+|z|^{q}}=\|f\|_\infty.
$$
We prove the latter one.
 Let $g_0 \in G$ be an element which attains $\|f\|_\infty=\sup_{g \in G}|f(g)|$. 
We fix $z \in \mathbb{Z}$. Then we have
$$
\begin{aligned}
\|\mathcal{P}_\ve f(\cdot, z)\|_{\infty} 
\geq |f(g_0)| - \inf_{x \in X}\Big| f(g_0) - f\Big( \tau_\ve \big(\Phi_0(x)
* \exp(-z\rho_{\mathbb{R}}(\gamma_p))\big)\Big)\Big|.
\end{aligned}
$$
On the other hand, we have
$$
\begin{aligned}
&\hspace{0.5cm}\inf_{x \in X}d_{\mathrm{CC}}\Big( g_0,  \tau_\ve \big(\Phi_0(x)
* \exp(-z\rho_{\mathbb{R}}(\gamma_p))\big)\Big)\\
&= \ve \inf_{x \in X}d_{\mathrm{CC}}\Big( \tau_{1/\ve}(g_0),  
\Phi_0(x)* \exp(-z\rho_{\mathbb{R}}(\gamma_p))\Big)<\ve M
\end{aligned}
$$
for some $M=M(z)>0$. From the continuity of $f$, for any $\delta>0$, there exists $\delta'>0$ such that 
$d_{\mathrm{CC}}(g_0, h)<\delta'$ implies $|f(g_0) - f(h)|<\delta$. 
By choosing a sufficiently small $\ve>0$, we have 
$$
d_{\mathrm{CC}}\Big(g_0, \tau_\ve\big( \Phi_0(x_*)
* \exp(-z\rho_{\mathbb{R}}(\gamma_p)) \big) \Big)<\delta'
$$
for some $x_* \in X$. Then we have
$$
\begin{aligned}
&\hspace{0.5cm}\inf_{x \in X}\Big| f(g_0) - f\Big( \tau_\ve \big(\Phi_0(x)
* \exp(-z\rho_{\mathbb{R}}(\gamma_p))\big)\Big)\Big|\\
&\leq \Big| f(g_0) - 
f\Big( \tau_\ve \big(\Phi_0(x_*)* \exp(-z\rho_{\mathbb{R}}(\gamma_p))\big)\Big)\Big|<\delta
\end{aligned}
$$
and this implies $\lim_{\ve \searrow 0}\|\mathcal{P}_\ve f(\cdot, z)\|_\infty=\|f\|_\infty$ 
for $z \in \mathbb{Z}$. 
By using the dominated convergence theorem, we obtain
$\lim_{\ve \searrow 0}\|\mathcal{P}_\ve f\|_{\infty, r}=\|f\|_\infty$.
This completes the proof. \qed

\subsection{Proof of Theorem \ref{CLT1}}

The following lemma is significant to prove Theorem \ref{CLT1}.

\begin{lm}\label{key-lem}
Let $f \in C_0^\infty(G)$ and $q > 4r+1$. Then we have 
$$
\Big\|\frac{1}{N\ve^2}\big(I - \mathcal{L}_p^N\big)\mathcal{P}_\ve f 
- \mathcal{P}_\ve \A f \Big\|_{\infty, q} \LA 0
$$
as $N \to \infty, \, \ve \searrow 0$ and $N^2\ve \searrow 0$, 
where $\mathcal{L}_p$ is the transition-shift operator defined by {\normalfont{(\ref{transition-shift})}}
and 
$\A$ is the sub-elliptic operator defined by \normalfont{(\ref{generator})}.
\end{lm}

\noindent
{\bf Proof.} 
We divide the proof into several steps. 

\vspace{1mm}
\noindent
{\bf Step~1}. 
We first apply Taylor's formula (cf.~Alexopoulos \cite[Lemma 5.3]{A2})  
for the ($*$)-coordinates of the second kind to $f \in C_0^\infty(G)$ 
at $\tau_\ve\big(\Phi_0(x)* e^{-z\rho}\big)  \in G$. 
By recalling that $(G, *)$ is a stratified Lie group, we have
\begin{align}\label{eq1}
& \frac{1}{N\ve^2} (I - \mathcal{L}_p^N)\mathcal{P}_\ve f(x, z) \nn \\
&=-\sum_{(i, k)}\frac{\ve^{k-2}}{N}X_{i*}^{(k)}f\big(\tau_\ve\big(\Phi_0(x) * e^{-z\rho}\big)\big)
\sum_{c \in \Omega_{x, N}(X)}p(c)\big( \mathcal{B}_N(x, z, c)\big)_{i*}^{(k)}\nn\\
&\hspace{0.45cm}-\Big(\sum_{(i_1, k_1)  \geq (i_2, k_2)}
\frac{\ve^{k_1+k_2-2}}{2N}X_{i_1*}^{(k_1)}X_{i_2*}^{(k_2)}+
\sum_{(i_2, k_2)  > (i_1, k_1)}
\frac{\ve^{k_1+k_2-2}}{2N} X_{i_2*}^{(k_2)}X_{i_1*}^{(k_1)} \Big)\nn\\
&\hspace{1.3cm}\times f\Big(\tau_\ve\big(\Phi_0(x) * e^{-z\rho}\big)\Big)
\sum_{c \in \Omega_{x, N}(X)}p(c) \big( \mathcal{B}_N(x, z, c)\big)_{i_1*}^{(k_1)}
\big( \mathcal{B}_N(x, z, c)\big)_{i_2*}^{(k_2)}\nn\\
&\hspace{0.45cm} - \sum_{(i_1, k_1), (i_2, k_2), (i_3, k_3)}\frac{\ve^{k_1+k_2+k_3-2}}{6N}
\frac{\del^3 f}{\del g_{i_1*}^{(k_1)}\del g_{i_2*}^{(k_2)} \del g_{i_3*}^{(k_3)}}(\theta) 
\sum_{c \in \Omega_{x, N}(X)}p(c) \big( \mathcal{B}_N(x, z, c)\big)_{i_1*}^{(k_1)}\nn\\
&\hspace{1.3cm}\times
\big( \mathcal{B}_N(x, z, c)\big)_{i_2*}^{(k_2)}
\big( \mathcal{B}_N(x, z, c)\big)_{i_3*}^{(k_3)} \qquad (x \in V, \, z \in \mathbb{Z}),
\end{align}
for some $\theta \in G$ with
$|\theta_{i_*}^{(k)}| \leq \ve^k\big|\big(\mathcal{B}_N(x, z, c)\big)_{i*}^{(k)}\big|$
for $i=1, 2, \dots, d_k$ and $ k=1, 2, \dots, r$,
where
the summation $\sum_{(i_1, k_1) \geq (i_2, k_2)}$ runs over
all $(i_1, k_1)$ and $(i_2, k_2)$ with $k_1>k_2$ or $k_1=k_2, \, i_1 \geq i_2$. 
Here we set
$$
\mathcal{B}_N(x, z, c):=e^{z\rho}*\Phi_0(x)^{-1} * \Phi_0\big(t(c)\big) * e^{-(z+N)\rho} 
\quad \big(N \in \mathbb{N}, \, x \in V, \, z \in \mathbb{Z}, \, c \in \Omega_{x, N}(X)\big). 
$$ 
We denote by $\mathrm{Ord}_\ve(k)$ the terms of 
the right-hand side of (\ref{eq1}) whose order of $\ve$ equals just 
$k$. Then (\ref{eq1}) is rewritten as
$$
\frac{1}{N\ve^2}(I-\mathcal{L}_p^N)\mathcal{P}_\ve f(x, z)
=\mathrm{Ord}_\ve(-1) + \mathrm{Ord}_\ve(0)+\sum_{k \geq 1} \mathrm{Ord}_\ve(k)
 \qquad (x \in V, \, z \in \mathbb{Z}), 
$$
where 
$$
\begin{aligned}
\mathrm{Ord}_\ve(-1) &= - \frac{1}{N\ve} \sum_{i=1}^{d_1}
X_{i*}^{(1)}f\big(\tau_\ve\big(\Phi_0(x) * e^{-z\rho}\big)\big)
\sum_{c \in \Omega_{x, N}(X)}p(c)\big( \mathcal{B}_N(x, z, c)\big)_{i*}^{(1)}, \\
\mathrm{Ord}_\ve(0) &= -\frac{1}{N} \sum_{i=1}^{d_2} X_{i*}^{(2)}f\big(\tau_\ve\big(\Phi_0(x) * e^{-z\rho}\big)\big)
\sum_{c \in \Omega_{x, N}(X)}p(c)\Big\{ \big( \mathcal{B}_N(x, z, c)\big)_{i*}^{(2)}\\
&\hspace{1cm}-\frac{1}{2}\sum_{1 \leq \lambda <\nu \leq d_1}\big( \mathcal{B}_N(x, z, c)\big)_{\lambda*}^{(1)}
\big( \mathcal{B}_N(x, z, c)\big)_{\nu*}^{(1)} [\![X_\lambda^{(1)}, X_\nu^{(1)}]\!]\big|_{X_i^{(2)}} \Big\}\\
&\hspace{0.5cm}-\frac{1}{2N} \sum_{1 \leq i, j \leq d_1}
X_{i*}^{(1)}X_{j*}^{(1)}f\big(\tau_\ve\big(\Phi_0(x) * e^{-z\rho}\big)\big)\\
 &\hspace{1cm}\times\sum_{c \in \Omega_{x, N}(X)}p(c) \big( \mathcal{B}_N(x, z, c)\big)_{i*}^{(1)}
\big( \mathcal{B}_N(x, z, c)\big)_{j*}^{(1)}
\end{aligned}
$$
and $\sum_{k \geq 1} \mathrm{Ord}_\ve(k)$ is given by the sum of the following three parts:
$$
\begin{aligned}
\mathcal{I}_1(\ve, N) &=- \sum_{k \geq 3} \sum_{i=1}^{d_k}
\frac{\ve^{k-2}}{N} X_{i*}^{(k)}f\big( \tau_\ve\big(\Phi_0(x) * e^{-z\rho}\big)\big)
 \sum_{c \in \Omega_{x, N}(X)}p(c)\big( \mathcal{B}_N(x, z, c)\big)_{i*}^{(k)},\\
\mathcal{I}_2(\ve, N) &=-\Big(\sum_{\substack{(i_1, k_1)  \geq (i_2, k_2) \\ k_1+k_2 \geq 3}}
\frac{\ve^{k_1+k_2-2}}{2N}X_{i_1*}^{(k_1)}X_{i_2*}^{(k_2)}+
\sum_{\substack{(i_2, k_2)  > (i_1, k_1) \\ k_1+k_2 \geq 3}}
\frac{\ve^{k_1+k_2-2}}{2N} X_{i_2*}^{(k_2)}X_{i_1*}^{(k_1)} \Big)\\
&\hspace{1cm}\times f\big( \tau_\ve\big(\Phi_0(x) * e^{-z\rho}\big)\big) \sum_{c \in \Omega_{x, N}(X)}p(c)
\big( \mathcal{B}_N(x, z, c)\big)_{i_1*}^{(k_1)}
\big( \mathcal{B}_N(x, z, c)\big)_{i_2*}^{(k_2)}, \\
\mathcal{I}_3(\ve, N) &=-\sum_{(i_1, k_1), (i_2, k_2), (i_3, k_3)}
\frac{\ve^{k_1+k_2+k_3-2}}{6N}
\frac{\del^3 f}{\del g_{i_1*}^{(k_1)}\del g_{i_2*}^{(k_2)}\del g_{i_3*}^{(k_3)}}(\theta)\\
&\hspace{1cm} \times  \sum_{c \in \Omega_{x, N}(X)}p(c)
\big( \mathcal{B}_N(x, z, c)\big)_{i_1*}^{(k_1)}
\big( \mathcal{B}_N(x, z, c)\big)_{i_2*}^{(k_2)}
\big( \mathcal{B}_N(x, z, c)\big)_{i_3*}^{(k_3)}.
\end{aligned}
$$
To complete the proof of Lemma \ref{key-lem}, 
it is sufficient to show the following items:

\noindent
{\normalfont (1)} $\mathrm{Ord}_\ve(-1)=0$. 

\noindent
{\normalfont (2)} We have
\begin{equation}\label{eq8-1}
\mathrm{Ord}_\ve(0) = -\A f\big( \tau_\ve\big(\Phi_0(x) * e^{-z\rho}\big)\big) + O\Big(\frac{1}{N}\Big).
\end{equation}
\noindent
{\normalfont (3)} As $N \to \infty, \ve \searrow 0$ and $N^2\ve \searrow 0$, we have
\begin{equation}\label{eq8-2}
 \|\mathcal{I}_i(\ve, N)\|_{\infty, q} \LA 0 \qquad (i=1, 2, 3).
\end{equation}

\noindent
{\bf Step~2}. We here show (1). 
We fix $i=1, 2, \dots, d_1$.  
By recalling  (\ref{m-harmonicity}) and (\ref{rel-1}), 
we have inductively 
\begin{align}\label{}
&\sum_{c \in \Omega_{x, N}(X)}p(c)\big(\mathcal{B}_N(x, z, c)\big)_{i*}^{(1)}\nn\\
&=\sum_{c' \in \Omega_{x, N-1}(X)}p(c') 
\sum_{e \in E_{t(c')}}p(e)\Big\{\log\Big(\Phi_0(x)^{-1}\cdot 
\Phi_0\big(t(c')\big)\cdot e^{-(N-1)\rho}\Big)\Big|_{X_{i}^{(1)}}\nn\\
&\hspace{1cm}+\log\Big( \Phi_0\big(o(e)\big)^{-1} \cdot \Phi_0\big(t(e)\big) \cdot e^{-\rho}\Big)
\Big|_{X_{i}^{(1)}}\Big\}\nn\\
&=\sum_{c' \in \Omega_{x, N-1}(X)}p(c') 
\log\Big(\Phi_0(x)^{-1}\cdot 
\Phi_0\big(t(c')\big)\cdot e^{-(N-1)\rho}\Big)\Big|_{X_{i}^{(1)}}
=0 \qquad (x \in V, \, z \in \mathbb{Z}). \nn
\end{align}

\noindent
{\bf Step~3}. We prove the item (2). 
First consider the coefficient of 
$X_{i_*}^{(2)}f\big( \tau_\ve\big(\Phi_0(x) * e^{-z\rho}\big)\big)$
which is given by
$$
\begin{aligned}
&-\frac{1}{N}\sum_{c \in \Omega_{x, N}(X)}p(c)
\Big\{ \big( \mathcal{B}_N(x, z, c)\big)_{i*}^{(2)}\nn\\
&\qquad  - \frac{1}{2}\sum_{1 \leq \lambda < \nu \leq d_1}
\big( \mathcal{B}_N(x, z, c)\big)_{\lambda*}^{(1)}
\big( \mathcal{B}_N(x, z, c)\big)_{\nu*}^{(1)}
 [\![X_\lambda^{(1)}, X_\nu^{(1)}]\!]\big|_{X_i^{(2)}} \Big\}\\
&=-\frac{1}{N}\sum_{c \in \Omega_{x, N}(X)}p(c) 
\Log \big( \mathcal{B}_N(x, z, c)\big)\big|_{X_{i}^{(2)}}
\qquad (x \in V, \, i=1, 2, \dots, d_2).
\end{aligned}
$$
Let us fix $i=1, 2, \dots, d_2$. 
We then deduce from (\ref{m-harmonicity}) and (\ref{rel-1}) that, 
for $x \in V$ and $z \in \mathbb{Z}$, 
$$
\begin{aligned}
&-\frac{1}{N}\sum_{c \in \Omega_{x, N}(X)}p(c) 
\log \big( \mathcal{B}_N(x, z, c)\big)\big|_{X_i^{(2)}}\\
&=-\frac{1}{N}\sum_{c' \in \Omega_{x, N-1}(X)}p(c')\sum_{e \in E_{t(c')}}p(e) 
\log\Big( \big(e^{z\rho}*\Phi_0(x)^{-1}*\Phi_0\big(t(c')\big)*e^{-(z+N-1)\rho} \big)\\
&\hspace{1cm} *\big(e^{(z+N-1)\rho}*\Phi_0\big(o(e)\big)^{-1}*\Phi_0\big(t(e)\big)* 
e^{-(z+N)\rho}\big)\Big)\Big|_{X_i^{(2)}}\\
&=-\frac{1}{N}\sum_{c' \in \Omega_{x, N-1}(X)}p(c') 
\log \Big( e^{z\rho} \cdot \Phi_0(x)^{-1} \cdot \Phi_0(t(c')) 
\cdot e^{-(z+N-1)\rho}\Big)\Big|_{X_i^{(2)}}\\
&\hspace{1cm}+\sum_{c' \in \Omega_{x, N-1}(X)}p(c') \sum_{e \in E_{t(c')}}p(e) 
\log \Big(e^{(z+N-1)\rho} \cdot d\Phi_0(e) \cdot e^{-(z+N)\rho}\Big)\Big|_{X_i^{(2)}}\\
&=-\frac{1}{N}\sum_{k=0}^{N-1}\sum_{c \in \Omega_{x, k}(X)}p(c)
\sum_{e \in E_{t(c)}}p(e)
\log \Big(e^{(z+k)\rho} \cdot d\Phi_0(e) \cdot e^{-(z+k+1)\rho}\Big)\Big|_{X_i^{(2)}}.
\end{aligned}
$$
For $g, h \in G$, we denote by $[g, h]:=g \cdot h \cdot g^{-1} \cdot h^{-1}$ 
the usual commutator of $g$ and $h$.
Then we have
$$
\begin{aligned}
&\sum_{e \in E_{t(c)}}p(e)
\log \Big(e^{(z+k)\rho} \cdot d\Phi_0(e)  \cdot  e^{-(z+k+1)\rho}\Big)\Big|_{X_i^{(2)}}\\
&=\sum_{e \in E_{t(c)}}p(e)
\log \Big(\big[e^{(z+k)\rho}, d\Phi_0(e)\big]  \cdot  d\Phi_0(e) \cdot e^{-\rho}\Big)\Big|_{X_i^{(2)}}\\
&=\sum_{e \in E_{t(c)}}p(e)
\log\Big(\big[e^{(z+k)\rho}, d\Phi_0(e)\big]\Big)\Big|_{X_i^{(2)}} 
+ \sum_{e \in E_{t(c)}}p(e)\log \big(d\Phi_0(e) \cdot e^{-\rho}\big)\big|_{X_i^{(2)}}\\
&=\sum_{e \in E_{t(c)}}p(e) \log \big( d\Phi_0(e) \cdot  e^{-\rho}\big)\big|_{X_i^{(2)}}
\qquad (z \in \mathbb{Z}, \, k=0, 1, \dots, N-1)
\end{aligned}
$$
by again using (\ref{m-harmonicity}).
Since the function 
$$
M_{i}(x):=\sum_{e \in E_x}p(e) 
\Log \big(d\Phi_0(e) \cdot e^{-\rho}\big)\big|_{X_{i}^{(2)}} 
\qquad (i=1, 2, \dots, d_2,  \, x \in V)
$$
satisfies $M_i(\gamma x)=M_i(x)$ for $\gamma \in \Gamma$ and $x \in V$
due to the $\Gamma$-invariance of $p$ and the $\Gamma$-equivariance of $\Phi_0$, 
there exists a function $\mathcal{M}_i : V_0 \LA \mathbb{R}$ such that 
$\mathcal{M}_i\big(\pi(x)\big)=M_i(x)$ for $i=1, 2, \dots, d_2$ and $x \in V$. 
Moreover, we have
$$
L^k \mathcal{M}_i\big(\pi(x)\big)
=L^k M_i(x)
\qquad (k \in \mathbb{N}, \, i=1, 2, \dots, d_2, \, x \in V)
$$
by using the $\Gamma$-invariance of $p$. 
Then the ergodic theorem 
(cf.~\cite[Theorem 3.2]{IKK}) for the transition operator $L$ gives
\begin{align}
&-\frac{1}{N}\sum_{c \in \Omega_{x, N}(X)}p(c) 
\log \big( \mathcal{B}_N(x, z, c)\big)\big|_{X_i^{(2)}}\nn\\
&=-\frac{1}{N}\sum_{k=0}^{N-1} L^k M_i (x)\nn\\
&=-\frac{1}{N}\sum_{k=0}^{N-1} L^k \mathcal{M}_i\big(\pi(x)\big)\nn\\
&=-\sum_{x \in V_0}m(x)\mathcal{M}_i(x)+O\Big(\frac{1}{N}\Big)
=-\beta(\Phi_0)\big|_{X_i^{(2)}} + O\Big(\frac{1}{N}\Big) 
\quad (x \in V, \, z \in \mathbb{Z}). \label{eq5}
\end{align}

We next consider the coefficient of 
$X_{i*}^{(1)}X_{j*}^{(1)}f\big( \tau_\ve\big(\Phi_0(x) * e^{-z\rho}\big)\big)$
which is given by
$$
-\frac{1}{2N}\sum_{c \in \Omega_{x, N}(X)}p(c)
\big( \mathcal{B}_N(x, z, c)\big)_{i*}^{(1)}\big( \mathcal{B}_N(x, z, c)\big)_{j*}^{(1)}
\qquad (x \in V, \, z \in \mathbb{Z}, \, i, j=1, 2, \dots, d_1).
$$
Fix $i, j=1, 2, \dots, d_1$. 
Then (\ref{m-harmonicity}) and (\ref{rel-1}) imply
$$
\begin{aligned}
&-\frac{1}{2N}\sum_{c \in \Omega_{x, N}(X)}p(c)
\big( \mathcal{B}_N(x, z, c)\big)_{i*}^{(1)}\big( \mathcal{B}_N(x, z, c)\big)_{j*}^{(1)}\\
&=-\frac{1}{2N}\sum_{c' \in \Omega_{x, N-1(X)}}p(c')\sum_{e \in E_{t(c')}}p(e)\\
&\hspace{1cm}\times \Big\{ \log\big( \mathcal{B}_{N-1}(x, z, c')\big)\big|_{X_{i}^{(1)}}
+\log\big(d\Phi_0(e)\cdot e^{-\rho}\big)\big|_{X_{i}^{(1)}}\Big\}\\
 &\hspace{1cm}\times \Big\{ \log\big( \mathcal{B}_{N-1}(x, z, c')\big)\big|_{X_{j}^{(1)}}
+\log\big(d\Phi_0(e) \cdot e^{-\rho}\big)\big|_{X_{j}^{(1)}}\Big\}\\
&=-\frac{1}{2N}\Big\{\sum_{c' \in \Omega_{x, N-1(X)}}p(c') 
\log\big( \mathcal{B}_{N-1}(x, z, c')\big)\big|_{X_{i}^{(1)}}
\log\big( \mathcal{B}_{N-1}(x, z, c')\big)\big|_{X_{j}^{(1)}}\\
&\hspace{1cm}+\sum_{e \in E_{t(c')}}p(e) 
\log\big(d\Phi_0(e)\cdot e^{-\rho}\big)\big|_{X_{i}^{(1)}}
\log\big(d\Phi_0(e)\cdot e^{-\rho}\big)\big|_{X_{j}^{(1)}}\Big\}\\
&=-\frac{1}{2N}\sum_{k=0}^{N-1}\sum_{c \in \Omega_{x, N}(X)}p(c)
\sum_{e \in E_{t(c)}}p(e)
\log\big(d\Phi_0(e)\cdot e^{-\rho}\big)\big|_{X_{i}^{(1)}}
\log\big(d\Phi_0(e)\cdot e^{-\rho}\big)\big|_{X_{j}^{(1)}}
\end{aligned}
$$
for $x \in V$ and $z \in \mathbb{Z}$.
In the same argument as above, the function $N_{ij} : V \LA \mathbb{R}$ defined by
$$
N_{ij}(x):=\sum_{e \in E_x}p(e)
\Log \big(d\Phi_0(e) \cdot e^{-\rho}\big)\big|_{X_i^{(1)}} 
\Log \big(d\Phi_0(e) \cdot e^{-\rho}\big)\big|_{X_j^{(1)}} \quad (i, j=1, 2, \dots, d_1, \, x \in V)
$$
is $\Gamma$-invariant and then there exists a function
$\mathcal{N}_{ij} : V_0 \LA \mathbb{R}$ such that 
$\mathcal{N}_{ij}\big(\pi(x)\big)=N_{ij}(x)$ for $x \in V$. 
We also have $$
L^k \mathcal{N}_{ij}\big(\pi(x)\big)
=L^k N_{ij}(x)
\qquad (k \in \mathbb{N}, \, i, j=1, 2, \dots, d_2, \, x \in V)
$$
by using the $\Gamma$-invariance of $p$. 
Hence, we obtain
\begin{align}
&-\frac{1}{2N}\sum_{c \in \Omega_{x, N}(X)}p(c)
\big( \mathcal{B}_N(x, z, c)\big)_{i*}^{(1)}\big( \mathcal{B}_N(x, z, c)\big)_{j*}^{(1)}\nn\\
&=-\frac{1}{2N}\sum_{k=0}^{N-1} L^k N_{ij}(x) \nn\\
&=-\frac{1}{2N}\sum_{k=0}^{N-1} L^k \mathcal{N}_{ij}\big(\pi(x)\big) \nn\\
&=-\frac{1}{2}\sum_{x \in V_0}
m(x)\mathcal{N}_{ij}(x)+O\Big(\frac{1}{N}\Big)\nn\\
&=-\frac{1}{2}\sum_{e \in E_0}\widetilde{m}(e)
\log\big(d\Phi_0(\widetilde{e})\cdot e^{-\rho}\big)\big|_{X_{i}^{(1)}}
\log\big(d\Phi_0(\widetilde{e})\cdot e^{-\rho}\big)\big|_{X_{j}^{(1)}}+O\Big(\frac{1}{N}\Big).\label{eq6}
\end{align}
by virtue of the ergodic theorem.
Recall that $\{V_1, V_2, \dots, V_{d_1}\}$ denotes an orthonormal basis of $(\g^{(1)}, g_0)$. 
We especially put $X_i^{(1)}=V_i$ for $i=1, 2, \dots, d_1$. 
Let $\{\omega_1, \omega_2, \dots, \omega_{d_1}\} 
\subset \Hom(\g^{(1)}, \mathbb{R}) \hookrightarrow \h^1(X_0, \mathbb{R})$
be the dual basis of $\{V_1, V_2, \dots, V_{d_1}\}$. 
Namely, $\omega_i(V_j)=\delta_{ij}$ for $i, j=1, 2, \dots, d_1$. 
It follows from (\ref{Alb map}) that
\begin{align}
&\sum_{e \in E_0}\widetilde{m}(e)
\Log\big(d\Phi_0(\widetilde{e}) \cdot e^{-\rho}\big)\big|_{V_i}
\Log\big(d\Phi_0(\widetilde{e}) \cdot e^{-\rho}\big)\big|_{V_j}\nn\\
&=\sum_{e \in E_0}\widetilde{m}(e)
\Log\big(d\Phi_0(\widetilde{e})\big)\big|_{V_i}\Log\big(d\Phi_0(\widetilde{e})\big)\big|_{V_j}
-\rho_{\mathbb{R}}(\gamma_p)\big|_{V_i}\rho_{\mathbb{R}}(\gamma_p)\big|_{V_j}\nn\\
&=\sum_{e \in E_0}\widetilde{m}(e)
{}^t\rho_{\mathbb{R}}(\omega_i)(e){}^t\rho_{\mathbb{R}}(\omega_j)(e)
-\omega_i(\rho_{\mathbb{R}}(\gamma_p))\omega_j(\rho_{\mathbb{R}}(\gamma_p))\nn\\
&=\sum_{e \in E_0}\widetilde{m}(e)\omega_i(e)\omega_j(e) 
- \la \gamma_p, \omega_i \ra \la \gamma_p, \omega_j \ra =\La \omega_i, \omega_j \Ra_p = \delta_{ij}.
\label{eq7} 
\end{align}
Hence, we obtain (\ref{eq8-1}) by combining (\ref{eq5}) with (\ref{eq6}) and (\ref{eq7}).

\vspace{2mm}
\noindent
{\bf Step~4}. 
We show (3) at the last step.  
We first discuss the estimate of $\mathcal{I}_1(\ve, N)$. 
By using (\ref{rel-5}) and (\ref{N^k}), we have
$$
\begin{aligned}
&\Big|\Big(  \Phi_0(x)^{-1}*\Phi_0\big(t(c)\big)\Big)_{i*}^{(k)}\Big|\nn\\
&\leq C\sum_{\substack{|K_1|+|K_2| \leq k\\ |K_2| > 0}}
\Big|\mathcal{P}_*^{K_1}\Big(  \Phi_0(x)^{-1}\Big)\Big|
\Big|\mathcal{P}^{K_2}\Big(  \Phi_0(x)^{-1}\cdot \Phi_0\big(t(c)\big)\Big)\Big|\nn\\
&\leq C\sum_{\substack{|K_1|+|K_2| \leq k \\ |K_2| > 0}}N^{|K_2|}
\Big|\mathcal{P}_*^{K_1}\Big( e^{-z\rho}*\big(\Phi_0(x) * e^{-z\rho}\big)^{-1}\Big)\Big|\nn\\
\end{aligned}
$$
for $i=1, 2, \dots, d_k$ and $k=1, 2, \dots, r$.
Then (\ref{CBH-formula}) implies that there is a continuous function 
$Q_1 : G \LA \mathbb{R}$ such that
\begin{align}
&\Big|\Big(  \Phi_0(x)^{-1}*\Phi_0\big(t(c)\big)\Big)_{i*}^{(k)}\Big|\nn\\
&\leq |z|^{k-1}Q_1\big( \tau_\ve\big(\Phi_0(x)*e^{-z\rho}\big)\big)
\sum_{\substack{|K_1|+|K_2| \leq k \\ |K_2| > 0}}\ve^{-|K_1|}N^{|K_2|}
\label{I_1-est-1}
\end{align}
for $i=1, 2, \dots, d_k$ and $k=1, 2, \dots, r$.
Thus, (\ref{CBH-formula}) and (\ref{I_1-est-1}) yields
\begin{align}
&\big|\big( \mathcal{B}_N(x, 0, c)\big)_{i*}^{(k)}\big|\nn\\
&\leq C\sum_{\substack{|L_1|+|L_2|=k \\ |L_1|, |L_2| \geq 0}} 
\Big|\mathcal{P}_*^{L_1}\Big(  \Phi_0(x)^{-1}*\Phi_0\big(t(c)\big)\Big)\Big|
\big|\mathcal{P}_*^{L_2}\big( e^{-N\rho}\big)\big|\nn\\
&\leq C|z|^kQ_2\big( \tau_\ve\big(\Phi_0(x)*e^{-z\rho}\big)\big)
\sum_{\substack{|L_1|+|L_2|=k \\ |L_1|, |L_2| \geq 0}} 
N^{|L_2|}\sum_{\substack{|K_1|+|K_2| \leq |L_1| \\ |K_2| > 0}}\ve^{-|K_1|}N^{|K_2|}\nn\\
&=C|z|^kQ_2\big( \tau_\ve\big(\Phi_0(x)*e^{-z\rho}\big)\big) F(\ve, N)
\label{I_1-est-2}
\end{align}
for some continuous function $Q_2 : G \LA \mathbb{R}$, 
where $F(\ve, N)$ denotes the polynomial of $\ve$ and $N$ which satisfies 
$\ve^{k-2}N^{-1}F(\ve, N) \to 0$ 
as $N \to \infty, \ve \searrow 0$ and $N^2\ve \searrow 0$. 

On the other hand, combining (\ref{I_1-est-2}) with  
$\rho_{\mathbb{R}}(\gamma_p) \in \g^{(1)}$ gives 
\begin{align}
&\frac{\ve^{k-2}}{N}\big|\big(\mathcal{B}_N(x, z, c)\big)_{i*}^{(k)}\big|\nn\\
&=\frac{\ve^{k-2}}{N}\Big| \Big( \big[ e^{z\rho},\mathcal{B}_N(x, 0, c)\big]_* 
* \mathcal{B}_N(x, 0, c)\Big)_{i*}^{(k)}\Big| \nn\\
&\leq C\frac{\ve^{k-2}}{N}\sum_{\substack{|K_1|+|K_2|=k \\ |K_1|, |K_2| \geq 0}}
\Big|\mathcal{P}_*^{K_1}\Big( \big[ e^{z\rho}, \mathcal{B}_N(x, 0, c)\big]_* \Big)\Big|
\Big|\mathcal{P}_*^{K_2}\big(\mathcal{B}_N(x, 0, c)\big)\Big| \nn\\
 &\leq C|z|^{2k}\frac{\ve^{k-2}}{N} Q_3\big( \tau_\ve\big(\Phi_0(x)*e^{-z\rho}\big)\big)
F(\ve, N)\nn\\
&\hspace{3cm}\big(i=1, 2, \dots, d_k, \, k=3, 4, \dots, r, \, x \in V, \, 
z \in \mathbb{Z}, \, c \in \Omega_{x, N}(X)\big)
\label{I_1-est-3}
\end{align}
for some continuous function $Q_3 : G \LA \mathbb{R}$. 
Hence, we obtain $\|\mathcal{I}_1(\ve, N)\|_{\infty, q} \to 0$
as $N \to \infty, \ve \searrow 0$ and $N^2\ve \searrow 0$ in 
$C_{\infty, q}(X \times \mathbb{Z})$
by using (\ref{I_1-est-3}).
This follows from $2k < 2r < q$. 
In the same argument as above, we also obtain
$\|\mathcal{I}_2(\ve, N)\|_{\infty, q} \to 0$ as $N \to \infty, \, \ve \searrow 0$ 
and $N^2\ve \searrow  0$ in $C_{\infty, q}(X \times \mathbb{Z})$-topology
since the order of $|z|$ in $\mathcal{I}_2(\ve, N)$ satisfies $2 \times 2k < 4r <q$. 

Finally, we study the estimate of the term $\mathcal{I}_3(\ve, N)$. 
We recall that $f \in C_0^\infty(G)$ and
$\supp \del^3 f/(\del g_{i_1*}^{(k_1)}\del g_{i_2*}^{(k_2)}\del g_{i_3*}^{(k_3)}) \subset \supp f.$
Therefore, it suffices to show by induction on $k=1, 2, \dots, r$ that, 
if $\ve N<1$,  
\begin{equation}\label{eq-10}
\ve^k \big|\big( \mathcal{B}_N(x, z, c)\big)_{i*}^{(k)}\big|
\leq |z|^kQ^{(k)}\big(\tau_\ve(\Phi_0(x) * e^{-z\rho})*\theta\big)
\times \ve N
\end{equation}
for some continuous function $Q^{(k)} : G \LA \mathbb{R}$, 
where $\theta \in G$ appears in the remainder term of (\ref{eq1}). 
The cases $k=1$ and $k=2$ are obvious. 
Suppose that (\ref{eq-10}) holds for less than $k$. 
Then we have 
$$
\ve^k \big|\big( \mathcal{B}_N(x, z, c)\big)_{i*}^{(k)}\big|
\leq C\ve^k \sum_{\substack{|K_1|+|K_2| \leq k \\ |K_2|>0}} 
\Big|\mathcal{P}_*^{K_1}\Big(  \Phi_0(x)^{-1}\Big)\Big|
\Big|\mathcal{P}^{K_2}\Big(  \Phi_0(x)^{-1}\cdot \Phi_0\big(t(c)\big)\Big)\Big|
$$ 
by using (\ref{rel-5}). Since
$$
\big( \Phi_0(x)^{-1}\big)_{i_1*}^{(k_1)}
=\Big( e^{-z\rho}*(\tau_{\ve^{-1}}\theta) 
* \big( \tau_{\ve^{-1}}(\tau_\ve(\Phi_0(x) * e^{-z\rho})*\theta)^{-1}\big)\Big)_{i_1*}^{(k_1)}
\qquad (k_1 \leq k-1), 
$$
we have inductively
$$
\big|\big( \Phi_0(x)^{-1}\big)_{i_1*}^{(k_1)}\big|
\leq |z|^{k_1}Q\big(\tau_\ve(\Phi_0(x) * e^{-z\rho})*\theta\big)
$$
for a continuous function $Q : G \LA \mathbb{R}$ and $k_1 \leq k-1$. 
We thus obtain
$$
\begin{aligned}
&\ve^k \big|\big( \mathcal{B}_N(x, z, c)\big)_{i*}^{(k)}\big|\nn\\
&\leq C\ve^k \sum_{\substack{|K_1|+|K_2| \leq k \\ |K_2| > 0}} N^{|K_2|}
\Big| \mathcal{P}_*^{K_1}\Big( e^{-z\rho}*(\tau_{\ve^{-1}}\theta) 
* \big( \tau_{\ve^{-1}}(\tau_\ve(\Phi_0(x) * e^{-z\rho})*\theta)^{-1}\big)\Big)\Big|\\
&\leq C|z|^kQ\big(\tau_\ve(\Phi_0(x) * e^{-z\rho})*\theta\big)
 \sum_{\substack{|K_1|+|K_2| \leq k \\ |K_2| > 0}}
 \ve^{k-|K_1|+1}N^{|K_2|+1}\\
 &\leq |z|^k Q^{(k)}\big(\tau_\ve(\Phi_0(x) * e^{-z\rho})*\theta\big)
\times \ve N
\end{aligned}
$$ 
for some continuous function $Q^{(k)} : G \LA \mathbb{R}$.
Therefore, (\ref{eq-10}) holds for $k=1, 2, \dots, r$ and 
this implies that $\|\mathcal{I}_3(\ve, N)\|_{\infty, q} \to 0$ as $N \to \infty, \, \ve \searrow 0$ 
and $N^2\ve \searrow  0$ in $C_{\infty, q}(X \times \mathbb{Z})$
since the order of $|z|$ in $\mathcal{I}_3(\ve, N)$ satisfies $3 k < 3r <q$. 
This completes the proof. \qed 

\vspace{2mm}
We now give the proof of Theorem \ref{CLT1} by using this lemma. 
We note that the infinitesimal operator $\A$ in Lemma \ref{key-lem} enjoys the following property.

\begin{lm}\label{operator-lem}
{\bf (cf.~Robinson \cite[page 304]{Rob})} 
The range of $\lambda -\A$ is dense in $C_\infty(G)$ for some $\lambda>0$.
Namely, $(\lambda -\A)\big(C_0^\infty(G)\big)$ is dense in $C_\infty(G)$. 
\end{lm}

\vspace{2mm}
\noindent
{\bf Proof of Theorem \ref{CLT1}}. 
(1) We follow the argument in Kotani \cite[Theorem 4]{Kotani}. 
Let $N=N(n)$ be the integer satisfying $n^{1/5} \leq N<n^{1/5}+1$
and $k_N$ and $r_N$ 
be the quotient and the remainder of $([nt]-[ns])/N(n)$, respectively.
Note that  $r_N < N$. 
We put $\ve_N:=n^{-1/2}$ and $h_N:=N\ve_N^2$. 
Then we have $N=N(n) \to \infty$,
$$
r_N^2\ve_N<N^2\ve_N \leq (1+n^{1/5})^2 \cdot n^{-1/2} \to 0,
$$
and $h_N \leq (1+n^{1/5}) \cdot n^{-1}\to 0$ as $n \to \infty$.
We also see that 
$$
r_N\ve_N^2 < N\ve_N^2 \leq (1+n^{1/5}) \cdot n^{-1} \to 0 \quad (n \to \infty).
$$
Hence, we have
$$
k_N h_N=\frac{[nt]-[ns]-r_N}{N} \cdot N\ve_N^2=\big([nt]-[ns]-r_N\big)\ve_N^2 \to t-s \quad (n \to \infty).
$$

Since $C_0^\infty(G) \subset \Dom(\A) \subset C_\infty(G)$ and 
$C_0^\infty(G)$ is dense in $C_\infty(G)$, the operator $\A$ is 
densely defined in $C_\infty(G)$. 
We use this fact and Lemma \ref{operator-lem} to  
apply Trotter's approximation theorem (cf.~Trotter \cite{Trotter} and Kurtz \cite{Kurtz}). 
We obtain, for $f \in C_0^\infty(G)$,
\begin{equation}\label{equ}
\lim_{n \to \infty}\Big\| \mathcal{L}_p^{Nk_N}\mathcal{P}_{n^{-1/2}}f - 
\mathcal{P}_{n^{-1/2}}\e^{-(t-s)\A}f\Big\|_{\infty, q}=0.
\end{equation}
Then Lemma \ref{key-lem} implies
\begin{equation}\label{eq8}
\lim_{n \to \infty}\Big\| \frac{1}{r_N\ve_N^2}\big(I-\mathcal{L}_p^{r_N}\big)\mathcal{P}_{n^{-1/2}}f 
- \mathcal{P}_{n^{-1/2}}\A f \Big\|_{\infty, q}=0
\end{equation}
for all $f \in C_0^\infty(G)$. We thus have
\begin{align}\label{eq9}
&\Big\| \mathcal{L}_p^{[nt]-[ns]}P_{n^{-1/2}}f - P_{n^{-1/2}}\e^{-(t-s)\A}f\Big\|_{\infty, q} \nn\\
&\leq \Big\| \big(I-\mathcal{L}_p^{r_N}\big)\mathcal{P}_{n^{-1/2}}f\Big\|_{\infty, q} + 
\Big\| \mathcal{L}_p^{Nk_N}\mathcal{P}_{n^{-1/2}}f 
- \mathcal{P}_{n^{-1/2}}\e^{-(t-s)\A}f\Big\|_{\infty, q}.
\end{align}
On the other hand, we have
\begin{align}\label{eq10}
&\Big\| \big(I - \mathcal{L}_p^{r_N}\big)\mathcal{P}_{n^{-1/2}}f\Big\|_{\infty, q} \nn\\
&\leq r_N\ve_N^2 \Big\|\frac{1}{r_N\ve_N^2}\big(I-\mathcal{L}_p^{r_N}\big)\mathcal{P}_{n^{-1/2}}f
-\mathcal{P}_{n^{-1/2}}\A f\Big\|_{\infty, q} + r_N\ve_N^2\big\| \mathcal{P}_{n^{-1/2}}\A f\big\|_{\infty, q} \nn \\
&\leq r_N\ve_N^2 \Big\|\frac{1}{r_N\ve_N^2}\big(I-\mathcal{L}_p^{r_N}\big)\mathcal{P}_{n^{-1/2}}f
-\mathcal{P}_{n^{-1/2}}\A f\Big\|_{\infty, q} + r_N\ve_N^2\big\| \A f\big\|_\infty^G.
\end{align}
We obtain (\ref{semigroup CLT1}) for $f \in C_0^\infty(G)$
by combining (\ref{eq8}), (\ref{eq9}) and (\ref{eq10}) 
with $r_N\ve_N^2 \to 0 \, (n \to \infty)$.
For $f \in C_\infty(G)$, we also obtain the convergence (\ref{semigroup CLT1})
by following the same argument as \cite[Theorem 2.1]{IKK}. 

\vspace{2mm}
\noindent
(2) For $t>0$ and $z \in \mathbb{Z}$, we have
$$
\begin{aligned}
&\big| \mathcal{L}_p^{[nt]}\mathcal{P}_{n^{-1/2}}f(x_n, z) - \e^{-t\A}f(g)\big|\\
&\leq \big| \mathcal{L}_p^{[nt]}\mathcal{P}_{n^{-1/2}}f(x_n, z) 
- \mathcal{P}_{n^{-1/2}}\e^{-t \A}f(x_n, z)\big|
+\big| \mathcal{P}_{n^{-1/2}}\e^{-t\A}f(x_n, z) - e^{-t\A}f(g)\big|\\
&\leq (1+|z|^q)\Big\| \mathcal{L}_p^{[nt]}\mathcal{P}_{n^{-1/2}}f 
- \mathcal{P}_{n^{-1/2}}\e^{-t\A}f\Big\|_{\infty, q} \\
&\hspace{1cm}
+ \Big| \e^{-t\A}f\Big(\tau_{n^{-1/2}}\big(\Phi_0(x_n) 
* \exp(-z \rho_{\mathbb{R}}(\gamma_p))\big)\Big) - \e^{-t\A}f(g)\Big|.
\end{aligned}
$$
We thus obtain (\ref{semigroup CLT1-2}) by (\ref{semigroup CLT1}) 
and the continuity of the function 
$\e^{-t \A}f : G \LA \mathbb{R}$.
This completes the proof of Theorem \ref{CLT1}. \qed

\vspace{2mm}
We now give several properties of $\beta(\Phi_0)$.

\begin{pr}\label{property-beta}
{\normalfont (1)} If the random walk on $X$ is $m$-symmetric, then $\beta(\Phi_0)=\bm{0}_{\g}$. 

\noindent
{\normalfont (2)}
Let $\Phi_0, \, \widehat{\Phi}_0 : X \LA G$ be two modified harmonic realizations. 
Then
$$
\beta(\Phi_0)=\beta(\widehat{\Phi}_0) - 
\big[\rho_{\mathbb{R}}(\gamma_p),
 \log\big(\Phi_0(x)^{-1} \cdot \widehat{\Phi}_0(x)\big)\big]\big|_{\g^{(2)}} \qquad (x \in V).
$$
In particular, if either 

$\bullet$ $\log \Phi_0(x_*)\big|_{\g^{(1)}}
=\log \widehat{\Phi}_0(x_*)\big|_{\g^{(1)}}$ for some reference point $x_* \in V$, or

$\bullet$ $\rho_{\mathbb{R}}(\gamma_p)=\bm{0}_{\g}$  \\
holds, then we have $\beta(\Phi_0)=\beta(\widehat{\Phi}_0)$.
\end{pr}

\noindent
{\bf Proof.}
Assertion (1) is easily obtained as follows:
$$
\begin{aligned}
\beta(\Phi_0)&=\frac{1}{2}\sum_{e \in E_0}
\Big\{ \widetilde{m}(e)\Log\big(d\Phi_0(\widetilde{e})\big)\big|_{\g^{(2)}}
+\widetilde{m}(\ol{e})\Log\big(d\Phi_0(\widetilde{\ol{e}})\big)\big|_{\g^{(2)}}\Big\}\\
&=\frac{1}{2}\sum_{e \in E_0}\big( \widetilde{m}(e)-\widetilde{m}(\ol{e})\big)
\Log\big(d\Phi_0(\widetilde{e})\big)\big|_{\g^{(2)}}=\bm{0}_{\g}.
\end{aligned}
$$
Next we show Assertion (2). 
We set $\Psi(x):=\Phi_0(x)^{-1} \cdot \widehat{\Phi}_0(x)$ for $x \in V.$ 
We note that the map $\Psi : X \LA G$ is $\Gamma$-invariant. 
Since the $\g^{(1)}$-components of $\Phi_0$ and $\widehat{\Phi}_0$ are 
uniquely determined up to $\g^{(1)}$-translation, there exists a constant vector $C \in \g^{(1)}$
such that $\log\big(\Psi(x)\big)\big|_{\g^{(1)}}=C$ for $x \in V.$
Define a function $F_i : V \LA \mathbb{R}$ by 
$F_i(x):=\log \big( \Psi(x)\big)\big|_{X_i^{(2)}}$ for $i=1, 2, \dots, d_2$ and $x \in V$. 
Then we see that the function $F_i$ is $\Gamma$-invariant. 
Hence, there is a function $\widehat{F}_i : V_0 \LA \mathbb{R}$ satisfying
$\widehat{F}_i\big(\pi(x)\big)=F_i(x)$ for $x \in V$. 
Then we obtain
$$
\begin{aligned}
\beta(\Phi_0)
&=\sum_{e \in E_0}\widetilde{m}(e) \log \Big( 
         \Psi\big(o(\widetilde{e})\big) \cdot 
         \big(d\widehat{\Phi}_0(\widetilde{e}) \cdot e^{-\rho}\big) \cdot e^\rho
         \cdot \Psi\big(t(\widetilde{e})\big)^{-1} \cdot e^{-\rho} \Big)\Big|_{\g^{(2)}}\\
&=\beta(\widehat{\Phi}_0)-\sum_{e \in E_0}\widetilde{m}(e)\Big\{
         \log\Big( \Psi\big(t(\widetilde{e})\big)\Big)\Big|_{\g^{(2)}}-
         \log\Big( \Psi\big(o(\widetilde{e})\big)\Big)\Big|_{\g^{(2)}}\Big\}
         -[\rho_{\mathbb{R}}(\gamma_p), C]\big|_{\g^{(2)}} \\
&=\beta(\widehat{\Phi}_0)
- \sum_{i=1}^{d_2}\big({}_{C_1(X_0, \mathbb{R})}\la \gamma_p, 
d\widehat{F}_i\ra_{C^1(X_0, \mathbb{R})}\big)
X_i^{(2)} -[\rho_{\mathbb{R}}(\gamma_p), C]\big|_{\g^{(2)}} \\
&= \beta(\widehat{\Phi}_0)
- \sum_{i=1}^{d_2}\big({}_{C_0(X_0, \mathbb{R})}\la \del(\gamma_p), 
\widehat{F}_i\ra_{C^0(X_0, \mathbb{R})}\big)X_i^{(2)} -[\rho_{\mathbb{R}}(\gamma_p), C]\big|_{\g^{(2)}}\\
&=\beta(\widehat{\Phi}_0)
-[\rho_{\mathbb{R}}(\gamma_p), C]\big|_{\g^{(2)}}, 
\end{aligned}
$$
where we used (\ref{CBH-formula}) for the second line and $\gamma_p \in \h_1(X_0, \mathbb{R})$
for the fourth line. \qed

\subsection{Proof of Theorem \ref{FCLT1}}

We now assume the centered condition {\bf (C)}: 
$\rho_{\mathbb{R}}(\gamma_p)=\bm{0}_{\g}$,
throughout this subsection. 
For $k=1, 2, \dots, r$, we denote by $(G^{(k)}, \cdot)$ and $(G^{(k)}, *)$
the connected and simply connected nilpotent Lie group of step $k$ 
and the corresponding limit group
whose Lie algebras are 
$\big(\g^{(1)} \oplus \g^{(2)} \oplus \cdots \oplus \g^{(k)}, [ \cdot, \cdot]\big)$ and 
$\big(\g^{(1)} \oplus \g^{(2)} \oplus \cdots \oplus \g^{(k)}, [\![ \cdot, \cdot]\!]\big),$
respectively. 
For the piecewise smooth stochastic process 
$(\mathcal{Y}_t^{(n)})_{0 \leq t \leq 1}
=(\mathcal{Y}_t^{(n), 1}, \mathcal{Y}_t^{(n), 2}, \dots, \mathcal{Y}_t^{(n), r})_{0 \leq t \leq 1}$
defined in Section 2, we define its truncated process by
$$
\mathcal{Y}_t^{(n; \,k)}=
\big( \mathcal{Y}_t^{(n), 1}, 
 \mathcal{Y}_t^{(n), 2}, \dots,
  \mathcal{Y}_t^{(n), k} \big) \in G^{(k)}
  \qquad (0 \leq t \leq 1, \, k=1, 2, \dots, r)
$$
in the $(\cdot)$-coordinate system.
To complete the proof of Theorem \ref{FCLT1}, it is sufficient to show 
the tightness of $\{\Prob^{(n)}\}_{n=1}^\infty$ (Lemma \ref{tightness1}) and 
the convergence of the finite dimensional distribution of $\{\mathcal{Y}_\cdot^{(n)}\}_{n=1}^\infty$
(Lemma \ref{CFDD1}).

In the former part of this subsection, we aim to show the following.
\begin{lm}\label{tightness1}
Under {\bf (C)}, 
the family $\{\Prob^{(n)}\}_{n=1}^\infty$ is tight in 
$C^{0, \alpha\text{\rm{-H\"ol}}}_{\bm{1}_G}([0, 1]; G),$
where $\alpha$ is an arbitrary real number less than $1/2$. 
\end{lm}
As the first step of the proof of Lemma \ref{tightness1},
we prepare the following lemma.

\begin{lm}\label{sublemma1}
Let $m, n$ be positive integers. Then there exists a constant $C>0$ 
which is independent of $n$ 
{\rm(}however, it may depend on $m${\rm)} such that 
\begin{equation}\label{tight1-1}
\mathbb{E}^{\mathbb{P}_{x_*}}
\Big[ d_{\mathrm{CC}}(\mathcal{Y}_{s}^{(n; \, 2)}, 
\mathcal{Y}_{t}^{(n; \, 2)})^{4m}\Big]
\leq C(t-s)^{2m}  \qquad (0 \leq s \leq  t \leq 1).
\end{equation}
\end{lm}
\noindent
{\bf Proof.} 
The proof is partially based on Bayer--Friz \cite[Proposition 4.3]{BF}. 
We split the proof into several steps. 

\vspace{2mm}
\noindent
{\bf Step\,1.} 
At the beginning, we show 
\begin{equation}\label{tight-1}
\mathbb{E}^{\mathbb{P}_{x_*}}
\Big[ d_{\mathrm{CC}}(\mathcal{Y}_{t_k}^{(n; \, 2)}, \mathcal{Y}_{t_\ell}^{(n; \, 2)})^{4m}\Big]
\leq C\Big(\frac{\ell-k}{n}\Big)^{2m}  \qquad \big(n, m \in \mathbb{N}, 
\, t_k, t_\ell \in \mathcal{D}_n \, (k \leq \ell)\big)
\end{equation}
for some $C>0$ independent of $n$ (depending on $m$).
By recalling the equivalence of two homogeneous norms 
$\|\cdot\|_{\mathrm{CC}}$ and $\|\cdot\|_{\mathrm{hom}}$
(cf.~Proposition \ref{homogeneous equiv}), we readily see that 
(\ref{tight-1}) is equivalent to the existence 
of positive constants $C^{(1)}$ and $C^{(2)}$ independent of $n$ such that 
\begin{equation}\label{expectation-1}
\mathbb{E}^{\mathbb{P}_{x_*}}\Big[ \big\| \Log\big((\mathcal{Y}_{t_k}^{(n)})^{-1}
\cdot \mathcal{Y}_{t_{\ell}}^{(n)}\big)\big|_{\g^{(1)}}\big\|_{\g^{(1)}}^{4m}\Big] 
 \leq C^{(1)}\Big(\frac{\ell-k}{n}\Big)^{2m},
\end{equation}
\begin{equation}\label{expectation-2}
\mathbb{E}^{\mathbb{P}_{x_*}}\Big[ \big\| \Log\big((\mathcal{Y}_{t_k}^{(n)})^{-1}
 \cdot \mathcal{Y}_{t_{\ell}}^{(n)}\big)\big|_{\g^{(2)}}\big\|_{\g^{(2)}}^{2m}\Big]
  \leq C^{(2)}\Big(\frac{\ell-k}{n}\Big)^{2m}.
\end{equation}

\vspace{2mm}
\noindent
{\bf Step\,2.} We now show (\ref{expectation-1}). 
We see
\begin{align}\label{step2-1}
 &
 \mathbb{E}^{\mathbb{P}_{x_*}}\Big[ \big\| \Log\big((\mathcal{Y}_{t_k}^{(n)})^{-1} \cdot
 \mathcal{Y}_{t_{\ell}}^{(n)}\big)\big|_{\g^{(1)}}\big\|_{\g^{(1)}}^{4m}\Big]\nn\\
&=  \Big( \frac{1}{\sqrt{n}}\Big)^{4m}  \mathbb{E}^{\mathbb{P}_{x_*}}
\Big[ \Big( \sum_{i=1}^{d_1} \Log(\xi_k^{-1} \cdot \xi_\ell) \big|_{X_i^{(1)}}^2\Big)^{2m}\Big]\nn\\
&\leq \Big( \frac{1}{\sqrt{n}}\Big)^{4m} \cdot d_1^{2m} \max_{i=1, 2, \dots, d_1} \max_{x \in \mathcal{F}} 
\Big\{\sum_{c \in \Omega_{x, \ell-k}(X)}p(c)
\Log\Big( \Phi_0(x)^{-1} \cdot \Phi_0\big(t(c)\big)\Big)\Big|_{X_i^{(1)}}^{4m}\Big\},
\end{align}
where $\mathcal{F}$ stands for the fundamental domain in $X$ 
containing the reference point $x_* \in V$. 
For $i=1, 2, \dots, d_1, \, x \in \mathcal{F}, \, N \in \mathbb{N}$
 and $c=(e_1, e_2, \dots, e_N) \in \Omega_{x, N}(X)$, we put 
$$
\mathcal{M}_N^{(i, x)}(c)=\mathcal{M}_N^{(i, x)}(\Phi_0; c):=\Log\Big(\Phi_0(x)^{-1} 
\cdot \Phi_0\big(t(c)\big)\Big) \Big|_{X_i^{(1)}}
=\sum_{j=1}^{N} \Log \big( d\Phi_0(e_j)\big)\big|_{X_i^{(1)}}.
$$
By Lemma \ref{martingale}, $\{\mathcal{M}_N^{(i, x)}\}_{N=1}^\infty$ 
is an $\mathbb{R}$-valued martingale
for every $i=1, 2, \dots, d_1$ and $x \in \mathcal{F}$. 
Therefore, we apply the Burkholder--Davis--Gundy inequality
 with the exponent $4m$ to obtain
\begin{align}\label{martingale-1}
\sum_{c \in \Omega_{x, N}(X)}p(c)\big(\mathcal{M}_N^{(i, x)}(c)\big)^{4m}&=
\sum_{c \in \Omega_{x, N}(X)}p(c) 
\Big( \sum_{j=1}^{N}  \Log \big( d\Phi_0(e_j)\big)\big|_{X_i^{(1)}}\Big)^{4m}\nn\\
&\leq \mathcal{C}_{(4m)}^{4m}\sum_{c \in \Omega_{x, N}(X)}p(c)
\Big( \sum_{j=1}^{N}\Log\big( d\Phi_0(e_j)\big)\big|_{X_i^{(1)}}^2\Big)^{2m} \nn\\
&\leq  \mathcal{C}_{(4m)}^{4m} \|d\Phi_0\|_\infty^{4m} N^{2m}
\end{align}
for $i=1, 2, \dots, d_1, \, x \in \mathcal{F}$ and $N \in \mathbb{N}$,
where $\mathcal{C}_{(4m)}$ stands for the positive constant 
which appears in the Burkholder--Davis--Gundy inequality with the exponent $4m$. 
In particular, by putting $N=\ell-k$, (\ref{martingale-1}) leads to
\begin{equation}\label{step2-2}
\sum_{c \in \Omega_{x, \ell-k}(X)}p(c)\Log\Big( \Phi_0(x)^{-1} 
\cdot \Phi_0\big(t(c)\big)\Big)\Big|_{X_i^{(1)}}^{4m} \leq 
\mathcal{C}_{(4m)}^{4m} \|d\Phi_0\|_\infty^{4m} (\ell-k)^{2m}.
\end{equation}
Thus, we obtain
$$
\mathbb{E}^{\mathbb{P}_{x_*}}
\Big[ \big\| \Log\big((\mathcal{Y}_{t_k}^{(n)})^{-1} \cdot
 \mathcal{Y}_{t_{\ell}}^{(n)}\big)\big|_{\g^{(1)}}\big\|_{\g^{(1)}}^{4m}\Big]
  \leq d_1^{2m}  \mathcal{C}_{(4m)}^{4m} \|d\Phi_0\|_\infty^{4m}
   \cdot \Big(\frac{\ell-k}{n}\Big)^{2m}
   =C^{(1)}\Big(\frac{\ell-k}{n}\Big)^{2m}
$$
by combining (\ref{step2-1}) with (\ref{step2-2}), 
which is the desired estimate (\ref{expectation-1}).

\vspace{2mm}
\noindent
{\bf Step\,3.} Next we prove (\ref{expectation-2}). 
In the similar way to (\ref{step2-1}), 
we also have
\begin{align}\label{step3-1}
&\mathbb{E}^{\mathbb{P}_{x_*}}
\Big[ \big\| \Log\big((\mathcal{Y}_{t_k}^{(n)})^{-1} \cdot
 \mathcal{Y}_{t_{\ell}}^{(n)}\big)\big|_{\g^{(2)}}\big\|_{\g^{(2)}}^{2m}\Big]\nn\\
&\leq \Big(\frac{1}{n}\Big)^{2m} \cdot d_2^{2m}
\max_{i=1, 2, \dots, d_2} \max_{x \in \mathcal{F}}
\Big\{ \sum_{c \in \Omega_{x, \ell-k}(X)}  p(c)
\Log\Big( \Phi_0(x)^{-1} \cdot \Phi_0\big(t(c)\big)\Big)\Big|_{X_{i}^{(2)}}^{2m}\Big\}.
\end{align}
An elementary inequality 
$(a_1+a_2+\dots+a_K)^{2m} \leq K^{2m-1}(a_1^{2m} + a_2^{2m}\dots + a_K^{2m})$
yields
\begin{align}\label{step3-2}
&\Log\Big( \Phi_0(x)^{-1} \cdot \Phi_0\big(t(c)\big)\Big)\Big|_{X_i^{(2)}}^{2m}\nn\\
&=\Log \Big( \Phi_0\big(o(e_1)\big)^{-1} \cdot \Phi_0\big(t(e_1)\big) \cdot \cdots
  \Phi_0\big(o(e_{\ell-k})\big)^{-1} 
  \cdot \Phi_0\big(t(e_{\ell-k})\big)\Big)\Big|_{X_i^{(2)}}^{2m}\nn\\
&=\Big(\sum_{j=1}^{\ell-k}\Log \big( d\Phi_0(e_j) \big)\big|_{X_i^{(2)}}
-\frac{1}{2}\sum_{1 \leq j_1 < j_2 \leq \ell-k}
\sum_{1 \leq \lambda < \nu \leq d_1}[\![X_\lambda^{(1)}, X_\nu^{(1)}]\!]\big|_{X_i^{(2)}} \nn\\
&\hspace{1cm}\times\Big\{ 
           \Log \big( d\Phi_0(e_{j_1})\big)\big|_{X_\lambda^{(1)}} 
           \Log \big( d\Phi_0(e_{j_2})\big)\big|_{X_\nu^{(1)}}\nn\\
&\hspace{1.5cm}-  \Log \big( d\Phi_0(e_{j_1})\big)\big|_{X_\nu^{(1)}} 
\Log \big( d\Phi_0(e_{j_2})\big)\big|_{X_\lambda^{(1)}} \Big\}\Big)^{2m}\nn\\
&\leq 3^{2m-1} \Big\{ \Big( \sum_{j=1}^{\ell-k}
\Log \big( d\Phi_0(e_j) \big)\big|_{X_{i}^{(2)}}\Big)^{2m}\nn\\
&\hspace{1cm}+L\max_{1 \leq \lambda<\nu \leq d_1}
\Big( \sum_{1 \leq j_1 < j_2 \leq \ell-k}
\Log \big( d\Phi_0(e_{j_1})\big)\big|_{X_\lambda^{(1)}} 
\Log \big( d\Phi_0(e_{j_2})\big)\big|_{X_\nu^{(1)}} \Big)^{2m} \nn\\
&\hspace{1cm}+L \max_{1 \leq \lambda<\nu \leq d_1}
\Big( \sum_{1 \leq j_1 < j_2 \leq \ell-k}\Log \big( d\Phi_0(e_{j_1})\big)\big|_{X_\nu^{(1)}} 
\Log \big( d\Phi_0(e_{j_2})\big)\big|_{X_\lambda^{(1)}} \Big)^{2m}\Big\},
\end{align}
where we put 
$$
L:=\frac{1}{2}\max_{i=1, 2, \dots, d_2}\max_{1 \leq \lambda<\nu \leq d_1}
\big|[\![X_\lambda^{(1)}, X_\nu^{(1)}]\!]\big|_{X_i^{(2)}}\big|.
$$
We fix $i=1, 2, \dots, d_2$. 
Then the Jensen inequality gives
\begin{align}\label{step3-3}
\Big( \sum_{j=1}^{\ell-k}\Log \big( d\Phi_0(e_j) \big)\big|_{X_i^{(2)}}\Big)^{2m}
&= (\ell - k)^{2m} \Big( \sum_{j=1}^{\ell-k} \frac{1}{\ell -k}  
\Log \big( d\Phi_0(e_j) \big)\big|_{X_i^{(2)}} \Big)^{2m} \nn\\
&\leq (\ell - k)^{2m} \sum_{j=1}^{\ell-k}\frac{1}{\ell-k} 
\Log \big( d\Phi_0(e_j) \big)\big|_{X_i^{(2)}}^{2m} \nn\\
&\leq (\ell-k)^{2m} \|d\Phi_0\|_\infty^{4m}.
\end{align}
For $1 \leq \lambda < \nu  \leq d_1, \, x \in \mathcal{F}, \, N \in \mathbb{N}$ 
and $c=(e_1, e_2, \dots, e_N) \in \Omega_{x, N}(X)$, we put
$$
\begin{aligned}
\widetilde{\mathcal{M}}_N^{(\lambda, \nu, x)}(c)=\widetilde{\mathcal{M}}_N^{(\lambda, \nu, x)}(\Phi_0; c)&:=
\sum_{1 \leq j_1 < j_2 \leq N}\Log \big( d\Phi_0(e_{j_1})\big)\big|_{X_\lambda^{(1)}} 
\Log \big( d\Phi_0(e_{j_2})\big)\big|_{X_\nu^{(1)}}\\
&=\sum_{j_2=2}^N \Log  \big( d\Phi_0(e_{j_2})\big)\big|_{X_\nu^{(1)}}\Big( \sum_{j_1=1}^{j_2-1}
\Log \big( d\Phi_0(e_{j_1})\big)\big|_{X_\lambda^{(1)}}\Big).
\end{aligned}
$$
We clearly observe that $\{ \widetilde{\mathcal{M}}_N^{(\lambda, \nu, x)}\}_{N=1}^\infty$ 
is an $\mathbb{R}$-valued martingale 
for every $1 \leq\lambda <\nu \leq d$ and $x \in \mathcal{F}$ due to Lemma \ref{martingale}. 
Hence, we apply the Burkholder--Davis--Gundy inequality 
with the exponent $2m$ to obtain
\begin{align}\label{step3-4}
&\sum_{c \in \Omega_{x, N}(X)}p(c)
\big(\widetilde{\mathcal{M}}_N^{(\lambda, \nu, x)}(c)\big)^{2m}\nn\\
&\leq \mathcal{C}_{(2m)}^{2m} \sum_{c \in \Omega_{x, N}(X)} p(c) \Big\{ \sum_{j_2=2}^N 
  \Log  \big( d\Phi_0(e_{j_2})\big)\big|_{X_\nu^{(1)}}^2\Big( \sum_{j_1=1}^{j_2-1}
  \Log \big( d\Phi_0(e_{j_1})\big)\big|_{X_\lambda^{(1)}}\Big)^2\Big\}^{m}\nn\\
&\leq  \mathcal{C}_{(2m)}^{2m} \sum_{c \in \Omega_{x, N}(X)} p(c) N^{m} 
\sum_{j_2=2}^{N} \frac{1}{N-1} \Log  \big( d\Phi_0(e_{j_2})\big)\big|_{X_\nu^{(1)}}^{2m}
\Big( \sum_{j_1=1}^{j_2-1}\Log \big( d\Phi_0(e_{j_1})\big)\big|_{X_\lambda^{(1)}}\Big)^{2m}\nn\\
&\leq  \mathcal{C}_{(2m)}^{2m} N^{m} 
\sum_{j_2=2}^{N} \frac{1}{N-1} \Big\{ \sum_{c \in \Omega_{x, N}(X)} p(c)
\Log  \big( d\Phi_0(e_{j_2})\big)\big|_{X_\nu^{(1)}}^{4m}\Big\}^{1/2}\nn\\
&\hspace{1cm}\times\Big\{  \sum_{c \in \Omega_{x, N}(X)} p(c) 
\Big( \sum_{j_1=1}^{j_2-1}\Log \big( d\Phi_0(e_{j_1})\big)\big|_{X_\lambda^{(1)}}\Big)^{4m}\Big\}^{1/2}\nn\\
&\leq  \mathcal{C}_{(2m)}^{2m} \|d\Phi_0\|_\infty^{2m} N^{m} 
\sum_{j_2=2}^{N} \frac{1}{N-1} \Big\{  \sum_{c \in \Omega_{x, N}(X)} p(c) 
\Big( \sum_{j_1=1}^{j_2-1}\Log \big( d\Phi_0(e_{j_1})\big)\big|_{X_\lambda^{(1)}}\Big)^{4m}\Big\}^{1/2},
\end{align}
where we used Jensen's inequality for the third line and Schwarz' inequality for the fourth line. 
Then we have
\begin{align}\label{step3-5}
 & \sum_{c \in \Omega_{x, N}(X)} p(c) \Big( \sum_{j_1=1}^{j_2-1}
 \Log \big( d\Phi_0(e_{j_1})\big)\big|_{X_\lambda^{(1)}}\Big)^{4m}\nn\\
 &\leq \mathcal{C}_{(4m)}^{4m}  \sum_{c \in \Omega_{x, N}(X)} p(c)\Big(\sum_{j_1=1}^{j_2-1} 
 \Log \big( d\Phi_0(e_{j_1})\big)\big|_{X_\lambda^{(1)}}^{2}\Big)^{2m}\nn\\
 &= \mathcal{C}_{(4m)}^{4m} (j_2-1)^{2m} \sum_{c \in \Omega_{x, N}(X)} p(c)
 \Big(\sum_{j_1=1}^{j_2-1}\frac{1}{j_2-1} \Log \big( d\Phi_0(e_{j_1})\big)\big|_{X_\lambda^{(1)}}^{2}\Big)^{2m}\nn\\
 &\leq \mathcal{C}_{(4m)}^{4m} j_2^{2m}  \sum_{c \in \Omega_{x, N}(X)} p(c)
 \sum_{j_1=1}^{j_2-1}\frac{1}{j_2-1}\Log \big( d\Phi_0(e_{j_1})\big)\big|_{X_\lambda^{(1)}}^{4m}
  \leq \mathcal{C}_{(4m)}^{4m} \|d\Phi_0\|_\infty^{4m}j_2^{2m}
\end{align}
by applying the Burkholder--Davis--Gundy inequality with the exponent $4m$.
It follows from (\ref{step3-4}) and (\ref{step3-5}) that 
\begin{align}\label{step3-6}
&\sum_{c \in \Omega_{x, N}(X)}p(c)
\big(\widetilde{\mathcal{M}}_N^{(\lambda, \nu, x)}(c)\big)^{2m} \nn\\
&\leq \mathcal{C}_{(2m)}^{2m} \|d\Phi_0\|_\infty^{2m} N^{m} 
\sum_{j_2=2}^{N} \frac{1}{N-1} \Big( \mathcal{C}_{(4m)}^{4m} \|d\Phi_0\|_\infty^{4m}j_2^{2m} \Big)^{1/2}\nn\\
&\leq  \mathcal{C}_{(2m)}^{2m}\mathcal{C}_{(4m)}^{2m} \|d\Phi_0\|_\infty^{4m} N^{m} 
\sum_{j_2=2}^{N} \frac{1}{N-1} \cdot N^m = \mathcal{C}_{(2m)}^{2m}\mathcal{C}_{(4m)}^{2m} 
\|d\Phi_0\|_\infty^{4m} N^{2m}.
\end{align}
We now put $N=\ell-k$. Then (\ref{step3-6}) implies
\begin{align}\label{step3-7}
&\sum_{c \in \Omega_{x, \ell-k}(X)}p(c)
\Big\{\Big( \sum_{1 \leq j_1 < j_2 \leq \ell-k}\Log \big( d\Phi_0(e_{j_1})\big)\big|_{X_\lambda^{(1)}} 
\Log \big( d\Phi_0(e_{j_2})\big)\big|_{X_\nu^{(1)}} \Big)^{2m} \nn\\
&\hspace{1cm}+\Big( \sum_{1 \leq j_1 < j_2 \leq \ell-k}
\Log \big( d\Phi_0(e_{j_1})\big)\big|_{X_\nu^{(1)}} 
\Log \big( d\Phi_0(e_{j_2})\big)\big|_{X_\lambda^{(1)}} \Big)^{2m}\Big\}\nn\\
&\leq 2 \mathcal{C}_{(2m)}^{2m}\mathcal{C}_{(4m)}^{2m}
 \|d\Phi_0\|_\infty^{4m} (\ell-k)^{2m} \qquad (1 \leq \lambda < \nu \leq d_1).
\end{align}
By combining (\ref{step3-1}) with (\ref{step3-2}), (\ref{step3-3}) 
and (\ref{step3-7}), we obtain 
$$
\begin{aligned}
&\mathbb{E}^{\mathbb{P}_{x_*}}
\Big[ \big\| \Log\big((\mathcal{Y}_{t_k}^{(n)})^{-1} \cdot 
\mathcal{Y}_{t_{\ell}}^{(n)}\big)\big|_{\g^{(2)}}\big\|_{\g^{(2)}}^{2m}\Big] \\
&\leq \Big(\frac{1}{n}\Big)^{2m}d_2^{2m} 3^{2m-1}\|d\Phi_0\|_\infty^{4m}
\Big\{ 1+2L\mathcal{C}_{(2m)}^{2m}\mathcal{C}_{(4m)}^{2m} 
\Big\}(\ell-k)^{2m} = C^{(2)}\Big( \frac{\ell-k}{n}\Big)^{2m}.
\end{aligned}
$$
This is the desired estimate (\ref{expectation-2}), 
and thus we have shown (\ref{tight-1}).

\vspace{2mm}
\noindent
{\bf Step\,4.} We finally prove   (\ref{tight1-1}). 
Suppose that $t_k \leq s \leq t_{k+1}$ and $t_{\ell} \leq t \leq t_{\ell+1}$
for some $1 \leq k \leq \ell \leq n$.
Since the stochastic process $\mathcal{Y}_\cdot^{(n)}$ is 
given by the $d_{\mathrm{CC}}$-geodesic interpolation,
we have
$$
\begin{aligned}
d_{\mathrm{CC}}(\mathcal{Y}_s^{(n; \,2)}, \mathcal{Y}_{t_{k+1}}^{(n; \,2)}) 
&= (k-ns)d_{\mathrm{CC}}(\mathcal{Y}_{t_k}^{(n; \,2)}, \mathcal{Y}_{t_{k+1}}^{(n; \,2)}),\\
d_{\mathrm{CC}}(\mathcal{Y}_{t_\ell}^{(n; \,2)}, \mathcal{Y}_{t}^{(n; \,2)}) 
&= (nt-\ell)d_{\mathrm{CC}}(\mathcal{Y}_{t_\ell}^{(n; \,2)}, \mathcal{Y}_{t_{\ell+1}}^{(n; \,2)}).
\end{aligned}
$$
By using (\ref{tight-1}) and the triangle inequality, we have 
$$
\begin{aligned}
&\mathbb{E}^{\mathbb{P}_{x_*}}
\Big[ d_{\mathrm{CC}}(\mathcal{Y}_{s}^{(n; \,2)}, \mathcal{Y}_{t}^{(n; \,2)})^{4m}\Big]\\
&\leq 3^{4m-1} \Big\{ (k+1-ns)^{4m} 
\cdot C\Big(\frac{1}{n}\Big)^{2m} + C \Big(\frac{\ell-k-1}{n}\Big)^{2m}
+  (nt-\ell)^{4m}\cdot C\Big(\frac{1}{n}\Big)^{2m} \Big\}\\
&\leq C\Big\{ (t_{k+1}-s)^{2m} + (t_\ell-t_{k+1})^{2m} 
+ (t-t_{\ell})^{2m}\Big\} \leq C(t-s)^{2m},
\end{aligned}
$$
which is the desired estimate (\ref{tight1-1}) and we have proved Lemma \ref{sublemma1}. \qed

\vspace{2mm}
In what follows, we write
$d\mathcal{Y}_{s, t}^{(n)*}:=(\mathcal{Y}_s^{(n)})^{-1} * \mathcal{Y}_t^{(n)}$
for $n \in \mathbb{N}$ and $0 \leq s \leq t \leq 1$.
By using Lemma \ref{sublemma1}, we obtain the following.

\begin{lm}\label{sublemma2}
For $m, n \in \mathbb{N}$, $k=1, 2, \dots, r$ and $\alpha<\frac{2m-1}{4m}$, 
there exist an $\mathcal{F}_\infty$-measurable set $\Omega_k^{(n)} \subset \Omega_{x_*}(X)$ and 
a non-negative random variable 
$\mathcal{K}_k^{(n)} \in L^{4m}(\Omega_{x_*}(X)
 \to \mathbb{R}; \,\mathbb{P}_{x_*})$ 
 such that 
$\mathbb{P}_{x_*}(\Omega_{k}^{(n)})=1$ and 
\begin{equation}\label{Holder-estimate}
d_{\mathrm{CC}}\big( \mathcal{Y}_s^{(n; \, k)}(c),  
\mathcal{Y}_t^{(n; \, k)}(c)\big)
\leq \mathcal{K}_k^{(n)}(c)(t-s)^\alpha 
\qquad (c \in \Omega_{k}^{(n)}, \, 0 \leq s \leq t \leq 1). 
\end{equation}
\end{lm}

\noindent
{\bf Proof.} We partially follow Lyons' original proof (cf.~\cite[Theorem 2.2.1]{Lyons})
for the extension theorem in rough path theory. 
We show (\ref{Holder-estimate}) by induction on the step number $k=1, 2, \dots, r$. 

\vspace{2mm}
\noindent
{\bf Step~1.} 
In the cases $k=1, 2$, we have already obtained (\ref{Holder-estimate})
in Lemma \ref{sublemma1}. 
Indeed, 
(\ref{Holder-estimate}) for $k=1, 2$ 
are readily obtained by a simple application of
the Kolmogorov--Chentov criterion 
with the bound
\begin{equation}\label{bound}
\|\mathcal{K}_k^{(n)}\|_{L^{4m}(\mathbb{P}_{x_*})} \leq \frac{5C}{(1-2^{-\theta})(1-2^{\alpha-\theta})} 
\qquad (n, m \in \mathbb{N}, \, k=1, 2),
\end{equation}
where $\theta=(2m-1)/4m$ and $C$ is a constant
independent of $n$, which appears in the right-hand side of (\ref{tight1-1}). 
See e.g., Stroock \cite[Theorem 4.3.2]{Stroock} for details. 

\vspace{2mm}
\noindent
{\bf Step~2.} Suppose that (\ref{Holder-estimate}) holds up to step $k$. 
Then, for $n \in \mathbb{N}$, there are $\mathcal{F}_\infty$-measurable
sets $\{\widehat{\Omega}_j^{(n)}\}_{j=1}^k$ and non-negative
random variables $\{\widehat{\mathcal{K}}_j^{(n)}\}_{j=1}^k$ such that
$\mathbb{P}_{x_*}(\widehat{\Omega}_{j}^{(n)})=1$ for  
$ j=1, 2, \dots, k$ and 
\begin{align}\label{Holder-estimate2}
&\big\| \big( d\mathcal{Y}_{s, t}^{(n)*}(c) \big)^{(j)}\big\|_{\mathbb{R}^{d_j}}
\leq \widehat{\mathcal{K}}_{j}^{(n)}(c)(t-s)^{j\alpha}\nn\\
&\hspace{6cm}(j=1, 2, \dots, k, \, 
c \in \widehat{\Omega}_j^{(n)}, \, 0 \leq s \leq t \leq 1)
\end{align}
with $\widehat{\mathcal{K}}_j^{(n)} \in L^{4m/j}(\Omega_{x_*}(X) \to \mathbb{R}; 
\,  \mathbb{P}_{x_*})$ for $m \in \mathbb{N}$ and $j=1, 2, \dots, k$.

We fix $0 \leq s \leq t \leq 1$ and $n \in \mathbb{N}$. 
Set $\widehat{\Omega}_{k+1}^{(n)}=\bigcap_{j=1}^k \widehat{\Omega}_j^{(n)}$. 
We denote by $\Delta$ the partition $\{s=t_0<t_1<\cdots<t_N=t\}$ 
of the time interval $[s, t]$ independent of $n \in \mathbb{N}$. 
We define two $G^{(k+1)}$-valued random variables $\mathcal{Z}_{s, t}^{(n)}$
and $\mathcal{Z}(\Delta)_{s, t}^{(n)}$ by
$$
\begin{aligned}
\big(\mathcal{Z}_{s, t}^{(n)}\big)^{(j)}&:=\begin{cases}
\big( d\mathcal{Y}_{s, t}^{(n)*}\big)^{(j)}, & (j=1, 2, \dots, k), \\
\bm{0} & (j=k+1),
\end{cases}\\
\mathcal{Z}(\Delta)_{s, t}^{(n)}&
:=\mathcal{Z}_{t_0, t_1}^{(n)} * \mathcal{Z}_{t_1, t_2}^{(n)}
* \cdots * \mathcal{Z}_{t_{N-1}, t_N}^{(n)},
\end{aligned}
$$
respectively. For $i=1, 2, \dots, d_{k+1}$, 
(\ref{CBH-formula}) and  (\ref{Holder-estimate2}) imply
$$
\begin{aligned}
&\Big| \big(\mathcal{Z}(\Delta)_{s, t}^{(n)}(c)\big)^{(k+1)}_{i*}
- \big(\mathcal{Z}(\Delta \setminus \{t_{N-1}\})_{s, t}^{(n)}(c)\big)^{(k+1)}_{i*}\Big|\\
&=\Big| \big( \mathcal{Z}_{t_{N-2}, t_{N-1}}^{(n)}(c) * 
\mathcal{Z}_{t_{N-1}, t_{N}}^{(n)}(c)\big)^{(k+1)}_{i*}
 - \big(\mathcal{Z}_{t_{N-2}, t_{N}}^{(n)}(c)\big)^{(k+1)}_{i*}\Big|\\
 &=\Bigg| \sum_{\substack{|K_1|+|K_2| = k+1 \\ |K_1|, |K_2| \geq 0}}C_{K_1, K_2}
 \mathcal{P}^{K_1}_*\big(\mathcal{Z}_{t_{N-2}, t_{N-1}}^{(n)}(c)\big)
 \mathcal{P}^{K_2}_*\big(\mathcal{Z}_{t_{N-1}, t_{N}}^{(n)}(c)\big)\Bigg|\\
 &\leq C \sum_{\substack{|K_1|+|K_2| = k+1 \\ |K_1|, |K_2| \geq 0}}
 \Big|\mathcal{P}^{K_1}_* \big(d\mathcal{Y}_{t_{N-2}, t_{N-1}}^{(n)*}(c)\big)\Big| 
 \Big| \mathcal{P}^{K_2}_*\big(d\mathcal{Y}_{t_{N-1}, t_N}^{(n)*}(c)\big)\Big|   \\
 &\leq \widehat{\mathcal{K}}_{k+1}^{(n)}(c)(t_N-t_{N-2})^{(k+1)\alpha} 
 \leq \widehat{\mathcal{K}}_{k+1}^{(n)}(c)
 \Big( \frac{2}{N-1}(t-s)\Big)^{(k+1)\alpha} \qquad (c \in \widehat{\Omega}_{k+1}^{(n)}),
\end{aligned}
$$
where the random variable $\widehat{\mathcal{K}}_{k+1}^{(n)} : \Omega_{x_*}(X) \LA \mathbb{R}$
is given by
$$
\begin{aligned}
\widehat{\mathcal{K}}_{k+1}^{(n)}(c)&:=C \sum_{\substack{|K_1|+|K_2| = k+1 \\ |K_1|, |K_2| \geq 0}}
\mathcal{Q}^{(n, K_1)}(c)\mathcal{Q}^{(n, K_2)}(c),\\
\mathcal{Q}^{(n, K)}(c)&:=\widehat{\mathcal{K}}_{k_1}^{(n)}(c) \widehat{\mathcal{K}}_{k_2}^{(n)}(c) 
 \cdots 
\widehat{\mathcal{K}}_{k_\ell}^{(n)}(c) \qquad 
\big( K=\big((i_1, k_1), (i_2, k_2), \dots, (i_\ell, k_\ell)\big) \big).
\end{aligned}
$$
Note that $\widehat{\mathcal{K}}_{k+1}^{(n)}$ is non-negative and 
it has the following integrability:
$$
\begin{aligned}
\mathbb{E}^{\mathbb{P}_{x_*}}\big[(\widehat{\mathcal{K}}_{k+1}^{(n)})^{4m/(k+1)}\big]
&\leq C\sum_{\substack{k_1, \dots, k_\ell>0 \\ k_1+k_2+\dots +k_\ell=k+1}}
\mathbb{E}^{\mathbb{P}_{x_*}}\Big[
\big(\widehat{\mathcal{K}}_{k_1}^{(n)}\widehat{\mathcal{K}}_{k_2}^{(n)}\cdots  
\widehat{\mathcal{K}}_{k_\ell}^{(n)}\big)^{4m/(k+1)}\Big]\\
&\leq C\sum_{\substack{k_1, \dots, k_\ell>0 \\ k_1+k_2+\dots +k_\ell=k+1}}
\prod_{\lambda=1}^\ell \mathbb{E}^{\mathbb{P}_{x_*}}
\Big[ \big( \widehat{\mathcal{K}}_{k_\lambda}^{(n)}\big)^{4m/k_\lambda}\Big]^{k_\lambda/(k+1)}<\infty,
\end{aligned}
$$
where we used the generalized H\"older inequality for the second line. 
By removing points in $\Delta$ successively until the partition $\Delta$ 
coincides with $\{s, t\}$, we have
\begin{align}\label{Holder2}
&\Big| \big(\mathcal{Z}(\Delta)_{s, t}^{(n)}(c)\big)_{i_*}^{(k+1)}\Big|\nn\\
&\leq\Big| \big(\mathcal{Z}(\Delta \setminus \{t_{N-1}\})_{s, t}^{(n)}(c)\big)_{i_*}^{(k+1)}\Big|
+\widehat{\mathcal{K}}_{k+1}^{(n)}(c)
 \Big( \frac{2}{N-1}(t-s)\Big)^{(k+1)\alpha} \nn\\
&\leq \Big| \big(\mathcal{Z}(\{s, t\})_{s, t}^{(n)}(c)\big)_{i_*}^{(k+1)}\Big|
+\sum_{\ell=1}^{N-2}\widehat{\mathcal{K}}_{k+1}^{(n)}(c)
\Big( \frac{2}{N-\ell}\Big)^{(k+1)\alpha}(t-s)^{(k+1)\alpha}\nn\\
&\leq  \Big| \big(\mathcal{Z}_{s, t}^{(n)}(c)\big)_{i_*}^{(k+1)}\Big|
+\widehat{\mathcal{K}}_{k+1}^{(n)}(c)  2^{(k+1)\alpha} \zeta\big( (k+1)\alpha\big) 
 (t-s)^{(k+1)\alpha} \nn\\
&\leq \widehat{\mathcal{K}}_{k+1}^{(n)}(c)(t-s)^{(k+1)\alpha} 
\qquad (i=1, 2, \dots, d_{k+1}, \, c \in \widehat{\Omega}_{k+1}^{(n)}),
\end{align}
where $\zeta(z)$ denotes the Riemann zeta function 
$\zeta(z):=\sum_{n=1}^\infty (1/n^z)$ for $z \in \mathbb{R}.$

We will show that the family $\{\mathcal{Z}(\Delta)_{s, t}^{(n)}\}$ 
satisfies the Cauchy convergence principle. 
Let $\delta>0$ and take two partitions $\Delta=\{s=t_0<t_1 \cdots <t_N=t\}$ 
and $\Delta'$ of $[s, t]$
independent of $n \in \mathbb{N}$ satisfying $|\Delta|, |\Delta'|<\delta$. 
We set $\widehat{\Delta}:=\Delta \cup \Delta'$ and write
$$
\widehat{\Delta}_\ell=\widehat{\Delta} \cap [t_\ell, t_{\ell+1}]
=\{t_\ell=s_{\ell0}<s_{\ell1}<\cdots<s_{\ell L_\ell}=t_{\ell+1}\} \qquad (\ell=0, 1, \dots, N-1).
$$
By using (\ref{Holder2}), we have
$$
\begin{aligned}
&\Big| \big(\mathcal{Z}(\Delta)_{s, t}^{(n)}(c)\big)_{i_*}^{(k+1)}
- \big(\mathcal{Z}(\widehat{\Delta})_{s, t}^{(n)}(c)\big)_{i_*}^{(k+1)}\Big|\\
&=\Big| \big(\mathcal{Z}_{t_0, t_1}^{(n)}(c) 
* \cdots * \mathcal{Z}_{t_{N-1}, t_N}^{(n)}(c)\big)^{(k+1)}_{i_*}
 - \big( \mathcal{Z}(\widehat{\Delta}_0)_{t_0, t_1}^{(n)}(c)
* \cdots * \mathcal{Z}(\widehat{\Delta}_{N-1})_{t_{N-1}, t_N}^{(n)}(c)\big)^{(k+1)}_{i_*}\Big|\\
&=\Big| \big(\mathcal{Z}_{t_0, t_1}^{(n)}(c)\big)^{(k+1)}_{i_*}
+\big(\mathcal{Z}_{t_1, t_2}^{(n)}(c) 
* \cdots * \mathcal{Z}_{t_{N-1}, t_N}^{(n)}(c)\big)^{(k+1)}_{i_*}\\
&\hspace{1cm} -\big(\mathcal{Z}(\widehat{\Delta}_0)_{t_0, t_1}^{(n)}(c)\big)^{(k+1)}_{i_*}
-\big(\mathcal{Z}(\widehat{\Delta}_1)_{t_1, t_2}^{(n)} (c)
* \cdots * \mathcal{Z}(\widehat{\Delta}_{N-1})_{t_{N-1}, t_N}^{(n)}(c)\big)^{(k+1)}_{i_*}\Big|\\
&\leq \widehat{\mathcal{K}}_{k+1}^{(n)}(c)(t_1-t_0)^{(k+1)\alpha}
+\Big| \big(\mathcal{Z}_{t_1, t_2}^{(n)}(c) 
* \cdots * \mathcal{Z}_{t_{N-1}, t_N}^{(n)}(c)\big)^{(k+1)}_{i_*}\\
&\hspace{1cm} - \big( \mathcal{Z}(\widehat{\Delta}_0)_{t_1, t_2}^{(n)}(c)
* \cdots * \mathcal{Z}(\widehat{\Delta}_{N-1})_{t_{N-1}, t_N}^{(n)}(c)\big)^{(k+1)}_{i_*}\Big|
\qquad (i=1, 2, \dots, d_{k+1}, \, c \in \widehat{\Omega}_{k+1}^{(n)}).
\end{aligned}
$$
Repeating this kind of estimate and recalling $(k+1)\alpha>1$ yield
\begin{align}\label{Holder3}
&\Big| \big(\mathcal{Z}(\Delta)_{s, t}^{(n)}(c)\big)_{i_*}^{(k+1)}
- \big(\mathcal{Z}(\widehat{\Delta})_{s, t}^{(n)}(c)\big)_{i_*}^{(k+1)}\Big|\nn\\
&\leq \sum_{\ell=0}^{N-1}\widehat{\mathcal{K}}_{k+1}^{(n)}(c)(t_{\ell+1}-t_{\ell})^{(k+1)\alpha}\nn\\
&\leq \widehat{\mathcal{K}}_{k+1}^{(n)}(c) \Big( \max_{\Delta}(t_{\ell+1}-t_{\ell})^{(k+1)\alpha-1}\Big)
\sum_{\ell=0}^{N-1}(t_{\ell+1}-t_{\ell})\nn\\
&\leq \widehat{\mathcal{K}}_{k+1}^{(n)}(c)(t-s) \times \delta^{(k+1)\alpha-1} 
\qquad (i=1, 2, \dots, d_{k+1}, \, c \in \widehat{\Omega}_{k+1}^{(n)}).
\end{align}
We thus obtain
$$
\begin{aligned}
&\Big| \big(\mathcal{Z}(\Delta)_{s, t}^{(n)}(c)\big)_{i_*}^{(k+1)}
- \big(\mathcal{Z}(\Delta')_{s, t}^{(n)}(c)\big)_{i_*}^{(k+1)}\Big|\\
&\leq \Big| \big(\mathcal{Z}(\Delta)_{s, t}^{(n)}(c)\big)_{i_*}^{(k+1)}
- \big(\mathcal{Z}(\widehat{\Delta})_{s, t}^{(n)}(c)\big)_{i_*}^{(k+1)}\Big|
+\Big| \big(\mathcal{Z}(\widehat{\Delta})_{s, t}^{(n)}(c)\big)_{i_*}^{(k+1)}
- \big(\mathcal{Z}(\widetilde{\Delta})_{s, t}^{(n)}(c)\big)_{i_*}^{(k+1)}\Big|\\
&\leq 2\widehat{\mathcal{K}}_{k+1}^{(n)}(c)(t-s) \times \delta^{(k+1)\alpha-1} \LA 0 
\qquad (i=1, 2, \dots, d_{k+1}, \, c \in \widehat{\Omega}_{k+1}^{(n)})
\end{aligned}
$$
as $\delta \searrow 0$ uniformly in $0 \leq s \leq t \leq 1$ by (\ref{Holder3}). 
Therefore, noting the estimate (\ref{Holder2}), there exists a random variable
$$
\ol{\mathcal{Z}}_{s, t}^{(n)}(c):=\begin{cases}
\dis \lim_{|\Delta| \searrow 0}\mathcal{Z}(\Delta)_{s, t}^{(n)}(c) &
 (c \in \widehat{\Omega}_{k+1}^{(n)}), \\
\bm{1}_G &  (c \in \Omega_{x_*}(X) \setminus \widehat{\Omega}_{k+1}^{(n)}),
\end{cases} \qquad (0 \leq s \leq t \leq 1)
$$
satisfying
$$
\big\| \big( \ol{\mathcal{Z}}_{s, t}^{(n)}(c) \big)^{(k+1)}\big\|_{\mathbb{R}^{d_{k+1}}}
\leq \widehat{\mathcal{K}}_{k+1}^{(n)}(c)(t-s)^{(k+1)\alpha} \qquad (c \in \widehat{\Omega}_{k+1}^{(n)}).
$$

Our final goal is to show 
$$
\ol{\mathcal{Z}}_{s, t}^{(n)}(c)
=\mathcal{Y}_{s}^{(n; \, k+1)}(c) *
\mathcal{Y}_{t}^{(n; \, k+1)}(c)
\qquad (0 \leq s \leq t \leq 1, \, c \in \widehat{\Omega}_{k+1}^{(n)}).
$$
However, it suffices to check that
\begin{align}\label{last aim}
\big( \ol{\mathcal{Z}}_{s, t}^{(n)}(c) \big)^{(k+1)}
&=\big(d\mathcal{Y}_{s, t}^{(n)*}(c)\big)^{(k+1)}
\qquad (0 \leq s \leq t \leq 1, \, c  \in \widehat{\Omega}_{k+1}^{(n)})
\end{align}
by the definition of $\ol{\mathcal{Z}}_{s, t}^{(n)}$.
We fix $i=1, 2, \dots, d_{k+1}$ and $c \in \widehat{\Omega}_{k+1}^{(n)}$. Put 
$$
\Psi^i_{s, t}(c):=\big(d\mathcal{Y}_{s, t}^{(n)*}(c)\big)^{(k+1)}_{i*} - 
\big( \ol{\mathcal{Z}}_{s, t}^{(n)}(c) \big)^{(k+1)}_{i_*} \quad (0 \leq s \leq t \leq 1).
$$
Then we easily see that $\Psi^i_{s, t}(c)$ is additive in the sense that 
\begin{equation}\label{additivity}
\Psi^i_{s, t}(c)=\Psi^i_{s, u}(c)+\Psi^i_{u, t}(c) \qquad (0 \leq s \leq u \leq t \leq 1).
\end{equation}
Since the piecewise smooth 
stochastic process $(\mathcal{Y}_t^{(n)})_{0 \leq t \leq 1}$ is defined by the $d_{\mathrm{CC}}$-
geodesic interpolation of $\{\mathcal{X}_{t_k}^{(n)}\}_{k=0}^n$, we know 
$$
\big\| \big(d\mathcal{Y}_{s, t}^{(n)*}(c)\big)^{(k+1)}\big\|_{\mathbb{R}^{d_{k+1}}}
\leq \widetilde{\mathcal{K}}^{(n)}_{k+1}(c)(t-s)^{(k+1)\alpha} 
\qquad (c \in \widetilde{\Omega}_{k+1}^{(n)})
$$
for some set $\widetilde{\Omega}_{k+1}^{(n)}$ 
with $\mathbb{P}_{x_*}(\widetilde{\Omega}_{k+1}^{(n)})=1$ and 
random variable $\widetilde{\mathcal{K}}^{(n)}_{k+1} : \Omega_{x_*}(X) \LA \mathbb{R}$. 
Then we have
$$
\Big| \Psi^i_{s, t}(c)\Big| \leq \big( \widetilde{\mathcal{K}}^{(n)}_{k+1}(c)
+\widehat{\mathcal{K}}^{(n)}_{k+1}(c)\big)(t-s)^{(k+1)\alpha}
 \qquad (0 \leq s \leq t \leq 1, \, 
 c \in \widetilde{\Omega}_{k+1}^{(n)} \cap \widehat{\Omega}_{k+1}^{(n)}).
$$
We may write $\widehat{\Omega}_{k+1}^{(n)}$ instead of
$\widetilde{\Omega}_{k+1}^{(n)} \cap \widehat{\Omega}_{k+1}^{(n)}$ by abuse of notation, because
its probability is equal to one. 
For any small $\ve>0$, there is a sufficiently large $N \in \mathbb{N}$ such that 
$1/N < \ve$. 
We obtain, as $\ve \searrow 0$, 
$$
\begin{aligned}
\Big|\Psi^i_{0, t}(c)\Big|
&= \Big| \Psi^i_{0, 1/N}(c) + \Psi^i_{1/N, 2/N}(c) + 
\cdots + \Psi^i_{[Nt]/N, t}(c)\Big|\\
&\leq \big( \widetilde{\mathcal{K}}^{(n)}_{k+1}(c)
+\widehat{\mathcal{K}}^{(n)}_{k+1}(c)\big)\ve^{(k+1)\alpha-1}
\Big\{ \underbrace{\frac{1}{N}+ \frac{1}{N}+\cdots 
+ \frac{1}{N}}_{[Nt]\text{-times}}+\Big(t-\frac{[Nt]}{N}\Big)\Big\}\\
&= \big( \widetilde{\mathcal{K}}^{(n)}_{k+1}(c)
+\widehat{\mathcal{K}}^{(n)}_{k+1}(c)\big)\ve^{(k+1)\alpha-1}t \LA 0 
\qquad (0 \leq t \leq 1, \, c \in \widehat{\Omega}_{k+1}^{(n)}) 
\end{aligned}
$$
by (\ref{additivity}) and $(k+1)\alpha-1>0$.
This implies that $\Psi^i_{0, t}(c)=0$ for
$0 \leq t \leq 1$ and $c \in \widehat{\Omega}_{k+1}^{(n)}$.  
Therefore, it follows from (\ref{last aim}) that
$$
\Psi^i_{s, t}(c)
=\Psi^i_{0, t}(c)-\Psi^i_{0, s}(c)=0 
\qquad (0 \leq s \leq t \leq 1, \, c \in \widehat{\Omega}_{k+1}^{(n)}),
$$
which means (\ref{Holder3}). 
Consequently,  there exist a $\mathcal{F}_\infty$-measurable set $\Omega_{k+1}^{(n)} \subset \Omega_{x_*}(X)$
with probability one and a non-negative random variable 
$\mathcal{K}_{k+1}^{(n)} \in L^{4m}(\Omega_{x_*}(X) \to \mathbb{R}; \,\mathbb{P}_{x_*})$
satisfying
$$
\begin{aligned}
d_{\mathrm{CC}}\big( \mathcal{Y}_s^{(n; \, k+1)}(c),  
\mathcal{Y}_t^{(n; \, k+1)}(c)\big)
&\leq \mathcal{K}_{k+1}^{(n)}(c)(t-s)^\alpha \quad (0 \leq s \leq t \leq 1, \, c \in \Omega_{k+1}^{(n)}).
\end{aligned}
$$
This completes the proof of Lemma \ref{sublemma2}. \qed

\vspace{2mm}
\noindent
{\bf Proof of Lemma \ref{tightness1}}. 
For $m, n \in \mathbb{N}$ and $\widehat{\alpha}<\frac{2m-1}{4m}$, Lemma \ref{sublemma2} implies that
$$
\begin{aligned}
\mathbb{E}^{\mathbb{P}_{x_*}}\Big[d_{\mathrm{CC}}\big( \mathcal{Y}_s^{(n; \, r)},  
\mathcal{Y}_t^{(n; \, r)}\big)^{4m}\Big]
&\leq \mathbb{E}^{\mathbb{P}_{x_*}}\big[\big(\mathcal{K}_r^{(n)}\big)^{4m}\big](t-s)^{4m\widehat{\alpha}} 
\end{aligned}
$$
for $0 \leq s \leq t \leq 1$. 
Thus, it follows from (\ref{bound}) that
$$
\mathbb{E}^{\mathbb{P}_{x_*}}\Big[d_{\mathrm{CC}}\big( \mathcal{Y}_s^{(n; \, r)},  
\mathcal{Y}_t^{(n; \, r)}\big)^{4m}\Big]
\leq C(t-s)^{4m\widehat{\alpha}} \qquad (0 \leq s \leq t \leq 1). 
$$ 
for a positive constant $C>0$ independent of $n \in \mathbb{N}$. 
By  applying
the Kolmogorov tightness criterion (cf.~Friz--Hairer \cite[Section 3.1]{FH}), we know that 
the family $\{\Prob^{(n)}\}_{n=1}^\infty$ is tight in 
$C^{0, \alpha\text{\rm{-H\"ol}}}_{\bm{1}_G}([0, 1]; G)$ for 
$\alpha <\frac{4m\widehat{\alpha}-1}{4m}<\frac{1}{2}-\frac{1}{2m}$. 
Since $m \in \mathbb{N}$ is taken arbitrarily, we complete the proof. 
\qed

\vspace{2mm}
We conclude Theorem \ref{FCLT1} by showing the following 
convergence of the finite dimensional distribution. 

\begin{lm}\label{CFDD1}
Let $\ell \in \mathbb{N}$. For fixed $0 \leq s_1 < s_2 < \cdots < s_\ell \leq 1$, we have
$$
(\mathcal{Y}_{s_1}^{(n)}, \mathcal{Y}_{s_2}^{(n)}, \dots, \mathcal{Y}_{s_\ell}^{(n)})
\overset{(d)}{\LA} (Y_{s_1}, Y_{s_2}, \dots, Y_{s_\ell}) \quad (n \to \infty).
$$
\end{lm}

\noindent
{\bf Proof.} We only prove the convergence for $\ell=2$. 
General cases $(\ell \geq 3)$ can be also proved by repeating the same argument. 
Put $s=s_1$ and $t=s_2$. 
Then, by applying Theorem \ref{CLT1}, we obtain
$(\mathcal{X}_{s}^{(n)}, 
\mathcal{X}_{t}^{(n)})
\overset{(d)}{\LA} (Y_{s}, Y_{t})$
as $n \to \infty$ in the same way as \cite[Lemma 4.2]{IKK}.
On the other hand, Lemma \ref{sublemma2} tells us that
there exists a non-negative random variable 
$\mathcal{K}_r^{(n)}\in 
L^{4m}(\Omega_{x_*}(X) \to \mathbb{R}; \,\mathbb{P}_{x_*})$
such that
$$
d_{\mathrm{CC}}
\big( \mathcal{Y}_s^{(n)}(c), \mathcal{Y}_t^{(n)}(c) \big)
\leq \mathcal{K}_r^{(n)}(c)(t-s)^\alpha \quad \mathbb{P}_{x_*}\text{-a.s.} 
\qquad (0 \leq s \leq t \leq 1).
$$
Now suppose that $t_k \leq t \leq t_{k+1}$ for some $k=0, 1, \dots, n-1$. 
For all $\ve>0$ and sufficiently large $m \in \mathbb{N}$, 
by using Chebyshev's inequality, we have
$$
\begin{aligned}
&
\mathbb{P}_{x_*}\Big(d_{\mathrm{CC}}
\big( \mathcal{X}_t^{(n)}, \mathcal{Y}_t^{(n)} \big) >\ve\Big)\\
&\leq \frac{1}{\ve^{4m}}\mathbb{E}^{\mathbb{P}_{x_*}}
\Big[ d_{\mathrm{CC}}
\big( \mathcal{X}_t^{(n)}, \mathcal{Y}_t^{(n)} \big)^{4m}\Big]\\
&\leq \frac{1}{\ve^{4m}}\mathbb{E}^{\mathbb{P}_{x_*}}
\Big[ d_{\mathrm{CC}}
\big( \mathcal{Y}_{t_k}^{(n)}, \mathcal{Y}_{t_{k+1}}^{(n)} \big)^{4m}\Big]\\
&\leq \frac{1}{\ve^{4m}}\mathbb{E}^{\mathbb{P}_{x_*}}
\Big[ (\mathcal{K}_r^{(n)})^{4m}(t_{k+1}-t_k)^{4m\alpha}\Big]=
\frac{1}{n^{2m-1}\ve^{4m}}\mathbb{E}^{\mathbb{P}_{x_*}}
\big[(\mathcal{K}_r^{(n)})^{4m}\big] \LA 0 \quad (n \to \infty).
\end{aligned}
$$
Thus, Slutzky's theorem (cf.~Klenke \cite[Theorem 13.8]{Klenke}) allows us to obtain 
the desired convergence
$(\mathcal{Y}_{s}^{(n)}, 
\mathcal{Y}_{t}^{(n)})
\overset{(d)}{\LA} (Y_{s}, Y_{t})$
as $n \to \infty$. 
This completes the proof. \qed

\subsection{Proof of Theorem \ref{FCLT1-general}}

In this section, we show Theorem \ref{FCLT1-general}, which is a generalization of Theorem \ref{FCLT1}
to non-harmonic cases. Our first aim is to show that 
the same pathwise H\"older estimate as Lemma \ref{sublemma2}
holds for the stochastic process $\{\ol{\mathcal{Y}}^{(n)}\}_{n=1}^\infty$. 

\begin{lm}\label{sublemma-general}
For $m, n \in \mathbb{N}$ and $\alpha<\frac{2m-1}{4m}$, 
there exist an $\mathcal{F}_\infty$-measurable set $\ol{\Omega}_r^{(n)} \subset \Omega_{x_*}(X)$ and 
a non-negative random variable 
$\ol{\mathcal{K}}_r^{(n)} \in L^{4m}\big(\Omega_{x_*}(X) \to \mathbb{R}; \mathbb{P}_{x_*}\big)$
such that 
$\mathbb{P}_{x_*}(\ol{\Omega}_r^{(n)})=1$ and 
\begin{equation}\label{pathwise-general}
d_{\mathrm{CC}}\big(\ol{\mathcal{Y}}_{s}^{(n)}(c), \ol{\mathcal{Y}}_{t}^{(n)}(c)\big)
\leq \ol{\mathcal{K}}_r^{(n)}(c)(t-s)^\alpha 
\qquad (c \in \ol{\Omega}_{r}^{(n)}, \, 0 \leq s <t \leq 1). 
\end{equation}
\end{lm}

\noindent
{\bf Proof.} Fix $n \in \mathbb{N}$ and $1 \leq k \leq \ell \leq n$. 
We then have
$$
d_{\mathrm{CC}}(\ol{\mathcal{Y}}_{t_k}^{(n)}, \ol{\mathcal{Y}}_{t_{\ell}}^{(n)})
\leq d_{\mathrm{CC}}(\ol{\mathcal{Y}}_{t_k}^{(n)}, \mathcal{Y}_{t_k}^{(n)})
+d_{\mathrm{CC}}(\mathcal{Y}_{t_k}^{(n)}, \mathcal{Y}_{t_\ell}^{(n)})
+d_{\mathrm{CC}}(\mathcal{Y}_{t_\ell}^{(n)}, \ol{\mathcal{Y}}_{t_\ell}^{(n)}). 
$$
Set set 
$\mathcal{Z}^{(n)}_t:=(\mathcal{Y}^{(n)}_t)^{-1} * \ol{\mathcal{Y}}^{(n)}_t$ 
for $0 \leq t \leq 1$ and $n \in \mathbb{N}$. 
We then have
$$
\log\big(\mathcal{Z}_{t_k}^{(n)})|_{\g^{(1)}}=\frac{1}{\sqrt{n}}\mathrm{Cor}_{\g^{(1)}}(w_{k}) 
\qquad (n \in \mathbb{N}, \, k=0, 1, \dots, n)
$$
so that there is a constant $C>0$ such that 
$\big\| \log\big(\mathcal{Z}_{t_k}^{(n)})|_{\g^{(1)}} \big\|_{\g^{(1)}} \leq Cn^{-1/2}$
for $n \in \mathbb{N}$ and $k=0, 1, 2, \dots, n$. 
It follows from the choice of the components of $\Phi_0(x) \, (x \in V)$ that 
$\big\| \log\big(\mathcal{Z}_{t_k}^{(n)})|_{\g^{(i)}} \big\|_{\g^{(i)}} \leq Cn^{-i/2}$
for $n \in \mathbb{N}$ and $k=0, 1, 2, \dots, n$. 
By using Proposition \ref{homogeneous equiv}, we have
\begin{equation}\label{corrector-est-1}
d_{\mathrm{CC}}(\ol{\mathcal{Y}}_{t_k}^{(n)}, \mathcal{Y}_{t_k}^{(n)})
\leq C\big\|\mathcal{Z}_{t_k}^{(n)}\big\|_{\Hom}
=C\sum_{i=1}^r \big\| \log\big(\mathcal{Z}_{t_k}^{(n)})|_{\g^{(i)}} \big\|_{\g^{(i)}}^{1/i}
\leq \frac{C}{\sqrt{n}}
\end{equation}
for $n \in \mathbb{N}$ and $k=0, 1, 2, \dots, n$. 
Then Lemma \ref{sublemma2} and (\ref{corrector-est-1}) imply that there exist
an $\mathcal{F}_\infty$-measurable set $\ol{\Omega}_r^{(n)} \subset \Omega_{x_*}(X)$ and  
a non-negative random variable 
$\ol{\mathcal{K}}_r^{(n)} \in L^{4m}\big(\Omega_{x_*}(X) \to \mathbb{R}; \mathbb{P}_{x_*}\big)$
such that 
$\mathbb{P}_{x_*}^{(n^{-1/2})}(\ol{\Omega}_r^{(n)})=1$ and 
\begin{align}
d_{\mathrm{CC}}\big(\ol{\mathcal{Y}}_{t_k}^{(n)}(c), \ol{\mathcal{Y}}_{t_\ell}^{(n)}(c)\big)
&\leq \frac{C}{\sqrt{n}} + \mathcal{K}_r^{(n)}(c)\Big(\frac{\ell-k}{n}\Big)^\alpha + \frac{C}{\sqrt{n}}\nn \\
&\leq \ol{\mathcal{K}}_r^{(n)}(c)\Big(\frac{\ell-k}{n}\Big)^\alpha
\qquad (c \in \ol{\Omega}_r^{(n)}, \, 0 \leq k \leq \ell \leq n). \label{general-est-partial}
\end{align}
For $0 \leq s < t \leq 1$, take $0 \leq k \leq \ell \leq n$ such that
$k/n \leq s < (k+1)/n$ and $\ell/n \leq t < (\ell+1)/n$. 
By the definition of $(\ol{\mathcal{Y}}_t^{(n)})_{0 \leq t \leq 1}$, we have
$$
\begin{aligned}
d_{\mathrm{CC}}\big(\ol{\mathcal{Y}}_{s}^{(n)}, \ol{\mathcal{Y}}_{t_{k+1}}^{(n)}\big)
&=(k-ns)d_{\mathrm{CC}}\big(\ol{\mathcal{Y}}_{t_k}^{(n)}, \ol{\mathcal{Y}}_{t_{k+1}}^{(n)}\big),\\
d_{\mathrm{CC}}\big(\ol{\mathcal{Y}}_{t_\ell}^{(n)}, \ol{\mathcal{Y}}_{t}^{(n)}\big)
&=(nt-\ell)d_{\mathrm{CC}}\big(\ol{\mathcal{Y}}_{t_\ell}^{(n)}, \ol{\mathcal{Y}}_{t_{\ell+1}}^{(n)}\big).
\end{aligned}
$$
We then use the triangular inequality and (\ref{general-est-partial}) to obtain
$$
\begin{aligned}
&d_{\mathrm{CC}}\big(\ol{\mathcal{Y}}_{s}^{(n)}(c), \ol{\mathcal{Y}}_{t}^{(n)}(c)\big)\\
&\leq (k-ns)  \ol{\mathcal{K}}_r^{(n)}(c)\Big(\frac{1}{n}\Big)^\alpha
+\ol{\mathcal{K}}_r^{(n)}(c)\Big(\frac{\ell-k-1}{n}\Big)^\alpha + (nt-\ell)\ol{\mathcal{K}}_r^{(n)}(c)
\Big(\frac{1}{n}\Big)^\alpha\\
&\leq \ol{\mathcal{K}}_r^{(n)}(c) \Big\{ \Big(\frac{k+1}{n}-s\Big)^\alpha + \Big(\frac{\ell-k-1}{n}\Big)^\alpha
+\Big(t - \frac{\ell}{n}\Big)^\alpha\Big\} 
\leq \ol{\mathcal{K}}_r^{(n)}(c) (t-s)^\alpha \qquad (c \in \ol{\Omega}_r^{(n)}).
\end{aligned}
$$
This completes the proof. \qed

\vspace{2mm}
\noindent
{\bf Proof of Theorem \ref{FCLT1-general}.}
The proof is split into two steps.

\vspace{2mm}
\noindent
{\bf Step~1.} We show that 
$\{\ol{\mathcal{Y}}^{(n)}\}_{n=1}^\infty$ converges in law
to $(Y_t)_{0 \leq t \leq 1}$ in 
$C_{\bm{1}_G}([0, 1]; G)$ as $n \to \infty$. 
For $0 \leq t \leq 1$, take an integer $0 \leq k \leq n$ such that $k/n \leq t <(k+1)/n$. 
Then
(\ref{Holder-estimate}), (\ref{pathwise-general}) and (\ref{corrector-est-1})
imply, $\mathbb{P}_{x_*}$-almost surely, 
\begin{align}
d_{\mathrm{CC}}(\mathcal{Y}_t^{(n)}, \ol{\mathcal{Y}}_t^{(n)})
&\leq d_{\mathrm{CC}}(\mathcal{Y}_{t_k}^{(n)}, \mathcal{Y}_t^{(n)})
+d_{\mathrm{CC}}(\mathcal{Y}_{t_k}^{(n)}, \ol{\mathcal{Y}}_{t_k}^{(n)})
+d_{\mathrm{CC}}(\ol{\mathcal{Y}}_{t_k}^{(n)}, \ol{\mathcal{Y}}_t^{(n)})\nn\\
&\leq \mathcal{K}_r^{(n)}\Big(t-\frac{k}{n}\Big)^\alpha + 
\frac{C}{\sqrt{n}} + \ol{\mathcal{K}}_r^{(n)}\Big(t-\frac{k}{n}\Big)^\alpha\nn\\
&\leq \big\{ \mathcal{K}_r^{(n)}+\ol{\mathcal{K}}_r^{(n)}+C\big\}\Big(\frac{1}{\sqrt{n}}\Big)^\alpha
\qquad \Big(m \in \mathbb{N}, \, \alpha<\frac{2m-1}{4m}\Big). \label{difference-y}
\end{align}
Let $\rho$ be a metric on $C_{\bm{1}_G}([0, 1]; G)$ defined by
$\rho(w^{(1)}, w^{(2)}):=
\max_{0 \leq t \leq 1}d_{\mathrm{CC}}(w_t^{(1)}, w_t^{(2)}).$
By applying the Chebyshev inequality and (\ref{difference-y}), we have, for 
$\ve >0$ and $m \in \mathbb{N}$, 
$$
\begin{aligned}
&\mathbb{P}_{x_*}\Big(\rho(\mathcal{Y}^{(n)}, \ol{\mathcal{Y}}^{(n)})>\ve\Big)\\
&\leq \Big(\frac{1}{\ve}\Big)^{4m}
\mathbb{E}^{\mathbb{P}_{x_*}}
\Big[\rho(\mathcal{Y}^{(n)}, \ol{\mathcal{Y}}^{(n)})^{4m}\Big]\\
&\leq \Big(\frac{1}{\ve}\Big)^{4m}
\mathbb{E}^{\mathbb{P}_{x_*}}\Big[\max_{0 \leq t \leq 1}
d_{\mathrm{CC}}(\mathcal{Y}_t^{(n)}, \ol{\mathcal{Y}}_t^{(n)})^{4m}\Big]\\
&\leq 3^{4m-1}\Big(\frac{1}{\ve}\Big)^{4m}\Big(\frac{1}{\sqrt{n}}\Big)^{4m\alpha}
\Big\{ \mathbb{E}^{\mathbb{P}_{x_*}}\big[(\mathcal{K}_r^{(n)})^{4m}\big]
+\mathbb{E}^{\mathbb{P}_{x_*}}\big[(\ol{\mathcal{K}}_r^{(n)})^{4m}\big]
+\mathbb{E}^{\mathbb{P}_{x_*}}\big[C^{4m}\big]\Big\} \to 0
\end{aligned}
$$
as $n \to \infty$. 
Therefore, by Slutzky's theorem, 
the convergence in law of $\{\ol{\mathcal{Y}}^{(n)}\}_{n=1}^\infty$ to the 
diffusion process $(Y_t)_{0 \leq t \leq 1}$ in $C_{\bm{1}_G}([0, 1]; G)$ as $n \to \infty$ is obtained. 

\vspace{2mm}
\noindent
{\bf Step~2.}
By the previous step, we see that the convergence of 
the finite-dimensional distribution of $\{\ol{\mathcal{Y}}^{(n)}\}_{n=1}^{\infty}$ holds. 
On the other hand, we can prove
 that the sequence of probability measures 
$\{\ol{\Prob}^{(n)}:=\mathbb{P}_{x_*} \circ (\ol{\mathcal{Y}}^{(n)})^{-1}\}_{n=1}^\infty$
is tight
in $C_{\bm{1}_G}^{0, \alpha{{\normalfont \hol}}}([0, 1]; G)$, 
by applying Lemma \ref{sublemma-general} and
 by following the same argument as the proof of Lemma \ref{tightness1}. 
We complete the proof by combining these two. \qed




\section{An explicit representation of the limiting diffusions and a relation with rough path theory}

\subsection{An explicit representation of the limiting diffusion}

Let us consider an SDE on $\mathbb{R}^N$
\begin{equation}\label{general SDE}
d\xi_t=\sum_{i=1}^d U_i(\xi_t) \circ dB_t^i + U_0(\xi_t) \, dt, \qquad \xi_0=x_0 \in \mathbb{R}^N,
\end{equation}
where $U_0, U_1, \dots, U_d$ are $C^\infty$-vector fields on $\mathbb{R}^d$ and
$(B_t)_{0 \leq t \leq 1}=(B_t^1, B_t^2, \dots, B_t^d)_{0 \leq t \leq 1}$
 is a $d$-dimensional standard Brownian motion. 
The symbol $\circ$ denotes the usual Stratonovich type stochastic integral. 
As is well-known, a number of authors have studied explicit
representations of the unique solution 
to (\ref{general SDE}) as a functional of 
It\^o/Stratonovich iterated integrals under some assumptions 
on vector fields $U_0, U_1, \dots, U_d$. 
In particular,  Kunita \cite{Kunita} has obtained 
the explicit formula by using the CBH formula 
in the case where the Lie algebra generated by $U_0, U_1, \dots, U_d$ 
is nilpotent or solvable. 
Castell \cite{Castell} 
gave a universal representation formula, 
which contains the above results in the nilpotent case and 
extends the study of Ben Arous \cite{BA} to more general diffusions. 

We now recall the result in \cite{Castell} when the Lie algebra
generated by $U_0, U_1, \dots, U_d$ is nilpotent of step $r$.
We first introduce several notations of multi-indices. 
Set $\mathcal{I}^{(k)}=\{0, 1, \dots, d\}^k$ and 
let $I=(i_1, i_2, \dots, i_k) \in \mathcal{I}^{(k)}$ be a multi-index of length $|I|=k$. 
For vector fields $U_0, U_1, \dots, U_d$ on $\mathbb{R}^d$ and 
$I=(i_1, i_2, \dots, i_k) \in \mathcal{I}^{(k)}$, we denote by $U^I$
the vector field of the form 
$U^I=[U_{i_1}, [U_{i_2}, \cdots, [U_{i_{k-1}}, U_{i_k}] \cdots ]].$
For a multi-index $I=(i_1, i_2, \dots, i_k) \in \mathcal{I}^{(k)}$,
 we define the Stratonovich iterated integral $B_t^I$ by
$$
B_t^I := \int_{\Delta^{(k)}[0, t]}\circ dB_{t_1}^{i_1} \circ dB_{t_2}^{i_2} \cdots \circ dB_{t_k}^{i_k},
$$
where $\Delta^{(k)}[0, t]:=\{(t_1, t_2, \dots, t_k) \in [0, t]^k \, | \, 0<t_1<t_2<\cdots<t_k<t\}$
for $0 \leq t \leq 1$
and $B_t^0=t$ for convention. 
Next we introduce notations of the permutations. 
Denote by $\frak{S}_k$ be 
the symmetric group of degree $k$. For a permutation $\sigma \in \frak{S}_k$, 
we write $e(\sigma)$ for the cardinality of the set 
$\{i \in \{1, 2, \dots, k-1\} \, | \, \sigma(i)>\sigma(i+1)\}$,
which we call the number of inversions of $\sigma$.
For $I=(i_1, i_2, \dots, i_k) \in \mathcal{I}^{(k)}$ and $\sigma \in \frak{S}_k$, we put
$I_\sigma:=(i_{\sigma(1)}, i_{\sigma(2)}, \dots, i_{\sigma(k)}) \in \mathcal{I}^{(k)}. $

\begin{pr}\label{Castell's thm}
{\bf (cf.~\cite{Castell})} 
Let $U_0, U_1, \dots, U_d$ be bounded $C^\infty$-vector fields on $\mathbb{R}^N$ such that 
the Lie algebra generated by $U_0, U_1, \dots, U_d$ is nilpotent of step $r$. 
We consider the solution $(\xi_t)_{0 \leq t \leq 1}$ 
of {\normalfont (\ref{general SDE})}. 
Then we have
$$
\xi_t=\exp\Big( \sum_{k=1}^{r} \sum_{I \in \mathcal{I}^{(k)}}c_t^I U^I\Big)(x_0) 
\qquad (0 \leq t \leq 1) \quad a.s.,
$$
where 
$$
c_t^I:=\sum_{\sigma \in \frak{S}_{|I|}}\frac{(-1)^{e(\sigma)}}{|I|^2 \dis \binom{|I|-1}{e(\sigma)}}B_t^{I_{\sigma^{-1}}}. 
$$
\end{pr}
Here we give several concrete computations of $c_t^IU^I$. 

\vspace{1mm}
\noindent
$\bullet$ If $I=(i) \in \mathcal{I}^{(1)}$, we see $c_t^I=B_t^i$ for $0 \leq t \leq 1$ and $i=0, 1, \dots, d$.
Therefore, we have 
$$
\sum_{I \in \mathcal{I}^{(1)}}c_t^IU^I=\sum_{i=0}^d B_t^iU_i=tU_0+\sum_{i=1}^d B_t^i U_i.
$$

\noindent 
$\bullet$ If $I=(i, j) \in \mathcal{I}^{(2)}$ with $i \neq j$, we also see
$$
c_t^I=\begin{cases}
\dis \frac{1}{4}\int_0^t \int_0^u (\circ dB_s^i \circ dB_u^j - \circ dB_s^j \circ dB_u^i) & (i<j),\\
\vspace{-2mm}\\
\dis -\frac{1}{4}\int_0^t \int_0^u (\circ dB_s^i \circ dB_u^j - \circ dB_s^j \circ dB_u^i) & (i>j).\\
\end{cases}
$$
Since $[U_i, U_j]=-[U_j, U_i]$ holds for $i \neq j$, we have
$$
\begin{aligned}
\sum_{I \in \mathcal{I}^{(2)}}c_t^IU^I&=\sum_{0 \leq i < j \leq d}
\frac{1}{2}\int_0^t \int_0^u (\circ dB_s^i \circ dB_u^j - \circ dB_s^j \circ dB_u^i)[U_i, U_j]\\
&=\sum_{0 \leq i < j \leq d}
\frac{1}{2}\int_0^t (B_s^i  dB_s^j - B_s^j dB_s^i)[U_i, U_j].
\end{aligned}
$$
The stochastic integral
$$
\frac{1}{2}\int_0^t (B_s^idB_s^j - B_s^jdB_s^i) \qquad (0 \leq t \leq 1, \, 1 \leq i < j \leq d)
$$
indicates the well-known {\it L\'evy's stochastic area} enclosed by the Brownian curve
$\{(B_s^i, B_s^j) \in \mathbb{R}^2 \, | \, 0 \leq s \leq t\}$ and its chord.

We now provide an explicit representation of $(Y_t)_{0 \leq t \leq 1}$, 
the solution to the SDE (\ref{SDE}). 
As mentioned in Section 3.1, 
since $G$ is identified with $\mathbb{R}^N \, (N=d_1+d_2+\cdots+d_r)$,
we may apply Proposition \ref{Castell's thm} by replacing 
$U_0, U_1, \dots, U_{d}$ by $V_0, V_1, \dots, V_{d_1}$, where  
$V_0=\beta(\Phi_0)_*$. Then we have

\begin{tm}\label{limiting diffusion}
The limiting diffusion process $(Y_t)_{0 \leq t \leq 1}$ is explicitly represented as
\begin{align}\label{rep}
Y_t
&=\exp\Big( t\beta(\Phi_0)_*+\sum_{i=1}^{d_1} B_t^i V_{i*} \nn\\
&\hspace{1cm}+ \sum_{0 \leq i < j \leq {d_1}}
\frac{1}{2}\int_0^t (B_s^i  dB_s^j - B_s^j dB_s^i) [\![V_{i*}, V_{j*}]\!]+
\sum_{k=3}^{r} \sum_{I \in \mathcal{I}^{(k)}}c_t^I V_*^I\Big)(\bm{1}_G),
\end{align}
where $V_*^I=[\![V_{i_1*}, [\![V_{i_2*}, \cdots,  [\![V_{i_k-1*}, V_{i_k*}]\!] \cdots ]\!] ]\!]$
for $I=(i_1, i_2, \dots, i_k) \in \mathcal{I}^{(k)}$. 
\end{tm}
We should note that some of $[\![V_{i*}, V_{j*}]\!] \, (0 \leq i < j \leq d_1)$
in (\ref{rep}) may vanish because $\{[\![V_{i*}, V_{j*}]\!]\}_{1 \leq i <j \leq d}$ is not always
linearly independent.

\begin{re}\label{diffusions}
We will discuss yet another functional CLT 
for non-symmetric random walks on a nilpotent covering graph 
in a subsequent paper {\normalfont \cite{IKN}}.
We also obtain a $G$-valued diffusion process 
whose infinitesimal generator differs from $-\A$ of $(Y_t)_{0 \leq t \leq 1}$
through another CLT. 
Precisely, the generator is the sub-Laplacian plus 
drift of the asymptotic direction $\rho_{\mathbb{R}}(\gamma_p) \in \g^{(1)}$, and 
the corresponding diffusion process $(\widehat{Y}_t)_{0 \leq t \leq 1}$ is given by
$$
\begin{aligned}
\widehat{Y}_t &=\exp\Big( t\rho_{\mathbb{R}}(\gamma_p)_*+\sum_{i=1}^{d_1} B_t^i V_{i*} \nn\\
&\hspace{1cm}+ \sum_{0 \leq i < j \leq {d_1}}
\frac{1}{2}\int_0^t (B_s^i  dB_s^j - B_s^j dB_s^i) [\![V_{i*}, V_{j*}]\!]+
\sum_{k=3}^{r} \sum_{I \in \mathcal{I}^{(k)}}c_t^I V_*^I\Big)(\bm{1}_G),
\end{aligned}
$$
where $V_{0}=\rho_{\mathbb{R}}(\gamma_p) \in \g^{(1)}$. 
We see that these two diffusions are completely same when the random walk on $X$
is $m$-symmetric. However, 
the difference between them appears in the case
where $\gamma_p \neq 0$ and $\rho_{\mathbb{R}}(\gamma_p) = \bm{0}_{\g}.$
Namely, $(Y_t)_{0 \leq t \leq 1}$ is still given by {\normalfont (\ref{rep})}, 
while $(\widehat{Y}_t)_{0 \leq t \leq 1}$ is nothing but the ``Brownian motion on $G$''
given by
$$
\widehat{Y}_t =\exp\Big(\sum_{i=1}^{d_1} B_t^i V_{i*} 
+ \sum_{1 \leq i < j \leq {d_1}}
\frac{1}{2}\int_0^t (B_s^i  dB_s^j - B_s^j dB_s^i) [\![V_{i*}, V_{j*}]\!]+ \cdots \Big)(\bm{1}_G). 
$$
\end{re}
Before closing this subsection,  
we prove the following, which was mentioned in Section 2. 
\begin{pr}\label{SDE-solution}
The $C_0$-semigroup $(e^{-t\A})_{0 \leq t \leq 1}$  
coincides with the $C_0$-semigroup $(T_t)_{0 \leq t \leq 1}$ on $C_\infty(G)$ defined by 
$T_t f(g)=\mathbb{E}[f(Y_t^g)]$ for $g \in G$, where 
$(Y_t^g)_{0 \leq t \leq 1}$ is a solution to the stochastic differential equation
\begin{equation}\label{g-SDE}
dY_t^g=\sum_{i=1}^{d_1}V_{i*}(Y_t^g) \circ dB_t^i + \beta(\Phi_0)_*(Y_t^g) \, dt,
\qquad Y_0^g=g \in G.
\end{equation}
\end{pr}

\noindent
{\bf Proof.} 
By recalling Lemma \ref{operator-lem}, the linear operator $\A$
satisfies the maximal dissipativity, that is, $\lambda - \A$ is surjective for some $\lambda>0$. 
Therefore, the Lumer--Fillips theorem implies that 
$(\e^{-t\A})_{0 \leq t \leq 1}$ is the unique Feller semigroup on $C_\infty(G)$
whose infinitesimal generator extends $\big(-\A, C_0^\infty(G)\big)$. 
By applying It\^o's formula to (\ref{g-SDE}), we easily see that
the generator of $(Y_t)_{0 \leq t \leq 1}$ coincides with
$-\A$ on $C_0^\infty(G)$. Therefore,
it suffices to show that the semigroup $(T_t)_{0 \leq t \leq 1}$
enjoys the Feller property, that is, $T_t\big(C_\infty(G)\big) \subset C_\infty(G)$ for $0 \leq t \leq 1$. 

Suppose $f \in C_\infty(G)$. 
For any $\ve >0$, we choose a sufficiently large $R>0$ such that 
$|f(g)|<\ve$ for $g \in B_R(\bm{1}_G)^c$, where
$B_R(\bm{1}_G):=\{g \in G \, | \, d_{\mathrm{CC}}(\bm{1}_G, g) < R\}$. 
Then, for $g \in B_{2R}(\bm{1}_G)^c$, we have
$$
\begin{aligned}
|T_tf(g)| &\leq \mathbb{E}\big[ |f(Y_t^g)| \, : \, d_{\mathrm{CC}}(g, Y_t^g) < R\big]
+ \mathbb{E}\big[ |f(Y_t^g)| \, : \, d_{\mathrm{CC}}(g, Y_t^g) \geq R\big]\\
&\leq \ve + \|f\|_\infty^G \mathbb{P}\big( d_{\mathrm{CC}}(g, Y_t^g) \geq R\big).
\end{aligned}
$$
By combining
Proposition \ref{homogeneous equiv} and the Chebyshev inequality 
with Theorem \ref{limiting diffusion}, 
$$
\begin{aligned}
\mathbb{P}\big( d_{\mathrm{CC}}(g, Y_t^g) \geq R\big)
&=\mathbb{P}\big( d_{\mathrm{CC}}(\bm{1}_G, Y_t) \geq R\big)\\
&\leq \mathbb{P}\big( C\|Y_t\|_{\mathrm{Hom}} \geq R\big)
\leq \frac{C}{R^2} \mathbb{E}
\Big[ \Big(\sum_{k=1}^r \Big\|\sum_{I \in \mathcal{I}^{(k)}}c_t^I V_*^I\Big\|_{\g^{(k)}}^{1/k}\Big)^2\Big].
\end{aligned}
$$
Now we recall the following fact (cf.~Friz--Riedel \cite[Lemma 2]{FR}): 
For a multi-index $I=(i_1, i_2, \dots, i_k) \in \mathcal{I}^{(k)}$, 
there exists a constant $C$ depending only on $k$ such that 
$$
\mathbb{E}\Big[\Big(\int_{\Delta^{(k)}[0, t]}
\circ dB_{t_1}^{i_1} \circ dB_{t_2}^{i_2} \cdots \circ dB_{t_k}^{i_k}\Big)^2\Big] 
\leq Ct^k \qquad (0 \leq t \leq 1). 
$$
In view of this bound, we obtain
$
\mathbb{P}\big( d_{\mathrm{CC}}(g, Y_t^g) \geq R\big) \leq CR^{-2}t. 
$
Taking a sufficiently large $R>0$ such that $C\|f\|_\infty^G tR^{-2}<\ve$, 
we conclude $|T_tf(g)|<2\ve$ for $g \in B_{2R}(\bm{1}_G)^c$. 
This implies that $T_t\big(C_\infty(G)\big) \subset C_\infty(G)$ for $0 \leq t \leq 1$. \qed

\subsection{The free case: a relation with rough path theory}

Consider the step-$r$ non-commutative tensor algebra
$T^{(r)}(\mathbb{R}^d)=\mathbb{R} 
\oplus \big( \bigoplus_{k=1}^r (\mathbb{R}^d)^{\otimes k}\big)$.
The tensor product on $T^{(r)}(\mathbb{R}^{d})$ is defined by
$$
(g_0, g_1, \dots, g_r) \otimes_r (h_0, h_1, \dots, h_r)
=\Big( g_0h_0, g_0h_1+g_1h_0, \dots, \sum_{k=0}^r g_k \otimes h_{r-k}\Big).
$$
An element $g=(g_0, g_1, \dots, g_r) \in T^{(r)}(\mathbb{R}^d)$ is occasionally 
written as $g=g_0+g_1+\cdots+g_r$.  
We define two subsets of $T^{(r)}(\mathbb{R}^d)$ by
$$
T_1^{(r)}(\mathbb{R}^d):=\{g \in T^{(r)}(\mathbb{R}^d) \, | \, g_0=1\}, \qquad
T_0^{(r)}(\mathbb{R}^d):=\{A \in T^{(r)}(\mathbb{R}^d) \, | \, A_0=0\},
$$
respectively. 
It is easy to see that $T_1^{(r)}(\mathbb{R}^d)$ is a Lie group under the tensor product $\otimes_r$. 
In fact, $\bm{1}=(1, 0, 0, \dots, 0)$ is the unit element of $T_1^{(r)}(\mathbb{R}^d)$
and the inverse element of $g \in T_1^{(r)}(\mathbb{R}^d)$ 
is given by $g^{-1}=\sum_{k=1}^r (-1)^k (g - \bm{1})^{\otimes_r k}$. 
The Lie bracket on
$T_0^{(r)}(\mathbb{R}^d)$ 
is defined by $[A, B]=A \otimes_r B - B \otimes_r A$ for $A, B \in T_0^{(r)}(\mathbb{R}^d)$.
Note that $T_0^{(r)}(\mathbb{R}^d)$ is the Lie algebra of the Lie group $T_1^{(r)}(\mathbb{R}^d)$, 
that is, $T_0^{(r)}(\mathbb{R}^d)$ is the tangent space of $T_1^{(r)}(\mathbb{R}^d)$
at $\bm{1}$. 
The diffeomorphism $\exp : T_0^{(r)}(\mathbb{R}^d) \LA T_1^{(r)}(\mathbb{R}^d)$ is defined by
$$
\exp(A):=1+\sum_{k=1}^r \frac{1}{k!}A^{\otimes_r k} \qquad \big(A \in T_0^{(r)}(\mathbb{R}^d)\big).
$$
Let $\{{\bf{e}}_1, {\bf{e}}_2, \dots, {\bf{e}}_d\}$ be the standard basis of $\mathbb{R}^d$. 
We introduce a discrete subgroup $\g^{(r)}(\mathbb{Z}^d) \subset T_0^{(r)}(\mathbb{R}^d)$ 
by the set of $\mathbb{Z}$-linear combinations of
${\bf{e}}_1, {\bf{e}}_2, \dots, {\bf{e}}_d$ 
together with $[{\bf{e}}_{i_1}, [{\bf{e}}_{i_2}, \cdots, [{\bf{e}}_{i_{k-1}}, {\bf{e}}_{i_k}] \cdots]]$ for
$i_1, i_2, \dots, i_k=1, 2, \dots, d$ and $k=2, 3, \dots, r$.

We now set $\Gamma=\mathbb{G}^{(r)}(\mathbb{Z}^d):=\exp\big( \g^{(r)}(\mathbb{Z}^d)\big).$
We also define $\g^{(r)}(\mathbb{R}^d)$ and $\mathbb{G}^{(r)}(\mathbb{R}^d)$ analogously. 
Then we see that $\big(\mathbb{G}^{(r)}(\mathbb{R}^d), \otimes_r\big)$
is the nilpotent Lie group in which $\Gamma$ is included as its cocompact lattice and 
the corresponding limit group coincides with $\big(\mathbb{G}^{(r)}(\mathbb{R}^d), \otimes_r\big)$
itself. 
We call $\big(\mathbb{G}^{(r)}(\mathbb{R}^d), \otimes_r\big)$ the {\it free nilpotent Lie group of step $r$} and 
$\big(\g^{(r)}(\mathbb{R}^d), [\cdot, \cdot]\big)$ the {\it free nilpotent Lie algebra of step $r$}.
Let $\g^{(1)}=\mathbb{R}^d$ and 
$\g^{(k)}=[\mathbb{R}^d, [\mathbb{R}^d, \cdots, [\mathbb{R}^d, \mathbb{R}^d] \cdots ]]$
($k$-times) for $k=2, 3, \dots, r$. 
Then we see 
that the Lie algebra $\g^{(r)}(\mathbb{R}^d)$ is decomposed into 
$\g^{(1)} \oplus \g^{(2)} \oplus \cdots \oplus \g^{(r)}$. 
The free nilpotent Lie group $\mathbb{G}^{(r)}(\mathbb{R}^d)$ 
is highly related to rough path theory, as is seen below (cf.~Friz--Victoir \cite{FV}). 
Let $(B_t)_{0 \leq t \leq 1}=(B_t^1, B_t^2, \, \dots, B_t^d)_{0 \leq t \leq 1}$
be a $d$-dimensional standard Brownian motion. 
We give the following two remarks.

\vspace{1mm}
\noindent
(1) Consider the case $r=2$. Then a $T_1^{(2)}(\mathbb{R}^d)$-valued path 
$({\bf{B}}_t)_{0 \leq t \leq 1}$ defined by 
$$
\begin{aligned}
{\bf{B}}_t&:=\exp\Big( \sum_{i=1}^{d} B_t^i {\bf{e}}_i + 
\sum_{1 \leq i<j \leq d}
\Big(\frac{1}{2}\int_0^t B_s^i  \circ dB_s^j - B_s^j \circ dB_s^i\Big) 
{\bf{e}}_i \otimes {\bf{e}}_j
\Big)\\
&=1+\sum_{i=1}^{d}B_t^i {\bf{e}}_i +
\sum_{i, j=1}^d 
\Big(\int_0^t \int_0^s \circ dB_u^i \circ dB_s^j \Big) {\bf e}_i \otimes {\bf e}_j
 \qquad (0 \leq t \leq 1)
\end{aligned}
$$
is regarded as a $\mathbb{G}^{(2)}(\mathbb{R}^d)$-valued path with probability one.
We call it {\it Stratonovich enhanced Brownian motion} or 
{\it standard Brownian rough path},  
which is a canonical lift of 
a sample path of the $d$-dimensional Brownian motion. We usually identify  
standard Brownian rough path $({\bf{B}}_t)_{0 \leq t \leq 1}$ with 
its increment $({\bf B}_{s, t}):=({\bf B}_s^{-1} \otimes_2 {\bf B}_t)_{0 \leq s \leq t \leq 1}$.

 \vspace{1mm}
\noindent
(2) Consider the case $r \geq 3$. 
The $T_1^{(r)}(\mathbb{R}^d)$-valued path $({\bf{B}}_t)_{0 \leq t \leq 1}$ defined by 
$$
{\bf{B}}_t:=1+\sum_{k=1}^r \sum_{i_1, i_2, \dots, i_k \in \{1, 2, \dots, d\}}
\Big(\int_{\Delta^{(k)}[0, t]}
\circ dB_{t_1}^{i_1} \circ dB_{t_2}^{i_2} \cdots 
\circ dB_{t_k}^{i_k}\Big) {\bf{e}}_{i_1} \otimes {\bf{e}}_{i_2} \otimes \cdots \otimes {\bf{e}}_{i_k}
$$
for $0 \leq t \leq 1$, is regarded as a $\mathbb{G}^{(r)}(\mathbb{R}^d)$-valued path with probability one, 
analogously in (1).  
Note that this path $({\bf{B}}_t)_{0 \leq t \leq 1}$ is nothing but the {\it Lyons extension} (or {\it lift})
of Stratonovich enhanced Brownian motion introduced in (1) to $\mathbb{G}^{(r)}(\mathbb{R}^d)$.

Let $\Gamma=\mathbb{G}^{(r)}(\mathbb{Z}^d)$
and $X$ be a $\Gamma$-nilpotent covering graph. Then
we see that $X$ is realized into the free nilpotent Lie group 
$G=\mathbb{G}^{(r)}(\mathbb{R}^d)$ through the modified harmonic realization $\Phi_0 : X \LA G$, 
because $\Gamma$ is a cocompact lattice in $G$. 
Then Theorem \ref{limiting diffusion} reads in terms of rough path theory.
Precisely speaking, the $\mathbb{G}^{(r)}(\mathbb{R}^d)$-valued diffusion process 
$(Y_t)_{0 \leq t \leq 1}$
which solves  (\ref{SDE}) is represented as 
the Lyons extension of the so-called {\it distorted Brownian rough path} of order $r$.

\begin{co}\label{RP-relation}
Let $\{V_1, V_2, \dots, V_d\}$ be an orthonormal basis of $\g^{(1)}$
with respect to the Albanese metric $g_0$. 
We write
$$
\beta(\Phi_0)=\sum_{1 \leq i < j \leq d} \beta(\Phi_0)^{ij}[V_i, V_j] \in \g^{(2)},
$$
where we note that 
$\{[V_i, V_j] \, : \, 1 \leq i < j \leq d\} \subset \g^{(2)}$ forms a basis of $\g^{(2)}$. 
Let $\overline{\bm{\beta}}(\Phi_0)
=\big(\overline{\bm{\beta}}(\Phi_0)^{ij}\big)_{i, j =1}^d$ 
be an anti-symmetric matrix defined by 
$$
\overline{\bm{\beta}}(\Phi_0)^{ij}:=
\begin{cases}
\beta(\Phi_0)^{ij}  & (1 \leq i < j \leq d), \\
-\beta(\Phi_0)^{ij}  & (1 \leq j < i \leq d), \\
0  & (i=j). 
\end{cases}
$$
Then the $\mathbb{G}^{(r)}(\mathbb{R}^d)$-valued diffusion process 
$(Y_t)_{0 \leq t \leq 1}$ 
coincides with the Lyons extension of the distorted Brownian rough path
$$
\ol{\bm{\mathrm{B}}}_{t}
=1+\ol{\bm{\mathrm{B}}}_{t}^1+\ol{\bm{\mathrm{B}}}_{t}^2
\in \mathbb{G}^{(2)}(\mathbb{R}^d) \qquad (0 \leq  s \leq t \leq 1)
$$
of order $r$, where
$$
\ol{\bm{\mathrm{B}}}_{t}^1:=\sum_{i=1}^{d}B_t^i V_i \in \mathbb{R}^d, \quad
\ol{\bm{\mathrm{B}}}_{t}^2:=\int_0^t\int_0^s \circ dB_{u} \otimes \circ dB_s
+ t \overline{\bm{\beta}}(\Phi_0) \in \mathbb{R}^d \otimes \mathbb{R}^d. 
$$
\end{co}

\section{Examples}

We discuss several examples of the modified standard realizations associated with 
non-symmetric random walks on nilpotent covering graphs. 
It goes without saying that the most typical but non-trivial example of nilpotent 
Lie groups of step 2 is the {\it {\rm 3}-dimensional Heisenberg group} defined by
$$
G=\mathbb{H}^3(\mathbb{R}):=\Bigg\{ 
\begin{bmatrix}
1 & x & z \\
0 & 1 & y \\
0 & 0 & 1
\end{bmatrix} \, \Bigg| \, x, y, z \in \mathbb{R}\Bigg\}
=(\mathbb{R}^3, \star),
$$
where the product $\star$ on $\mathbb{R}^3$ is given by
$$
(x, y, z)\star (x', y', z')=(x+x', y+y', z+z'+xy').
$$
This Lie group naturally appears in a lot of parts of mathematics
including Fourier analysis, geometry, topology and so on. 
First of all, we give a quick review of the basics of $G=\mathbb{H}^3(\mathbb{R})$.
Let $\Gamma=\mathbb{H}^3(\mathbb{Z})$ 
be the 3-dimensional discrete Heisenberg group.
Then, $G=\mathbb{H}^3(\mathbb{R})$ is the corresponding 
connected and simply connected nilpotent Lie group of step 2
such that $\Gamma$ is isomorphic to a cocompact lattice in $G$.
Furthermore, the corresponding Lie algebra $\g$ is given by
$$
\g=\Bigg\{ 
\begin{bmatrix}
0 & x & z \\
0 & 0 & y \\
0 & 0 & 0
\end{bmatrix} \, \Bigg| \, x, y, z \in \mathbb{R}\Bigg\}.
$$
Let $\{X_1, X_2, X_3\}$ be the standard basis of $\g$, that is, 
$$
X_1:=\begin{bmatrix}
0 & 1 & 0 \\
0 & 0 & 0 \\
0 & 0 & 0
\end{bmatrix}, \quad X_2:=\begin{bmatrix}
0 & 0 & 0 \\
0 & 0 & 1 \\
0 & 0 & 0
\end{bmatrix}, \quad X_3:=\begin{bmatrix}
0 & 0 & 1 \\
0 & 0 & 0 \\
0 & 0 & 0
\end{bmatrix}.
$$
We then see that the Lie algebra $\g$ is decomposed as
$\g=\g^{(1)} \oplus \g^{(2)}$, where
$\g^{(1)}:=\Span\{X_1, X_2\}$ and $\g^{(2)}:=\Span\{X_3\}$,
due to the algebraic relations $[X_1, X_2]=X_3$ and $[X_1, X_3]=[X_2, X_3]
=\bm{0}_{\g}$ under the matrix bracket $[X, Y]:=XY-YX$ for $X, Y \in \g$. 

\subsection{The 3D Heisenberg triangular lattice}
\begin{figure}[htbp]
\begin{center}
\includegraphics[width=10cm]{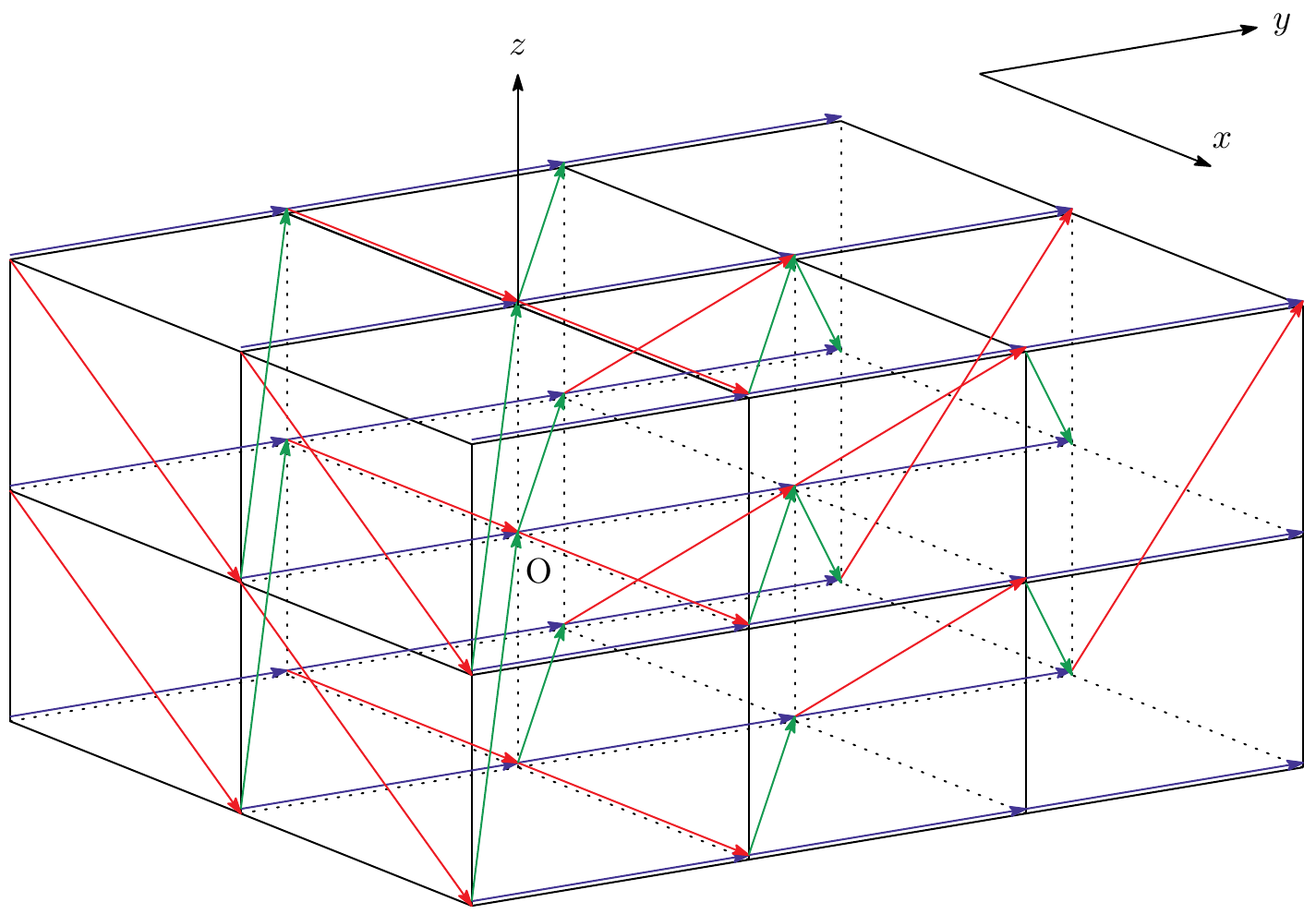}
\caption{A part of the 3-dimensional Heisenberg triangular lattice}
\label{Figure2}
\end{center}
\vspace{-0.5cm}
\end{figure}

\begin{figure}[htbp]
\begin{center}
\includegraphics[width=12cm]{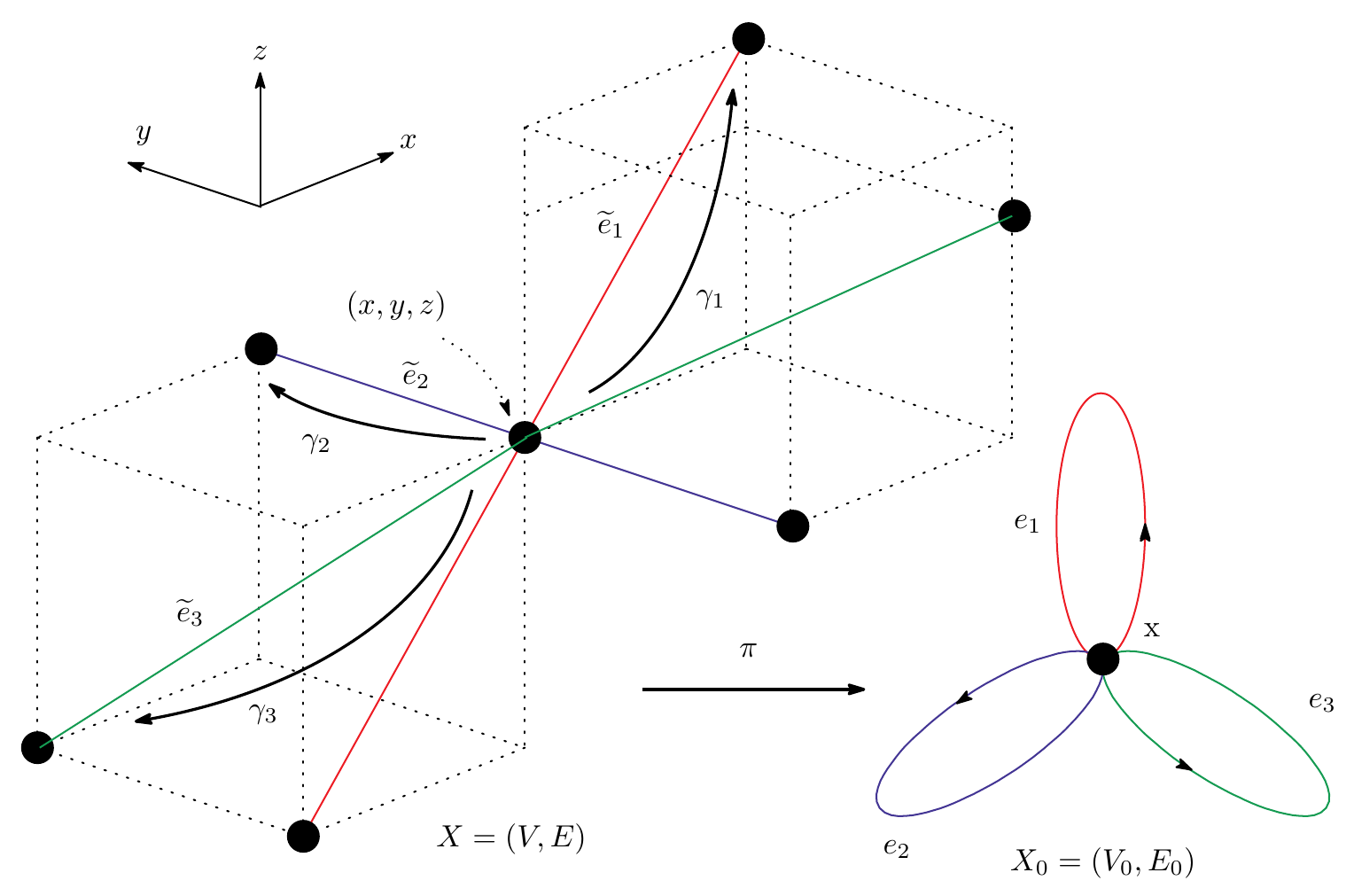}
\caption{The quotient $X_0=(V_0, E_0)=\Gamma \backslash X$ 
and the nearest neighbor vertices of $(x, y, z)$}
\label{Figure-nearest}
\end{center}
\vspace{-0.5cm}
\end{figure}

Let $\Gamma$ be generated by $\gamma_1=(1, 0, 0)$, $\gamma_2=(0, 1, 0)$
and $\gamma_3=(-1, 1, 0)$. 
We consider the {\it Cayley graph} $X=(V, E)$ of $\Gamma$ with the {\it generating set} 
$\mathcal{S}:=\{\gamma_1, \gamma_2, \gamma_3, 
\gamma_1^{-1}, \gamma_2^{-1}, \gamma_3^{-1}\}.$
Namely, $V=\mathbb{Z}^3$ and 
$E=\{(g, h) \in V \times V \, | \, h \cdot g^{-1} \in \mathcal{S}\}$ (see Figure \ref{Figure2}). 
If $e \in E$ is represented as $(g, h)$ for some $g, h \in V$, 
then its inverse edge $\ol{e}$ is equal to $(h, g)$. 
Moreover, the left action $\Gamma$ on the Cayley graph $X$ is given by
\begin{align}
\gamma_1g &= (x+1, y, z+y), & \gamma_2g&=(x, y+1, z),  & \gamma_3g&=(x-1, y+1, z-y), \nn\\
\gamma_1^{-1}g &= (x-1, y, z-y),  & \gamma_2^{-1}g&=(x, y-1, z), & \gamma_3^{-1}g&=(x+1, y-1, z+y-1),\nn
\end{align}
for $g=(x, y, z) \in G$. 
In view of the algebraic relation $\gamma_3 \star \gamma_1=\gamma_2$, 
we may call this Cayley graph $X$ a {\it {\rm 3}-dimensional Heisenberg triangular lattice}. 
The quotient graph of $X$ by the action $\Gamma$
 is the 3-bouquet graph $X_0=(V_0, E_0)$, where
$V_0=\{\x\}$ and $E_0=\{e_1, e_2, e_3\} \cup \{\ol{e}_1, \ol{e}_2, \ol{e}_3\}$ 
(see Figure \ref{Figure-nearest}).

Now we define a non-symmetric random walk on $X$. 
We introduce a transition probability $p : E \LA (0, 1]$ on $X$ by setting
\begin{align}
p\big((g, \gamma_1 g)\big)&:=\xi, & 
p\big((g, \gamma_2 g)\big)&:=\eta', & 
p\big((g, \gamma_3 g)\big)&:=\zeta, \nn\\
p\big((g, \gamma_1^{-1} g)\big)&:=\xi', & 
p\big((g, \gamma_2^{-1} g)\big)&:=\eta, & 
p\big((g, \gamma_3^{-1} g)\big)&:=\zeta',\nn
\end{align}
where $\xi, \xi', \eta, \eta', \zeta, \zeta' > 0$, $\xi+\xi'+\eta+\eta'+\zeta+\zeta'=1$ 
and 
\begin{equation}\label{non-symmetry}
\xi-\xi'=\eta-\eta'=\zeta-\zeta'=:\ve \geq 0.
\end{equation}
In what follows, we write 
$$
\hat{\xi}:=\xi+\xi',\quad \check{\xi}:=\xi-\xi', \quad  \hat{\eta}:=\eta+\eta', \quad \check{\eta}:=\eta-\eta', \quad 
\hat{\zeta}:=\zeta+\zeta', \quad \check{\zeta}:=\zeta-\zeta'
$$
for brevity. 
The invariant measure on $V_0=\{\x\}$ is given by $m(\x)=1$. 
The quantity $\ve$ in (\ref{non-symmetry})
indicates the intensity of the non-symmetry of this random walk and 
it is clear that the random walk is $m$-symmetric if and only if $\ve=0$.

The first homology group of $X_0$ 
is given by $\h_1(X_0, \mathbb{R})=\{[e_1], [e_2], [e_3]\}$.
Since $X_0$ is a bouquet graph, 
the difference operator $d : C^0(X_0, \mathbb{R}) \LA C^1(X_0, \mathbb{R})$
is the zero-map. Then we have 
$\h^1(X_0, \mathbb{R}) \cong 
\big(\mathcal{H}^1(X_0), \La \cdot, \cdot \Ra_p \big)=C^1(X_0, \mathbb{R}).$
Moreover, we obtain
\begin{equation}\label{h-direction}
\gamma_p=\sum_{e \in E_0}p(e)[e]=\ve\big([e_1]-[e_2]+[e_3]\big) \in \h_1(X_0, \mathbb{R})
\end{equation}
by definition. 
The canonical surjective linear map $\rho_{\mathbb{R}}: \h_1(X_0, \mathbb{R}) \LA \g^{(1)}$ is
given by
$$
\rho_{\mathbb{R}}([e_1])=X_1, \quad \rho_{\mathbb{R}}([e_2])=X_2, \quad \rho_{\mathbb{R}}([e_3])=X_2-X_1.
$$
Then we easily see that 
$\rho_{\mathbb{R}}(\gamma_p)=\bm{0}_{\g}$.
We introduce a basis $\{u_1, u_2\}$ in $\Hom(\g^{(1)}, \mathbb{R})$ by
$$
u_1(X)=x, \quad u_2(X)= y \quad \big(X=xX_1+yX_2 \in \g^{(1)}, \, x, y \in \mathbb{Z}\big).
$$
It should be noted that  $\{u_1, u_2\}$ is the dual basis of $\{X_1, X_2\}$ in $\g^{(1)}$. 
We write $\{\omega_1, \omega_2, \omega_3\} 
\subset \big(\h^1(X_0, \mathbb{R}), \La \cdot, \cdot \Ra_p\big)$
for the dual basis of $\{[e_1], [e_2], [e_3]\} \subset \h_1(X_0, \mathbb{R})$. 
By direct computation, we obtain
\begin{align}\label{inner-products-1}
\La \omega_1, \omega_1 \Ra_p &= \hat{\xi}-\check{\xi}^2=\hat{\xi}-\ve^2, 
& \La \omega_1, \omega_2 \Ra_p &=\check{\xi}\check{\eta}=\ve^2, \nonumber\\
\La \omega_2, \omega_2 \Ra_p &= \hat{\eta}-\check{\eta}^2=\hat{\eta}-\ve^2, 
& \La \omega_2, \omega_3 \Ra_p &=\check{\eta}\check{\zeta}=\ve^2,\\
\La \omega_3, \omega_3 \Ra_p &= \hat{\zeta}-\check{\zeta}^2=\hat{\zeta}-\ve^2, 
& \La \omega_1, \omega_3 \Ra_p &=-\check{\xi}\check{\zeta}=-\ve^2. \nonumber
\end{align}
We know that
$
u_1={}^t \rho_{\mathbb{R}}(u_1)=\omega_1 - \omega_3, 
\, u_2={}^t \rho_{\mathbb{R}}(u_2)=\omega_2 + \omega_3
$
form a $\mathbb{Z}$-basis in $\Hom(\g^{(1)}, \mathbb{R})$
by noting that $\Hom(\g^{(1)}, \mathbb{R})$ is regarded
as a 2-dimensional subspace of $\h^1(X_0, \mathbb{R})$ 
through the injective map ${}^t\rho_{\mathbb{R}}$.
It follows from  (\ref{inner-products-1}) that
\begin{equation}\label{inner-products-2}
\La u_1, u_1\Ra_p=\hat{\xi}+\hat{\zeta}, \quad \La u_1, u_2 \Ra_p=-\hat{\zeta}, 
\quad \La u_2, u_2 \Ra_p = \hat{\eta}+\hat{\zeta}.
\end{equation}
Then the volume of the Albanese torus 
 is computed as
$$
\mathrm{vol}(\mathrm{Alb}^{\Gamma})^{-1}:=\sqrt{\det \big( \La u_i, u_j \Ra_p \big)_{i, j=1}^2}
=(\hat{\xi}\hat{\eta}+\hat{\eta}\hat{\zeta}+\hat{\zeta}\hat{\xi})^{1/2}.
$$
Moreover, the Albanese metric $g_0$ on $\g^{(1)}$ is given by the following:
$$
\begin{aligned}
\la X_1, X_1 \ra_{g_0}
&=\frac{\hat{\eta}+\hat{\zeta}}{\hat{\xi}\hat{\eta}+\hat{\eta}\hat{\zeta}+\hat{\zeta}\hat{\xi}}
=(\hat{\eta}+\hat{\zeta})\mathrm{vol}(\mathrm{Alb}^{\Gamma})^2, \\
\la X_1, X_2 \ra_{g_0}
&=\frac{\hat{\zeta}}{\hat{\xi}\hat{\eta}+\hat{\eta}\hat{\zeta}+\hat{\zeta}\hat{\xi}}
=\hat{\zeta}\mathrm{vol}(\mathrm{Alb}^{\Gamma})^2, \\
\la X_2, X_2 \ra_{g_0}
&=\frac{\hat{\xi}+\hat{\zeta}}{\hat{\xi}\hat{\eta}+\hat{\eta}\hat{\zeta}+\hat{\zeta}\hat{\xi}}
=(\hat{\xi}+\hat{\zeta})\mathrm{vol}(\mathrm{Alb}^{\Gamma})^2.
\end{aligned}
$$

We are now in a position to determine the modified standard realization
$\Phi_0 : X \LA G$. 
Let $\widetilde{e}_i \, (i=1, 2, 3)$ be a lift of $e_i \in E_0$ to $X$ 
and put $\Phi_0\big( o(\widetilde{e}_i)\big)=\bm{1}_G=(0, 0, 0)$. 
Then we easily see that the realization satisfying
$
\Phi_0\big( t(\widetilde{e}_1)\big)=\gamma_1, \, \Phi_0\big( t(\widetilde{e}_2)\big)=\gamma_2$ and 
$\Phi_0\big( t(\widetilde{e}_3)\big)=\gamma_3$
is the modified harmonic realization. Let $\{v_1, v_2\}$ be the Gram--Schmidt orthonormalization of $\{u_1, u_2\}$,
and $\{V_1, V_2\}$ be the dual basis of $\{v_1, v_2\}$ in $\g^{(1)}$. 
We put $V_3:=[V_1, V_2]=V_1V_2-V_2V_1.$
We then have
$$
v_1=(\hat{\xi}+\hat{\zeta})^{-1/2}u_1, \quad v_2=(\hat{\xi}+\hat{\zeta})^{1/2}
\mathrm{vol}(\mathrm{Alb}^{\Gamma})\Big( \frac{\hat{\zeta}}{\hat{\xi}+\hat{\zeta}}u_1+u_2\Big)
$$
by (\ref{inner-products-2})
and hence we obtain
\begin{align}
V_1 &=(\hat{\xi}+\hat{\zeta})^{1/2}X_1 
                     - \hat{\zeta}(\hat{\xi}+\hat{\zeta})^{-1/2}X_2, \nonumber  \\
V_2 &=(\hat{\xi}+\hat{\zeta})^{-1/2}\mathrm{vol}(\mathrm{Alb}^{\Gamma})^{-1}X_2, \nn\\
V_3 &=\mathrm{vol}(\mathrm{Alb}^{\Gamma})^{-1}X_3. \nonumber
\end{align}
Finally, $\beta(\Phi_0) \in \g^{(2)}$ and the infinitesimal generator 
$\A$ in Theorem \ref{CLT1} are calculated as
$$
\beta(\Phi_0)=\frac{\ve}{2}\mathrm{vol}(\mathrm{Alb}^{\Gamma})V_3, \qquad
\A=-\frac{1}{2}(V_1^2+V_2^2)-\frac{\ve}{2}\mathrm{vol}(\mathrm{Alb}^{\Gamma})V_3,
$$
respectively. 

\subsection{The 3D Heisenberg dice lattice}

\begin{figure}[htbp]
\begin{center}
\includegraphics{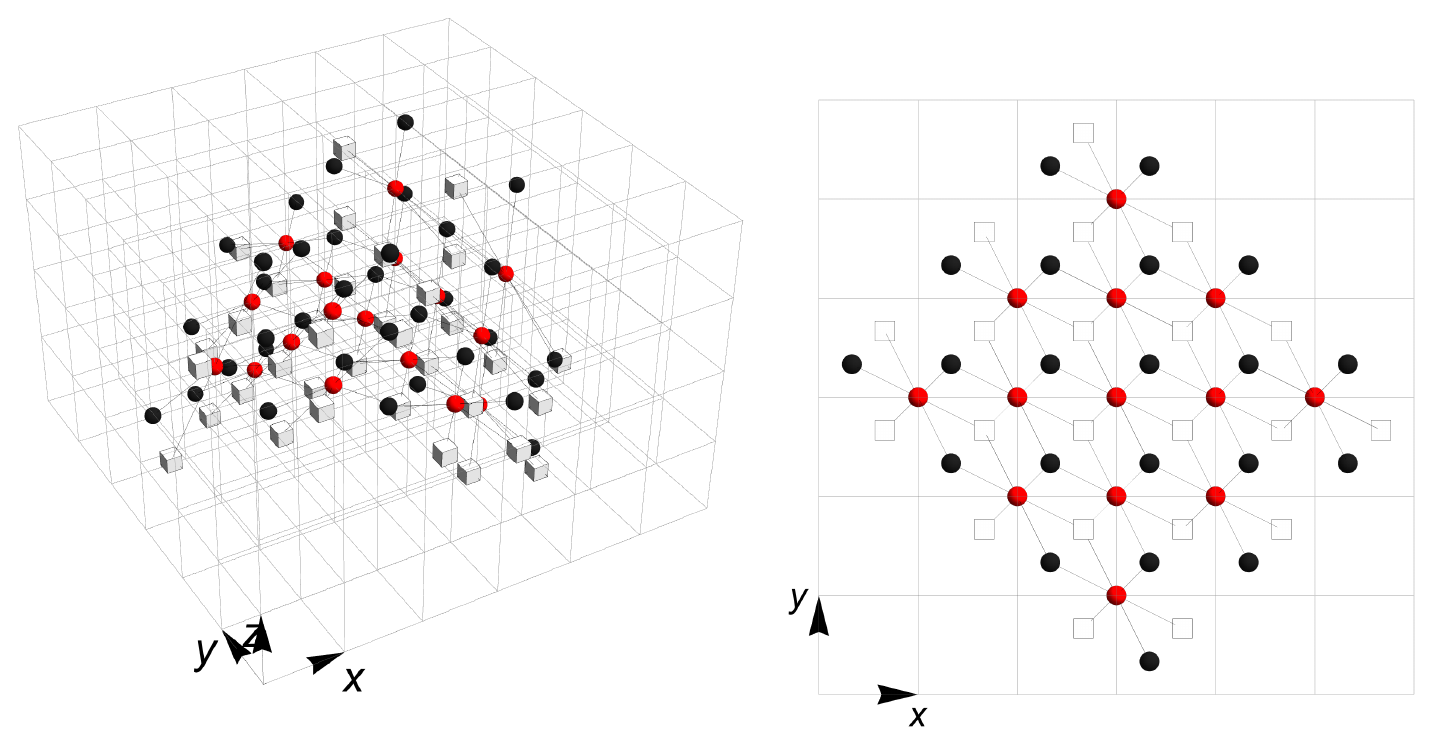}
\caption{A part of 3-dimensional Heisenberg dice lattice and 
the projection of it on the $xy$-plane}
\label{dices}
\end{center}
\vspace{-0.5cm}
\end{figure}

As another example of nilpotent covering graphs, 
we introduce the {\it {\rm 3}-dimensional Heisenberg dice lattice}.
This graph is defined by a covering graph of a finite graph 
consisting of three vertices with a covering transformation group 
$\Gamma=\mathbb{H}^3(\mathbb{Z})$ (see Figure \ref{dices}). 
We emphasize that it is regarded as an extension of the {\it dice graph} 
discussed in Namba \cite{Namba} to the nilpotent case. 

\begin{figure}[htbp]
\begin{center}
\includegraphics[width=7cm]{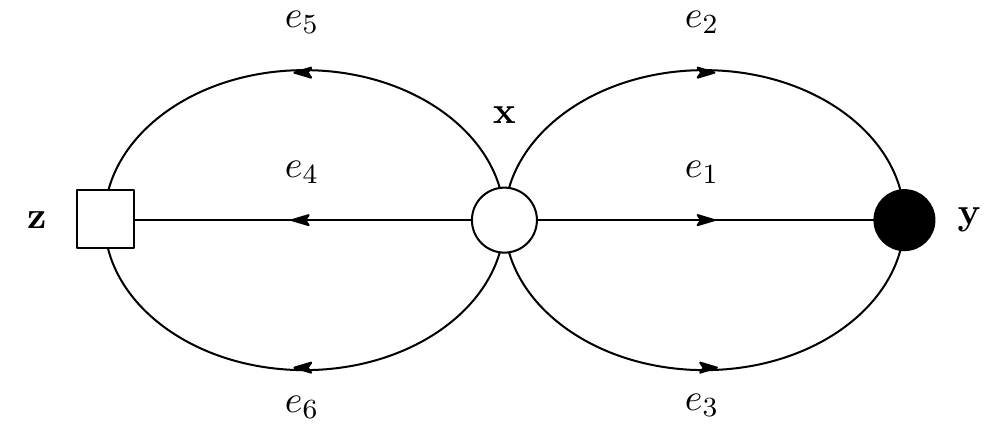}
\caption{The quotient $X_0=(V_0, E_0)$ of the 3D-Heisenberg dice graph $X=(V, E)$}
\label{Figure-dice-quotient}
\end{center}
\vspace{-0.5cm}
\end{figure}

Suppose that $\Gamma=\mathbb{H}^3(\mathbb{Z})$ is generated by two elements
$\gamma_1=(1, 0, 0)$ and $\gamma_2=(0, 1, 0)$. We also set two elements
$\bg_1:=(1/3, 1/3, 1/3)$, $\bg_2:=(-1/3, -1/3, -1/3)$ in 
$G=\mathbb{H}^3(\mathbb{R})$. We put
$$
\begin{aligned}
V_{1}&:= \big\{ g=\gamma_{i_1}^{\ve_1} \star \dots \star \gamma_{i_\ell}^{\ve_\ell} 
                                \star \bm{1}_G \, \big| \, i_k \in \{1, 2\}, \, \ve_k=\pm1 \, 
                                (1 \leq k \leq \ell), \, \ell \in \mathbb{N} \cup \{0\}\big\},\\
V_{2}&:= \big\{ g=\gamma_{i_1}^{\ve_1} \star\dots \star \gamma_{i_\ell}^{\ve_\ell} 
                                \star \bg_1 \, \big| \, i_k \in \{1, 2\}, \, \ve_k=\pm1 \, 
                                (1 \leq k \leq \ell), \, \ell \in \mathbb{N} \cup \{0\}\big\}, \\
V_{3}&:=\big\{ g=\gamma_{i_1}^{\ve_1} \star \dots \star  \gamma_{i_\ell}^{\ve_\ell} 
                                \star \bg_2 \, \big| \, i_k \in \{1, 2\}, \, \ve_k=\pm1 \, 
                                (1 \leq k \leq \ell), \, \ell \in \mathbb{N} \cup \{0\}\big\}.
\end{aligned}
$$
We consider an $\mathbb{H}^3(\mathbb{Z})$-nilpotent covering graph $X=(V, E)$ defined by
$V=V_1 \sqcup V_2 \sqcup V_3$ and $E=E_1 \sqcup E_2$, where
$$
\begin{aligned}
E_1 &:= \big\{ (g, h) \in V_1 \times V_2\, | \, g^{-1} \star h=\bg_1, \,
                                  \gamma_1^{-1} \star \bg_1, \, \gamma_2^{-1}  \star \bg_1 \big\}, \\
E_2 &:= \big\{ (g, h) \in V_1 \times V_3\, | \, g^{-1} \star h= \bg_2, \,
                                  \gamma_1  \star \bg_2, \, \gamma_2 \star \bg_2 \big\}.
\end{aligned}
$$
We note that  $X$ is invariant under the actions $\gamma_1$ and $\gamma_2$. 
Its quotient graph $X_0=(V_0, E_0)=\Gamma \backslash X$ is given by $V_0=\{\x, \y, \z\}$ 
and $E_0=\{e_i, \ol{e}_i \, | \, 1 \leq i \leq 6\}$ 
(cf. Figure \ref{Figure-dice-quotient}).

From now on we define a non-symmetric random walk on $X$. 
We define the transition probability $p : E \LA (0, 1]$ by
\begin{align}
p\big((g, g \star \bg_1)\big)&=\xi, & 
p\big((g, g \star \gamma_1^{-1}  \star \bg_1)\big)&=\eta, & 
p\big((g, g \star \gamma_2^{-1}  \star \bg_1)\big)&=\zeta, \nn\\
p\big((g, g \star \bg_2)\big)&=\zeta, & 
p\big((g, g \star \gamma_1 \star \bg_2)\big)&=\eta, &
p\big((g, g \star \gamma_2 \star \bg_2)\big)&=\xi, \nn\\
p\big(\ol{(g, g \star \bg_1)}\big)&=\gamma, &
p\big(\ol{(g, g \star \gamma_1^{-1}  \star \bg_1)}\big)&=\beta, &
p\big(\ol{(g, g \star \gamma_2^{-1}  \star \bg_1)}\big)&=\alpha, \nn\\
p\big(\ol{(g, g \star \bg_2)}\big)&=\alpha, &
p\big(\ol{(g, g \star \gamma_1  \star \bg_2)}\big)&=\beta, &
p\big(\ol{(g, g \star \gamma_2  \star \bg_2)}\big)&=\gamma, \nn
\end{align}
for every $g \in V_1$, where $\xi, \eta, \zeta, \alpha, \beta, \gamma>0$,
$2(\xi+\eta+\zeta)=1$ and $\alpha+\beta+\gamma=1$.
The invariant measure $m : V_0=\{\x, \y, \z\} \LA (0, 1]$ is given by
$m(\x)=1/2$ and $m(\y)=m(\z)=1/4$.
Note that this random walk is ($m$-)symmetric if and only if
$\alpha=2\zeta, \, \beta=2\eta$ and $\gamma=2\xi$.

The first homology group $\h_1(X_0, \mathbb{R})$ is spanned 
by the four 1-cycles
$$
[c_1]:=[e_1 * \ol{e}_2], \quad [c_2]:=[e_1*\ol{e}_3], \quad 
[c_3]:=[e_4*\ol{e}_5], \quad [c_4]:=[e_4*\ol{e}_6].
$$
Then the homological direction is calculated as
$$
\gamma_p=\frac{\beta-2\eta}{4}[c_1]+\frac{\alpha-2\zeta}{4}[c_2]
                   +\frac{\beta-2\eta}{4}[c_3]+\frac{\gamma-2\xi}{4}[c_4].
$$
The canonical surjective linear map
$\rho_{\mathbb{R}} : \h_1(X_0, \mathbb{R}) \LA \g^{(1)}$ is given by
$$
\rho_{\mathbb{R}}([c_1])=X_1, \quad 
\rho_{\mathbb{R}}([c_2])=X_2, \quad 
\rho_{\mathbb{R}}([c_3])=-X_1, \quad 
\rho_{\mathbb{R}}([c_4])=-X_2.
$$
Then we obtain
\begin{equation}\label{dice-2}
\rho_{\mathbb{R}}(\gamma_p)=\frac{(\alpha-\gamma)-2(\zeta-\xi)}{4}X_2.
\end{equation}

Let $\{u_1, u_2\} \subset \Hom(\g^{(1)}, \mathbb{R})$ 
be the dual basis of $\{X_1, X_2\} \subset \g^{(1)}$. 
We also denote by $\{\omega_1, \omega_2, \omega_3, \omega_4\} \subset
 \big(\h^1(X_0, \mathbb{R}), \La \cdot , \cdot \Ra_p\big)$ the dual basis
 of $\{[c_1], [c_2], [c_3], [c_4]\} \subset \h_1(X_0, \mathbb{R})$. 
 Namely, $\omega_i([c_j])=\delta_{ij}$ for $1 \leq i, j \leq 4$. 
Then the modified harmonicity (\ref{form-harmonicity}) yields
\begin{align}
\omega_1(e_1)&=\beta-\frac{\beta-2\eta}{4}, &
\omega_1(e_2)&=-(1-\beta)-\frac{\beta-2\eta}{4}, &
\omega_1(e_3)&=\beta-\frac{\beta-2\eta}{4}, \nn \\
\omega_1(e_4)&=-\frac{\beta-2\eta}{4}, &
\omega_1(e_5)&=-\frac{\beta-2\eta}{4}, &
\omega_1(e_6)&=-\frac{\beta-2\eta}{4}, \nn \\
\omega_2(e_1)&=\alpha-\frac{\alpha-2\zeta}{4}, &
\omega_2(e_2)&=\alpha-\frac{\alpha-2\zeta}{4}, &
\omega_2(e_3)&=-(1-\alpha)-\frac{\alpha-2\zeta}{4}, \nn \\
\omega_2(e_4)&=-\frac{\alpha-2\zeta}{4}, &
\omega_2(e_5)&=-\frac{\alpha-2\zeta}{4}, &
\omega_2(e_6)&=-\frac{\alpha-2\zeta}{4}, \nn \\
\omega_3(e_1)&=-\frac{\beta-2\eta}{4}, &
\omega_3(e_2)&=-\frac{\beta-2\eta}{4}, &
\omega_3(e_3)&=-\frac{\beta-2\eta}{4}, \nn \\
\omega_3(e_4)&=\beta-\frac{\beta-2\eta}{4}, &
\omega_3(e_5)&=-(1-\beta)-\frac{\beta-2\eta}{4}, &
\omega_3(e_6)&=\beta-\frac{\beta-2\eta}{4}, \nn \\
\omega_4(e_1)&=-\frac{\gamma-2\xi}{4}, &
\omega_4(e_2)&=-\frac{\gamma-2\xi}{4}, &
\omega_4(e_3)&=-\frac{\gamma-2\xi}{4}, \nn \\
\omega_4(e_4)&=\gamma-\frac{\gamma-2\xi}{4}, &
\omega_4(e_5)&=\gamma-\frac{\gamma-2\xi}{4}, &
\omega_4(e_6)&=-(1-\gamma)-\frac{\gamma-2\xi}{4}. \nn 
\end{align}
By direct computation, we have
\begin{align}\label{dice-3}
\La \omega_1, \omega_1 \Ra_p &= \frac{\beta+2\eta}{4} - \frac{(\beta+2\eta)^2}{8}, &
\La \omega_1, \omega_2 \Ra_p &= -\frac{(\alpha+2\zeta)(\beta+2\eta)}{8}, \nn \\
\La \omega_1, \omega_3 \Ra_p &= -\frac{(\beta - 2\eta)^2}{8}, &
\La \omega_1, \omega_4 \Ra_p &= -\frac{(\beta-2\eta)(\gamma-2\xi)}{8}, \nn \\
\La \omega_2, \omega_2 \Ra_p &= \frac{\alpha+2\zeta}{4} - \frac{(\alpha+2\zeta)^2}{8}, &
\La \omega_2, \omega_3 \Ra_p &= -\frac{(\alpha-2\zeta)(\beta-2\eta)}{8},  \\
\La \omega_2, \omega_4 \Ra_p &= -\frac{(\alpha - 2\zeta)(\gamma - 2\xi)}{8}, &
\La \omega_3, \omega_3 \Ra_p &= \frac{\beta+2\eta}{4} - \frac{(\beta+2\eta)^2}{8}, \nn \\
\La \omega_3, \omega_4 \Ra_p &= -\frac{(\beta+2\eta)(\gamma+2\xi)}{8}, &
\La \omega_4, \omega_4 \Ra_p &= \frac{\gamma+2\xi}{4} - \frac{(\gamma+2\xi)^2}{8}. \nn
\end{align}
Since the linear space $\Hom(\g^{(1)}, \mathbb{R})$ can be seen as a 2-dimensional subspace
of $\h^1(X_0, \mathbb{R})$ through the injection ${}^t \rho_{\mathbb{R}}$,
we see that 
$u_1={}^t \rho_{\mathbb{R}}(u_1)=\omega_1 - \omega_3$ and
$u_2={}^t \rho_{\mathbb{R}}(u_2)=\omega_2 - \omega_4$
form a $\mathbb{Z}$-basis in $\Hom(\g^{(1)}, \mathbb{R})$. We then obtain
$$
\La u_1, u_1 \Ra_p = \frac{\beta+2\eta-4\beta\eta}{2},\quad
\La u_1, u_2 \Ra_p = -\frac{\beta+2\eta-4\beta\eta}{4},
$$
$$
\La u_2, u_2 \Ra_p = \frac{(\beta+2\eta)(2-\beta-2\eta)+4\alpha \gamma+16\xi\zeta}{8}.
$$
by (\ref{dice-3}).
Thus the volume of the Albanese torus is computed as
$$
\mathrm{vol}(\mathrm{Alb}^{\Gamma})^{-1}=\frac{1}{4}\sqrt{(\beta+2\eta-4\beta\eta)
\big\{ (\beta+2\eta)-(\beta^2+4\eta^2)+4\alpha\gamma+16\xi\zeta\big\}}.
$$
Furthermore, the Albanese metric $g_0$ on $\g^{(1)}$ is given by
$$
\la X_1, X_1 \ra_{g_0} = \frac{(\beta+2\eta)(2-\beta-2\eta)+4\alpha \gamma+16\xi\zeta}{8}
\mathrm{vol}(\mathrm{Alb}^{\Gamma}),
$$
$$
\la X_1, X_2 \ra_{g_0} = \frac{\beta+2\eta-4\beta\eta}{4}
\mathrm{vol}(\mathrm{Alb}^{\Gamma}), \quad 
\la X_2, X_2 \ra_{g_0} = \frac{\beta+2\eta-4\beta\eta}{2}
\mathrm{vol}(\mathrm{Alb}^{\Gamma}).
$$

We now determine the modified standard realization 
$\Phi_0 : X \LA G=\mathbb{H}^3(\mathbb{R})$. 
Let $\widetilde{e}_i \, (i=1, 2, 3, 4, 5, 6)$ be a lift of $e_i \in E_0$ to $X$ 
and put $\Phi_0\big( o(\widetilde{e}_i)\big)=\bm{1}_G$. 
Then it follows from (\ref{m-harmonicity}) and (\ref{dice-2}) that
the $\Gamma$-equivariant realization $\Phi_0 : X \LA G$ satisfying
$$
\begin{aligned}
\Phi_0\big( t(\widetilde{e}_1)\big)&=\Big( \beta, 
                     \frac{(3\alpha+\gamma)+2(\zeta-\xi)}{4}, \kappa_1\Big), \\
\Phi_0\big( t(\widetilde{e}_2)\big)&=\Big( \beta-1, 
                     \frac{(3\alpha+\gamma)+2(\zeta-\xi)}{4}, 
                     \kappa_1 -\frac{(3\alpha+\gamma)+2(\zeta-\xi)}{4}\Big), \\
\Phi_0\big( t(\widetilde{e}_3)\big)&=\Big( \beta, 
                     \frac{(3\alpha+\gamma)+2(\zeta-\xi)}{4}-1, \kappa_1\Big), \\             
\Phi_0\big( t(\widetilde{e}_4)\big)&=\Big( -\beta, 
                     \frac{-(\alpha+3\gamma)+2(\zeta-\xi)}{4}, -\kappa_2\Big), \\
\Phi_0\big( t(\widetilde{e}_5)\big)&=\Big( 1-\beta, 
                     \frac{-(\alpha+3\gamma)+2(\zeta-\xi)}{4}, -\kappa_2
                     +\frac{-(\alpha+3\gamma)+2(\zeta-\xi)}{4}\Big), \\
\Phi_0\big( t(\widetilde{e}_6)\big)&=\Big( -\beta, 
                     \frac{-(\alpha+3\gamma)+2(\zeta-\xi)}{4}+1, -\kappa_2\Big)
\end{aligned}
$$
is the modified harmonic realization, where $\kappa_1, \kappa_2$ is two real parameters 
which indicates the ambiguity of the realization corresponding to $\g^{(2)}$. 
Let $\{v_1, v_2\}$ be the Gram--Schmidt orthonormalization 
of the basis $\{u_1, u_2\}$, that is, 
$$
v_1=\La u_1, u_1 \Ra_p^{-1/2}u_1, \quad 
v_2=\La u_1, u_1 \Ra_p^{1/2} \mathrm{vol}(\mathrm{Alb}^\Gamma)
               \Big( u_2 - \frac{\La u_1, u_2 \Ra_p}{\La u_1, u_1 \Ra_p}u_1 \Big),
$$
and $\{V_1, V_2\} \subset \g^{(1)}$ its dual basis. 
We write $V_3:=[V_1, V_2]=V_1V_2-V_2V_1$. 
Then we obtain
$$
v_1=\Big(\frac{\beta+2\eta-4\beta\eta}{2}\Big)^{-1/2}u_1,\quad
v_2=\Big(\frac{\beta+2\eta-4\beta\eta}{2}\Big)^{1/2} \mathrm{vol}(\mathrm{Alb}^\Gamma)
                 \Big( u_2 + \frac{1}{2}u_1\Big)
$$
by (\ref{dice-3}). Moreover, we have
$$
\begin{aligned}
V_1 &= \Big(\frac{\beta+2\eta-4\beta\eta}{2}\Big)^{1/2}X_1
                      - \frac{1}{2}\Big(\frac{\beta+2\eta-4\beta\eta}{2}\Big)^{1/2}X_2,\\
V_2 &= \Big(\frac{\beta+2\eta-4\beta\eta}{2}\Big)^{-1/2} 
                      \mathrm{vol}(\mathrm{Alb}^\Gamma)^{-1}X_2, \\
V_3 &= \mathrm{vol}(\mathrm{Alb}^\Gamma)^{-1}X_3.             
\end{aligned}
$$
Finally, we see that $\beta(\Phi_0) \in \g^{(2)}$ and the infinitesimal generator 
$\A$ are calculated as
$$
\begin{aligned}
\beta(\Phi_0) &= \sum_{i=1}^6 \big( \widetilde{m}(e_i) - \widetilde{m}(\ol{e}_i)\big) 
                          \Log \Big( d\Phi_0(\widetilde{e}_i) \cdot \exp\big(
                        -\rho_{\mathbb{R}}(\gamma_p)\big)\Big)\Big|_{\g^{(2)}}
     = \frac{\beta - 2\eta}{8}\mathrm{vol}(\mathrm{Alb}^\Gamma)V_3, \\
     \A&=-\frac{1}{2}(V_1^2+V_2^2)-\frac{\beta - 2\eta}{8}\mathrm{vol}(\mathrm{Alb}^\Gamma)V_3.
\end{aligned}
$$
We should observe that the coefficient of $\beta(\Phi_0)$ does not include
the parameters $\kappa_1$ and $\kappa_2$, though the realization $\Phi_0$ has the ambiguity 
of $\g^{(2)}$-components.

\begin{appendix}

\section{A comment on CLTs in the non-centered case}

As was already mentioned, the centered condition {\bf (C)} is crucial 
to establish the functional CLT (Theorem \ref{FCLT1}). 
We present a method to reduce 
the non-centered case $\rho_{\mathbb{R}}(\gamma_p) \neq \bm{0}_{\g}$
to the centered case by employing a 
measure-change technique based on Alexopoulos \cite{A2}.
See also Namba \cite{Namba} for this kind of technique in the case where $X$ is a crystal lattice. 
 
We consider a {\it positive} transition probability 
$p : E \LA (0, 1]$  to avoid several technical difficulties.
Then the random walk on $X$ associated with $p$ is automatically irreducible. 
Let  $\Phi_0 : X \LA G$ be the ($p$-)modified harmonic realization.
We define a function 
$F=F_x(\lambda) : V_0 \times \Hom(\g^{(1)}, \mathbb{R}) \LA \mathbb{R}$
by 
\begin{equation}\label{twist}
F_x(\lambda):=\sum_{e \in (E_0)_x}p(e)
\exp\Big( {}_{\Hom(\g^{(1)}, \mathbb{R})}\big\la \lambda, 
\Log\big( d\Phi_0(\widetilde{e})\big)\big|_{\g^{(1)}}\big\ra_{\g^{(1)}}\Big)
\end{equation}
for $x \in V_0$ and $\lambda \in \Hom(\g^{(1)}, \mathbb{R})$. 
Since the lemma below is obtained by
following the argument in \cite[Lemma 3.1]{Namba}, we omit the proof. 
\begin{lm}\label{lem-hessian}
For every $x \in V_0$, the function $F_x(\cdot) : \Hom(\g^{(1)}, \mathbb{R}) \LA (0, \infty)$
has a unique minimizer $\lambda_*=\lambda_*(x) \in \Hom(\g^{(1)}, \mathbb{R})$. 
\end{lm}

\vspace{2mm}
We now define a positive function $\p : E_0 \LA (0, 1]$ by
\begin{equation}\label{twist-p}
\p(e):=\frac{\exp\Big( {}_{\Hom(\g^{(1)}, \mathbb{R})}\big\la \lambda_*\big(o(e)\big), 
\Log\big( d\Phi_0(\widetilde{e})\big)\big|_{\g^{(1)}}\big\ra_{\g^{(1)}}\Big)}
{F_{o(e)}\big(\lambda_*\big(o(e)\big)\big)}p(e) \qquad (e \in E_0).
\end{equation}
It is straightforward to check that the function $\p$ also gives a positive
transition probability on $X_0$ and it yields an irreducible Markov chain
$(\Omega_x(X), \widehat{\mathbb{P}}_x, \{w_n^{(\p)}\}_{n=0}^\infty)$
with values in $X$. 
We then find a unique positive normalized
invariant measure $\m : V_0 \LA (0, 1]$ by applying the Perron-Frobenius theorem again.
We set $\widetilde{\m}(e):=\p(e)\m\big(o(e)\big)$ for $e \in E_0$. 
We also denote by $\p : E \LA (0, 1]$ and $\m : V \LA (0, 1]$
the $\Gamma$-invariant lifts of $\p : E_0 \LA (0,1]$ and
$\m : V_0 \LA (0, 1]$ to $X$, respectively. 
The Albanese metric on $\g^{(1)}$
associated with the transition probability $\p$ is denoted by $g_0^{(\p)}$.
We write  $\{V_1^{(\p)}, V_2^{(\p)}, \dots, V_{d_1}^{(\p)}\}$ for
an orthonormal basis of $(\g^{(1)}, g_0^{(\p)})$. 

Let $L_{(\p)} : C_\infty(X) \LA C_\infty(X)$  be the transition
operator associated with the transition probability $\p$.
By virtue of Lemma \ref{lem-hessian}, 
we have 
$$
\sum_{e \in (E_0)_x}p(e) 
\exp\Big( {}_{\Hom(\g^{(1)}, \mathbb{R})}\big\la \lambda_*,
 \Log \big(d\Phi_0(\widetilde{e})\big)\big|_{\g^{(1)}}\big\ra_{\g^{(1)}} \Big)
\Log \big(d\Phi_0(\widetilde{e})\big)\big|_{\g^{(1)}}=\bm{0}_{\g} \qquad (x \in V_0).
$$
Hence, we conclude
\begin{equation}\label{dF}
(L_{(\p)}-I)(\Log \Phi_0\big|_{\g^{(1)}})(x)=
\sum_{e \in E_x}\p(e) \Log \big(d\Phi_0(e)\big)\big|_{\g^{(1)}}=\bm{0}_\g
\qquad (x \in V).
\end{equation}
This means that the ($p$-)modified harmonic realization $\Phi_0 : X \LA G$ 
in the sense of (\ref{m-harmonicity})
is regarded as the ($\p$-)harmonic realization and
 $\rho_{\mathbb{R}}(\gamma_\p)=\bm{0}_\g$.

We fix a reference point $x_* \in V$ such that $\Phi_0(x_*)=\bm{1}_G$ and put 
$$
\xi_n^{(\p)}(c):=\Phi_0\big(w_n^{(\p)}(c)\big) \qquad \big(n \in \mathbb{N} \cup\{0\}, \, c \in \Omega_{x_*}(X)\big).
$$
This yields a $G$-valued random walk 
$(\Omega_{x_*}(X), \widehat{\mathbb{P}}_{x_*}, \{\xi_n^{(\p)}\}_{n=0}^\infty)$.
We define 
$$
\mathcal{Y}_{t_k}^{(n; \p)}(c):=\tau_{n^{-1/2}}\big( \xi_{nt_k}^{(\p)}(c)\big)
=\tau_{n^{-1/2}}\big( \Phi_0(w_{k}^{(\p)}(c))\big)
$$
for $k=0, 1, \dots, n, \, t_k \in \mathcal{D}_n$ and $c \in \Omega_{x_*}(X)$. 
We consider a $G$-valued stochastic process $(\mathcal{Y}_t^{(n;\, \p)})_{0 \leq t \leq 1} $
defined by the $d_{\mathrm{CC}}$-geodesic interpolation of
$\{\mathcal{Y}_{t_{k}}^{(n; \, \p)}\}_{k=0}^{n}$.
Let $(\widetilde{Y}_t)_{0 \leq t \leq 1}$ be the $G$-valued 
diffusion process which solves the SDE
$$
d\widetilde{Y}_t = \sum_{i=1}^{d_1} V_{i*}^{(\p)}(\widetilde{Y}_t) \circ dB_t^i 
+ \beta^{(\p)}(\Phi_0)_*(\widetilde{Y}_t) \, dt,\qquad \widetilde{Y}_0=\bm{1}_G,
$$
where 
$$
\beta^{(\p)}(\Phi_0)
:=\sum_{e \in E_0}\widetilde{\m}(e)
\Log \Big( \Phi_0\big(o(\widetilde{e})\big)^{-1} \cdot \Phi_0\big(t(\widetilde{e})\big)\Big)\Big|_{\frak{g}^{(2)}}.
$$
The following two theorems are CLTs for non-symmetric random walks 
associated with the changed transition probability $\p$. 
We remark that the proofs of these theorems below are done 
by combining the ones of Theorems \ref{CLT1} and \ref{FCLT1} with
the argument in \cite[Theorem 1.3]{Namba}.

\begin{tm}\label{Another CLT}
Let $P_\ve : C_\infty(G) \LA C_\infty(X)$ be the approximation operator defined 
by $P_\ve f(x):=f\big( \tau_\ve\big( \Phi_0(x)\big)\big)$ for $0 \leq \ve \leq 1$ and $x \in V$. 
Then we have, for $0 \leq s \leq t$ and $f \in C_\infty(G)$, 
\begin{equation}
\lim_{n \to \infty}\Big\| L_{(\p)}^{[nt]-[ns]}P_{n^{-1/2}} f - P_{n^{-1/2}} \e^{-(t-s)\A_{(\p)}}f\Big\|_\infty^X=0,
\end{equation}
where $(\e^{-t\A_{(\p)}})_{t \geq 0}$ is the $C_0$-semigroup with the infinitesimal generator 
$\A_{(\p)}$ on $C_0^\infty(G)$ defined by
\begin{equation}
\A_{(\p)}:=-\frac{1}{2}\sum_{i=1}^{d_1}(V_{i*}^{(\p)})^2 - \beta^{(\p)}(\Phi_0)_*.
\end{equation}
\end{tm}

\begin{tm}\label{Another FCLT}
The sequence $(\mathcal{Y}_t^{(n;\, \p)})_{0 \leq t \leq 1}$ converges in law 
to the $G$-valued diffusion process $(\widetilde{Y}_t)_{0 \leq t \leq 1}$ 
in $C^{0, \alpha\text{\normalfont{-H\"ol}}}_{\bm{1}_G}([0, 1]; G)$ as $n \to \infty$ 
for all $\alpha<1/2$.
\end{tm}
We emphasize that 
the transition probability $\p$
coincides with the given one $p$ under the centered condition {\bf (C)}. 
Therefore, Theorems \ref{Another CLT} and \ref{Another FCLT} are regarded 
as extensions of Theorems \ref{CLT1} 
(under the centered condition {\bf (C)}) and \ref{FCLT1} 
to the non-centered case.
We might prove Theorem \ref{FCLT1} without the centered condition {\bf (C)}
via Theorem \ref{Another FCLT}. 
We will discuss this problem in the future. 

\end{appendix}

\vspace{4mm}
\noindent
{\bf Acknowledgement.} 
The authors are grateful to Professor Shoichi Fujimori for making pictures 
of the 3-dimensional Heisenberg dice lattice 
and kindly allowing them to use these pictures in the present paper. 
They would also like to thank Professors Takafumi Amaba, Takahiro Aoyama, Peter K. Friz, 
Naotaka Kajino,
Atsushi Katsuda, Takashi Kumagai, Seiichiro Kusuoka, Kazumasa Kuwada, Laurent Saloff-Coste and Ryokichi Tanaka 
for helpful discussions and encouragement.
A part of this work was done during the stay of the third named author at 
Hausdorff Center for Mathematics, Universit\"at Bonn
in March 2017 with the support of research fund of 
Research Institute for Interdisciplinary Science, Okayama University. 
He would like to thank Professor Massimiliano Gubinelli for 
warm hospitality and helpful discussions.




\end{document}